\documentclass[final]{dmtcs-episciences}
\publicationdata{vol. 25:2 }{2023}{3}{10.46298/dmtcs.9293}{2022-04-03;
  2022-04-03; 2023-03-23; 2023-06-28}{2023-06-29}
\usepackage[utf8]{inputenc}
\usepackage{subfigure}
\usepackage[round]{natbib}

\usepackage{epsfig}
\usepackage{latexsym}
\usepackage{amsfonts,amsmath,amssymb,textcomp,dsfont}
\usepackage{stackengine}
\usepackage{algorithm}
\usepackage{algorithmic}
\usepackage{xcolor}
\usepackage{hyperref}
\usepackage[title]{appendix}
\allowdisplaybreaks

\newtheorem{thm}{Theorem}[section]
\newtheorem{lem}[thm]{Lemma}
\newtheorem{cor}[thm]{Corollary}
\newtheorem{prop}[thm]{Proposition}
\newtheorem{rem}[thm]{Remark}
\newtheorem{deff}[thm]{Definition}
\newtheorem{conj}[thm]{Conjecture}
\newtheorem{note}[thm]{Note}
\newtheorem{key}[thm]{Keywords}
\newtheorem{prob}[thm]{Problem}
\newenvironment{probb}{ \begin{prob} \rm}{ \end{prob} }
\newcommand{\bprobb}{\begin{probb}}
\newcommand{\eprobb}{\end{probb}}
\newcommand{\bth}{\begin{thm}}
\newcommand{\ethGL}{\end{thm}}
\newcommand{\bconj}{\begin{conj}}
\newcommand{\econj}{\end{conj}}
\newcommand{\bkey}{\begin{key}}
\newcommand{\ekey}{\end{key}}
\newcommand{\bl}{\begin{lem}}
\newcommand{\el}{\end{lem}}
\newcommand{\bdeff}{\begin{deff}}
\newcommand{\edeff}{\end{deff}}
\newcommand{\bcor}{\begin{cor}}
\newcommand{\ecor}{\end{cor}}
\newcommand{\bprop}{\begin{prop}}
\newcommand{\eprop}{\end{prop}}
\newcommand{\brem}{\begin{rem}}
\newcommand{\erem}{\end{rem}}
\newcommand{\beq}{\begin{equation}}
\newcommand{\eeq}{\end{equation}}
\newcommand{\beqn}{\begin{eqnarray}}
\newcommand{\eeqn}{\end{eqnarray}}
\newcommand{\beqns}{\begin{eqnarray}}
\newcommand{\eeqns}{\end{eqnarray}}
\newcommand{\ba}{\begin{array}}
\newcommand{\ea}{\end{array}}
\newcommand{\bit}{\begin{itemize}}
\newcommand{\eit}{\end{itemize}}
\newcommand{\bpr}{\noindent{\bf Proof.\hspace{1em}}}
 \newcommand{\epr}{\hfill\rule{3mm}{3mm}\vspace{\baselineskip}\\}
 \newcommand{\BO}{\mathcal{O}}

\newcommand{\BMC}{\mathcal{M}}
 \newcommand{\non}{\nonumber}
\newcommand{\bal}{\begin{align}}
\newcommand{\bals}{\begin{align*}}
\newcommand{\bs}{\begin{skip}}
 \newcommand{\eal}{ \end{align*}}
\newcommand{\eals}{\end{align*}}
\newcommand{\es}{\end{skip}}

 \newcommand{\ra}{\rightarrow}

\newcommand{\lf}{\left \lfloor}
\newcommand{\rf}{\right \rfloor}

\newcommand{\kl}{[\hspace{-0.5 mm}[}
 \newcommand{\kr}{]\hspace{-0.5 mm}]}

\newcommand{\lp}{\left (}
 \newcommand{\rp}{\right )}
\newcommand{\lb}{\left [}
 \newcommand{\rb}{\right ]}

  \def\ii{\mathbf{i}}

 \newcommand{\bx}{\mathbf{x}}
 \newcommand{\bp}{\mathbf{p}}
\newcommand{\bw}{\mathbf{w}}
 \newcommand{\bu}{\mathbf{u}}
\newcommand{\bv}{\mathbf{v}}
\newcommand{\bbe}{\mathbf{e}}
\newcommand{\R}{\mathbb{R}}

\newcommand{\ns}{n^*}
\newcommand{\is}{i^*}

\newcommand{\bet}{\tilde{\be}}

\newcommand{\ut}{\tilde{u}}

\newcommand{\Pt}{\tilde{P}}

\newcommand{\Xt}{\tilde{X}}

\newcommand{\Pb}{\bar{P}}

\newcommand{\Kb}{\bar{\K}}
\newcommand{\mub}{\bar{\mu}}

\newcommand{\K}{\kappa}
\newcommand{\la}{\lambda}

\newcommand{\al}{\alpha}
\newcommand{\be}{\beta}

\newcommand{\eps}{\varepsilon}
\newcommand{\II}{\infty}
\newcommand{\tet}{\theta}

\newcommand{\sig}{\sigma}
\newcommand{\sigd}{\sigma^2}

\newcommand{\gam}{\gamma}
\newcommand{\Gam}{\Gamma}

\def\1{{\ifmmode 1\mskip-1.5\thinmuskip\mathrm{l}%
        \else\textrm{1\hskip -.23em l}\fi}}

\newcommand{\Z}{\mbox{$\mathbb Z$}}
\newcommand{\E}{\mbox{$\mathbb E$}}

\def\P{{\mathbb {P}}}
\def\Pr{{\mathbb {P}}}

\definecolor{dgreen}{rgb}{0,.7,0}

\newcommand{\ignore}[1]{}
\newcommand{\Var}{\textup{Var}\,}
\newcommand{\Cov}{\textup{Cov}\,}
\newcommand{\Id}{\mathbb{I}}
\stackMath
\newcommand\tsup[2][2]{%
 \def\useanchorwidth{T}%
  \ifnum#1>1%
    \stackon[-.5pt]{\tsup[\numexpr#1-1\relax]{#2}}{\scriptscriptstyle\sim}%
  \else%
    \stackon[.5pt]{#2}{\scriptscriptstyle\sim}%
  \fi%
}

\newcommand{\bin}[2]
{
{#1\choose #2}
}

\def\la{\lambda}
\def\bu{\mathbf{u}}
\usepackage{url}

\author[Louchard et. al]{Guy Louchard\affiliationmark{1}
  \and Werner Schachinger\affiliationmark{2}
  \and Mark Daniel Ward\affiliationmark{3}}
\title[The number of distinct adjacent pairs in geometrically distributed words]{The number of distinct adjacent pairs in geometrically distributed words: a probabilistic and combinatorial analysis}
\affiliation{
  Universit\'e Libre de Bruxelles, Belgium\\
  University of Vienna, Austria\\
  Purdue University, USA}

\begin{document}

\maketitle
\begin{abstract}
The analysis of strings of $n$ random variables with geometric
distribution has recently attracted renewed interest: Archibald et~al.\ consider the number of distinct adjacent pairs in geometrically distributed words. They obtain the asymptotic ($n\ra\infty$) mean of
this number in the cases of different and identical  pairs. In this
paper we are interested in all asymptotic moments  in the identical
case, in the asymptotic variance in the different case and in  the
asymptotic distribution in both cases. We use two approaches: the
first one, the probabilistic approach, leads to variances in both
cases and to some conjectures on all  moments in the identical case
and  on the distribution in both cases. The second approach, the
combinatorial  one, relies on multivariate pattern matching
techniques, yielding exact formulas for first and second moments. We use such tools as Mellin transforms, Analytic Combinatorics, Markov Chains.
\end{abstract}

\section{Introduction}
We follow the notation and setup of~\cite{ABCKWW21}.
In this earlier work, the authors derived results about the asymptotic
mean of the numbers of different and identical pairs, in a sequence of
geometric random variables.  \cite{ABCKWW21} give a
broad selection of references to the literature, including
applications to leader election algorithms, pattern matching in
randomly generated words and permutations, gaps in sequences, the design
of codes, etc.  In the present work, we go far beyond the analysis of
the mean numbers of different and identical pairs.  We use two
approaches, namely, a probabilistic approach and also a combinatorial
approach.  We are able to derive results about the asymptotic
variance and distribution, and to make conjectures about higher
moments.  We also derive exact results, using multivariate pattern
matching, for the first and second moments.

As motivated by \cite{ABCKWW21},
we consider a string of $n$ independent random variables 
${Z_1,Z_2,\ldots,Z_n}$, with geometric distribution
$\mathbb{P}(Z_k=i)=P_i:=p\,q^{i-1}$ for $i\ge1$. 
Our eventual aim is to study the consecutive pairs of geometric random
variables in this sequence, with a goal of characterizing the
asymptotic behavior, as $n \ra \II$.

We use Iverson's notation, namely, for an event $A$, we write $\kl A
\kr = 1$ if event $A$ occurs, and $\kl A \kr = 0$ otherwise.
We want to precisely characterize the 
distribution of the number of times that $(i,j)$ appears as a
consecutive pair in ${Z_1,Z_2,\ldots,Z_n}$, i.e., the number
of $k$'s such that $X_{k}=i$ and $X_{k+1}=j$.  So we define
$X^{(n)}_{i,j}(m)$ as a Bernoulli random variable that indicates
whether the pair $(i,j)$ appears $m$ times in a sequence of $n$
geometric random variables:
$$
X^{(n)}_{i,j}(m):= \kl \mbox{pair } (i,j) \mbox{ appears } m \mbox{
  times in the string of size } n\kr.$$
It is useful to have a succinct notation for the Bernoulli random
variable $X_{i,j}^{(n)}$ that indicates that $(i,j)$ appears at least
one time in a sequence of $n$ geometric random variables:
$$
X_{i,j}^{(n)} :=1-X_{i,j}^{(n)}(0):= \kl \mbox{pair } (i,j) \mbox{
  appears at least once  in the string of size } n\kr.$$
Finally, we define $X_1^{(n)}$ as the number of types of matching
consecutive pairs (we say ``types'' because we only pay attention to
whether a pair $(i,i)$ occurs or does not occur, i.e., whether it
never occurs, or whether it occurs one or more times):
$$
X_1^{(n)}:=\sum_{i\ge1}  X_{i,i}^{(n)}.$$
Similarly, $X_2^{(n)}$ is the number of types of any matching
consecutive pairs (different or matching):
$$
X_2^{(n)}:=\sum_{i,j\ge1}  X_{i,j}^{(n)},$$
and finally $X_3^{(n)}$ is the number of types of different
consecutive pairs that occur:
$$
X_3^{(n)}:=\sum_{i\neq j}  X_{i,j}^{(n)}.
$$
Our methodology is to derive asymptotic expressions for the moments,
utilizing Mellin transforms applied to harmonic sums.  For context and
an in-depth explanation of such techniques, see the nice exposition in
\cite{FGD95}. 

One highlight of the precision of this analytic method is that we are
able to derive the dominant part of moments as well as the (tiny) periodic part, in the form of a Fourier series.
\ignore{ Rapid decrease  property: $\Gam(s)$ decreases exponentially in the direction  $\ii\II$:
\[|\Gam( \sig+\ii t )| \sim \sqrt{2\pi}|t|^{\sig-1/2}e^{-\pi |t|/2}.\]
Also, this property is true for all other functions we encounter. So inverting the Mellin transforms is easily justified  (see \cite[Sec.4]{LP04} for example).
All our expressions are given with an error term $n^{-\be}$, with  $\be > \Re(\chi)$, where $\chi$ is the rightmost pole used in our
singular expansions.}

The paper is organized as follows: In Section~2 we present our main
results, that is, asymptotic expressions for the variances of
$X^{(n)}_k,\,1\le k\le3$, and a result concerning the asymptotic independence of
the variables $X^{(n)}_{i,i},\,i\in\mathbb{N}$. In Section~3 we conjecture some stronger forms of asymptotic independence, based on which we are able to derive the limiting distribution and asymptotics of higher moments of $X^{(n)}_1$. Section~4 is devoted to the proofs of these results,
and to some considerations in support of a conjectured Gaussian limiting distribution of $X^{(n)}_3$.
In Section~5 we use a combinatorial approach to derive exact expressions
for first and second moments of $X^{(n)}_k,\,1\le k\le3$.
In the Appendix, we collect our results pertaining to Mellin transforms.
\section{Main results}

\begin{figure}[htbp]
	\centering
		\includegraphics[width=0.55\textwidth,angle=0]{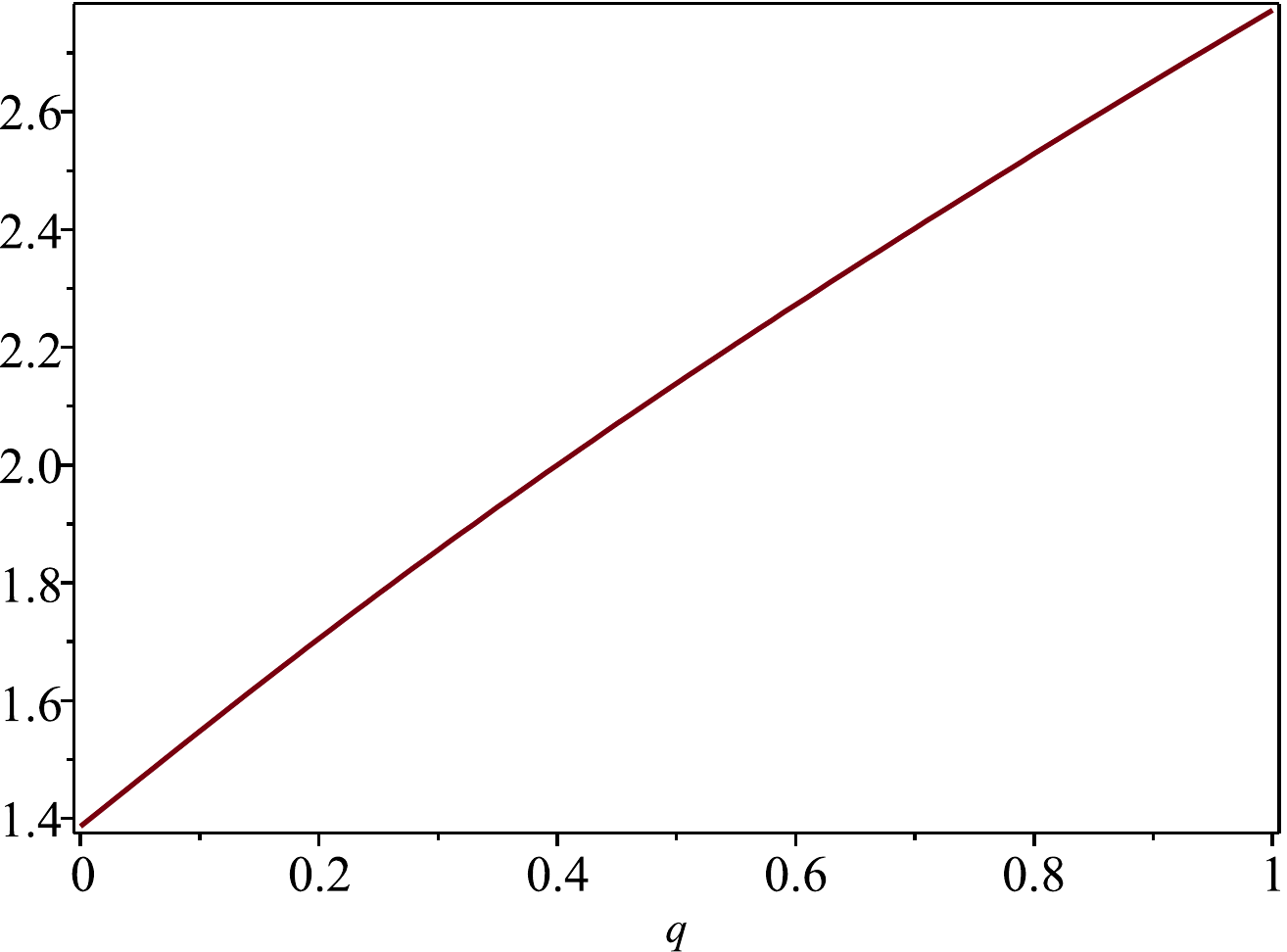}
	\caption{ Plot of $2(1-q)F'_1(0)$, showing the dependence of the constant term
	$\frac{2}{L}F'_1(0)$ on $q$.
	We leave it as an exercise to show that, for $q\to0$
        (resp.~$q\to1$), the limit is $2\ln 2$ (resp.~$4\ln 2$).}
	\label{F2}
\end{figure}

In a private communication, B.~Pittel observed that
the asymptotic distribution of $X^{(n)}_{i,j}(m)$ is Poisson,
\[\P[X^{(n)}_{i,j}(m)=1]\sim e^{-\la}\frac{\la^m}{m!},~\mbox{where}~\la=nP_i P_j.\]
Asymptotics of $\E X^{(n)}_1$, $\E X^{(n)}_2$ and $\E X^{(n)}_3$ have
also recently been obtained by \cite{ABCKWW21}, using generating
functions of the sequences of expectations.
One of our main results deals with asymptotics of Var$\,X^{(n)}_i$, $1\le i\le 3$, as $n\to\infty$. Our approach simply consists in using
$$\Var X^{(n)}_1=\sum_{i\ge1}\Var X_{i,i}^{(n)}+\sum_{i\ne j}\Cov(X_{i,i}^{(n)},X_{j,j}^{(n)}),$$
and similarly for $\Var X^{(n)}_2$ and $\Var X^{(n)}_3$. This necessitates thorough investigation of the involved covariances.
As it turns out, the main term of $\Var X^{(n)}_1$ is given by a term $S^{(n)}_1\sim\sum_{i\ge1}\Var X_{i,i}^{(n)}$, the double sum of covariances only contributing $\BO(\frac1n)$. This is different for $\Var X^{(n)}_2$, whose main term is a sum of
$S^{(n)}_2\sim\sum_{i,j\ge1}\Var X_{i,j}^{(n)}$ and another contribution $T^{(n)}_2$, stemming from the quadruple sum of covariances of different pairs, of order $\Theta(1)$. All of $S^{(n)}_1$, $S^{(n)}_2$, and $T^{(n)}_2$ are expressed in terms of Fourier series in $\ln (np^2)$.
A plot of the constant term of $T^{(n)}_2$ is provided in Figure~\ref{F2}.
\bth \label{thm:var123} Let $L:=\ln(1/q)$ and $\chi:=2\pi\ii/L$, where $\ii$ denotes the imaginary unit. We
also define
\begin{align}
S^{(n)}_1&:=\frac{\ln2}{2L}+\frac1{2L}\sum_{\ell\ne0}\Gamma\lp\frac{\ell\chi}2\rp
(np^2)^{-\frac{\ell\chi}2}\lp1-2^{-\frac{\ell\chi}2}\rp,\label{S1}\\
S^{(n)}_2&:=\frac{\ln 2}{L^2}\ln (np^2)+\frac{\ln 2}{2L^2}(2\gamma+\ln 2+2L)\notag\\
&~~~~+\frac{\ln (np^2)}{L^2}\sum_{\ell\ne0}\Gamma(\ell\chi)(np^2)^{-\ell\chi}
(1-2^{-\ell\chi})\label{S2}\\
&~~~~-\frac{1}{L^2}\sum_{\ell\ne0}\Gamma(\ell\chi)(np^2)^{-\ell\chi}
\left[
(1-2^{-\ell\chi})\left(\frac{\Gamma'(\ell\chi)}{\Gamma(\ell\chi)}-L\right)+2^{-\ell\chi}\ln2
\right]\notag\\
T^{(n)}_2&:=\frac{2}{L}F'_1(0)+\frac2L\sum_{\ell\ne0}
  \Gam(\ell\chi) F_1(\ell\chi) (np^2)^{-\ell\chi},
\label{T2}\end{align}
where $F_1(s)=\sum_{i,k\ge1}\lb(q^i+q^k-pq^{i+k-1})^{-s}-(q^i+q^k)^{-s}\rb$,
and the constant term of $T^{(n)}_2$ simplifies to
$$\frac{2}{L}F'_1(0)=-\frac2L\ln\Bigg(\prod_{i,j\ge1}\lp1-\frac pq\frac{q^{i+j}}{q^i+q^{j}}\rp\Bigg).$$
Then, as $n\to\infty$, the variances of $X^{(n)}_i$, $1\le i\le 3$, satisfy
\begin{align}\Var X^{(n)}_1&=S^{(n)}_1+\BO\lp\frac1{\sqrt{n}}\rp,\label{var1}\\
\Var X^{(n)}_2&=S^{(n)}_2+T^{(n)}_2+\BO\lp\frac{\ln n}{\sqrt{n}}\rp,\label{var2}\\
\Var X^{(n)}_3&=S^{(n)}_2-S^{(n)}_1+T^{(n)}_2+\BO\lp\frac{\ln n}{\sqrt{n}}\rp.\label{var3}
\end{align}
\ethGL
A question triggered by the observation that $\sum_{i\ne j}\Cov(X_{i,i}^{(n)},X_{j,j}^{(n)})=\BO\lp\frac1n\rp$ is: How ``close to being independent''
are $\big( X_{i,i}^{(n)}\big)_{i\in\mathbb{N}}$ ? The following theorem provides a partial answer in that regard.

\bth \label{thm:asyindep}
The random variables $X_{i,i}^{(n)},i\in\mathbb{N}$ are asymptotically independent, in the sense that, for any $k\in\mathbb{N}$, any subset $I\subseteq\mathbb{N}$ of size $k$, and any $(x_i)_{i\in I}\in\{0,1\}^k$
we have
\beq\label{eq:asyindep}
\P(X_{i,i}^{(n)}=x_i,i\in I)-\prod_{i\in I}\P(X_{i,i}^{(n)}=x_i)=\BO\left(\frac1n\right),
\eeq
with implied constant depending on $I$ only via $k$.
\ethGL

\brem
The random variables $\big(X_{i,i}^{(n)}\big)_{i\ge1}$ are negatively correlated: For finite $I\subseteq\mathbb{N}$ we have
$$\Pr(X_{i,i}^{(n)}=1,i\in I)\le\prod_{i\in I}\Pr(X_{i,i}^{(n)}=1),$$
as can easily be deduced from the following theorem.\\
{\bf Theorem (\cite{MD92}):} {\rm Let $V$ and $I$ be finite non-empty sets. Let $(Z_v:v\in V)$ be a family of independent random variables, each taking values in some set containing $I$; and for each $i\in I$, let $S_i=\{v\in V: Z_v=i\}$. Let $(\mathcal{F}_i:i\in I)$ be a family of collections of subsets of $V$ such that each collection is increasing (meaning that every superset of a set in $\mathcal{F}_i$ is also in $\mathcal{F}_i$) or each is decreasing
(meaning that every subset of a set in $\mathcal{F}_i$ is also in $\mathcal{F}_i$).
Then $\Pr\left(\bigcap_{i\in I}\{S_i\in \mathcal{F}_i\}\right)\le
\prod_{i\in I}\Pr\left(\{S_i\in \mathcal{F}_i\}\right)$.}\\
We just have to choose $V:=\{1,\ldots,n\}$, and all $\mathcal{F}_i$ equal to
$\mathcal{F}:=\{f\subseteq V:\exists k\in V:\{k,k+1\}\subseteq f\}$.

Cases like the following for $n=5$ and $i\ne j$,
{\small\begin{align*}\Pr(&X_{i,i}^{(5)}=1)\Pr(X_{j,j}^{(5)}=1)-\Pr(X_{i,i}^{(5)}=X_{j,j}^{(5)}=1)\\
&=\!P_i^2(4-2P_i-2P_i^2+P_i^3)P_j^2(4-2P_j-2P_j^2+P_j^3)-
P_i^2P_j^2(6-2P_i-2P_j)\\
&=\!P_i^2P_j^2\lb(1\!-\!P_i\!-\!P_j)(10\!+\!4P_i(1\!-\!P_i)\!+\!4P_j(1\!-\!P_j))+P_iP_j(12\!+\!P_i(2\!-\!P_i)P_j(2\!-\!P_j)\!-\!2P_i^2\!-\!2P_j^2)\rb\!>\!0
\end{align*}}\noindent
suggest
that the inequality may be strict for $|I|\ge2$. This is different for the
array $\big(X_{i,j}^{(n)}\big)_{i,j\ge1}$, where both strictly positive and strictly negative correlations can be observed: For $n=3$ and $i\ne j$,
{\small$$\Pr(X_{i,j}^{(3)}=X_{j,i}^{(3)}=1)-\Pr(X_{i,j}^{(3)}=1)\Pr(X_{j,i}^{(3)}=1)
=P_iP_j(P_i+P_j)-(2P_iP_j)^2=P_iP_j(P_i+P_j-4P_iP_j)>0$$
}\noindent
holds for $P_i,P_j$ small enough, and for different pairs $\big((k_i,m_i)\big)_{i\in I}$, with $|I|\ge n$, we clearly have
$$\Pr(X_{k_i,m_i}^{(n)}=1,i\in I)=0<\prod_{i\in I}\Pr(X_{k_i,m_i}^{(n)}=1).$$
\erem
\section{Further conjectures and results for pairs of identical letters}
\subsection{Higher moments} \label{S31}
The proof of Theorem~\ref{thm:var123} (see
Lemma~\ref{varsumeqsumvar}) shows that
$$\lim_{n\to\infty}(\textup{Var}\,X^{(n)}_1 - \textup{Var}\,\xi^{(n)})=0,$$
where $\xi^{(n)}:=\sum_{i\ge1}\kl\xi_i^{(n)}\ge1\kr$ is a sum of independent random variables,
with $\xi_i^{(n)}$ distributed as Poisson$(nP_i^2)$. Note that $\P[X_{i,i}^{(n)}=0]\sim \P[\xi_i^{(n)}=0]$ and
$\P[X_{i,i}^{(n)}=1]\sim \P[\xi_i^{(n)}\ge1]$. This leads us to the following conjecture.

\begin{conj}\label{conj:moments} For any $k\in \mathbb{N}$ we have
$\lim_{n\to\infty}(\E|X^{(n)}_1- \E X^{(n)}_1|^k - \E|\xi^{(n)}-\E\xi^{(n)}|^k)=0$.\end{conj}

\begin{thm}\label{thm:moments} If Conjecture~\ref{conj:moments} holds,
  the asymptotics of cumulants $\kappa_m^{(n)}$ of $X^{(n)}_1$   are given by
\beq\kappa_m^{(n)}=m!\sum_{j=1}^mV^{(n)}_j\frac{(-1)^{j+1}}{j}[\tet^m](e^\tet -1)^j,\label{E40}\eeq
where, using $L=\ln(1/q)$ again, asymptotics of $V^{(n)}_j,\,j\ge1$, are given by
\beq\begin{aligned}\label{E41}
V^{(n)}_j\sim\frac{\ln(np^2)}{2L}&+\frac{\gam}{2L}+\frac12+\frac1{2L}\sum_{k=2}^j  (-1)^{k+1} \bin{j}{k}\ln k\\
&+\frac1{2L}\sum_{\ell\ne0}\lp\sum_{k=1}^j  (-1)^{k} \bin{j}{k}k^{-\frac{\ell\chi}2}\rp\Gamma\lp\frac{\ell\chi}2\rp
(np^2)^{-\frac{\ell\chi}2}.
\end{aligned}\eeq
\end{thm}
\bpr We proceed as in \cite{HL00} and \cite{LP04}.

Let $S_{n}(\tet):=\ln(E ( e^{\tet \xi^{(n)}}))=\sum_{m=1}^\II\kappa_m^{(n)}\frac{\tet^m}{m!}$ be the cumulant generating function of $\xi^{(n)}$. Furthermore let
$n_2:=np^2/q^2$, and observe
$\E e^{\tet\kl\xi_i^{(n)}\ge1\kr}=1+(e^\tet-1)(1-e^{-n_2q^{2i}})$. By independence of $\big(\xi_i^{(n)}\big)_{i\ge1}$, we get
\begin{align*}
 S_{n}(\tet)&=\sum_{i=1}^\II \ln \lb  1+\lp e^\tet-1\rp \lp 1-e^{-n_2q^{2i}}\rp \rb\\
 &=\sum_{j=1}^\II \frac{(-1)^{j+1}}{j}(e^\tet -1)^j \lb \sum_{i=1}^\II \lp 1-e^{-n_2 q^{2i}}\rp^j  \rb.\end{align*}
 Now let
\begin{align*}V^{(n)}_j&:=\sum_{i=1}^\II \lp 1-e^{-n_2 q^{2i}}\rp^j
=\sum_{i=1}^\II \lb  \sum_{k=0}^j  (-1)^k \bin{j}{k} e^{-kn_2 q^{2i}} \rb\\
 &=\sum_{i=1}^\II \lb  \sum_{k=0}^j  (-1)^k  \bin{j}{k} e^{-kn_2 q^{2i}} - \sum_{k=0}^j  (-1)^k  \bin{j}{k}\rb
=\sum_{i=1}^\II \lb  \sum_{k=1}^j  (-1)^{k+1} \bin{j}{k}\lp 1- e^{-kn_2 q^{2i}}\rp \rb\\
&=\sum_{k=1}^j  (-1)^{k+1} \bin{j}{k}\sum_{i=1}^\II\lp 1- e^{-kn_2 q^{2i}}\rp,
\end{align*}
where the asymptotics of the inner sum can be obtained using
$G(knp^2)$ from Appendix~\ref{A1}, leading to \eqref{E41}.
Finally the cumulants $\kappa_m^{(n)}$ are found by extracting coefficients of $\tet^m$ from $S_n(\tet)$, and are given by finite linear combinations of the $(V^{(n)}_j)_{j\ge1}$, as stated in \eqref{E40}.\epr
\begin{rem}\label{rem:OEIS}
Explicit expressions for \eqref{E40} for small $m$ are
\beq\notag
\kappa_1^{(n)}\!\!=V^{(n)}_1,\quad \kappa_2^{(n)}\!=\!V^{(n)}_1\!-\!V^{(n)}_2,\quad \kappa_3^{(n)}\!=\!V^{(n)}_1\!-\!3V^{(n)}_2\!+\!2V^{(n)}_3,\quad\kappa_4^{(n)}\!=\!V^{(n)}_1\!-\!7V^{(n)}_2\!+\!12V^{(n)}_3\!-\!6V^{(n)}_4.
\eeq
The fact that $\frac1{j!}(e^x -1)^j$ is the generating function for the Stirling numbers of the second kind, see e.g.\ \cite[p.\,736]{FS09}, establishes that
the sequence of (absolute values of) the coefficients,
$(1,1,1,1,3,2,1,7,12,6,\ldots)$, is equal to OEIS sequence A028246 in \cite{OEIS}.
\end{rem}

The cumulants now allow for computation of moments: The mean of   $X^{(n)}_1$  is given by
\[\E X^{(n)}_1\sim V^{(n)}_1.\]
This is identical to \cite[Thm.\,2]{ABCKWW21}, see also \eqref{EX1}. Our approach here is simple and general.
Note that the mean does not rely on the state of
Conjecture~\ref{conj:moments}: the mean computation actually depends only on Lemma~\ref{Piiir}.
Similarly, the variance of $X^{(n)}_1$  is given by
\[\Var X^{(n)}_1\sim V^{(n)}_1-V^{(n)}_2.\]
After some algebra, we verify that this is identical to Thm \ref{thm:var123} .

\subsection{Limiting distribution}
A conjecture weaker than Conjecture~\ref{conj:moments} is

\begin{conj}\label{conj:limdist} For any $t\in \mathbb{R}$ we have
$\lim_{n\to\infty}[\mathbb{P}(X^{(n)}_1\le t) - \mathbb{P}(\xi^{(n)}\le t)]=0$.\end{conj}

\begin{figure}[htbp]
	\centering
		\includegraphics[width=0.72\textwidth,angle=0]{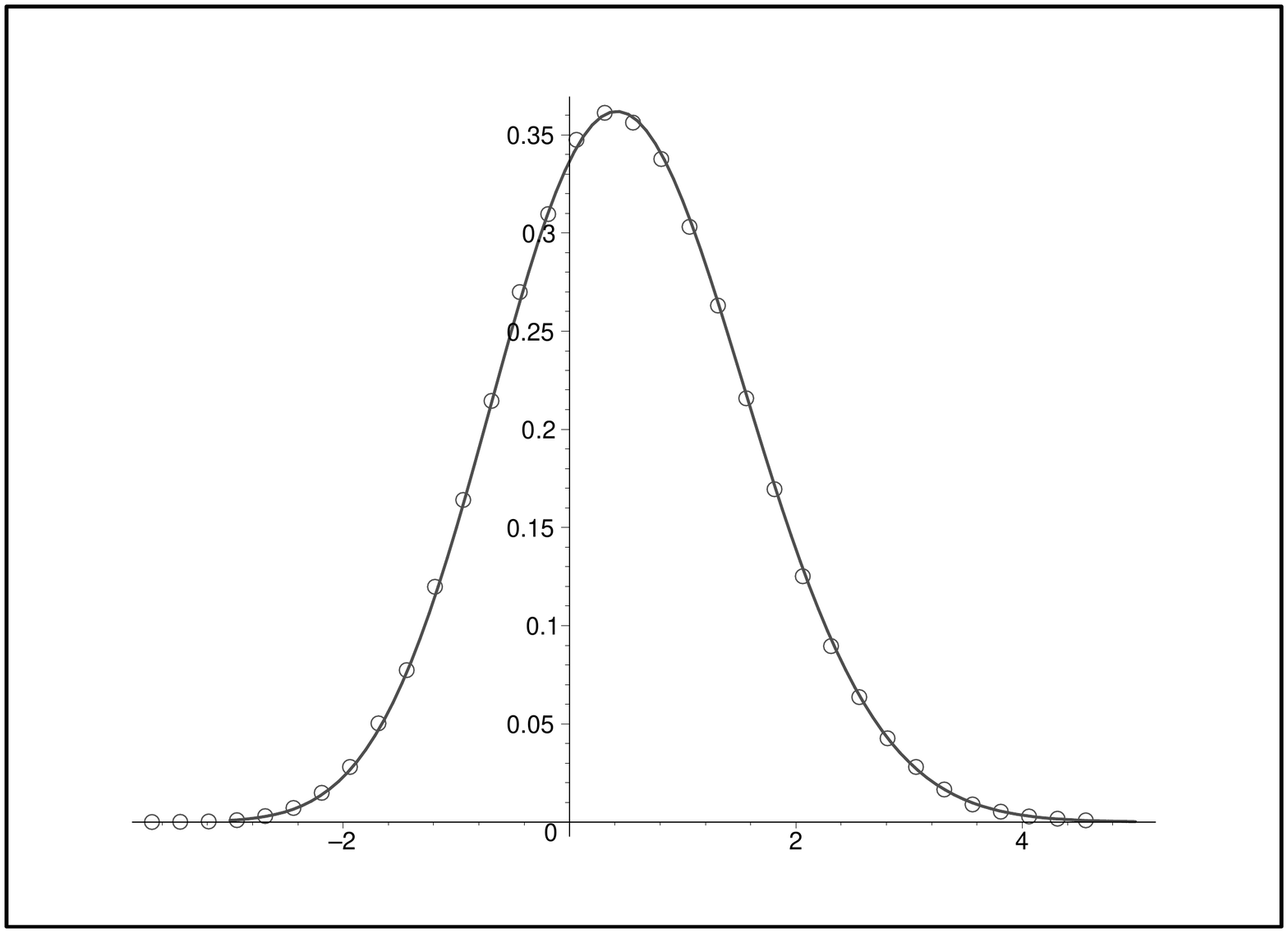}
	\caption{ Comparison between $f(\eta)$ (line) and the simulation
of $X^{(n)}_1$ (circles), $p=1/4$, number of simulated words $=50000$ for each $n\in\{10000,11547,13333,15396\}$.}
	\label{F1}
\end{figure}

\begin{thm}\label{thm:asydist} If Conjecture~\ref{conj:limdist} holds, the asymptotic distribution $f(\eta)$ of $X^{(n)}_1$ is given by (\ref{E42}).
\end{thm}
Set again $L=\ln(1/q)$ and $n_2=np^2/q^2$, set $\is=\ln(n_2)/(2L)$ (implying $q^{2\is}=1/n_2$), define $\eta:=i-\is$,
and use $\P(\xi_i^{(n)}=0)=e^{-n_2 q^{2i}}=\exp(-e^{-2L\eta})$.
This leads to
$$\P(\xi_k^{(n)}=0,k>i) = \exp\lp-\al e^{-2L\eta}\rp,
\mbox{ where }\al:=\frac{q^2}{1-q^2}.$$
As in \cite{HL00} and \cite{ALP05}, we proceed by defining
$$\Psi(\eta):=e^{-e^{-2L \eta}} \prod_{i=1} ^\II \lb 1-e^{-e^{-2L(\eta-i)}}\rb,$$
and observing that, as $n\to\II$, we have
\bal
\P(X^{(n)}_1=i^*+\eta)\sim &\, f(\eta):=\sum_{v=0}^\II \Psi(\eta-v+1)e^{-\al e^{-2L(\eta+1-v)}}\sum_{\overset{r_1<\cdots<r_v}{r_j\ \geq 2-v}}
\prod_{i=1}^v \frac{1-e^{-e^{-2L(\eta+r_i)}} }{e^{-e^{-2L(\eta+r_i)}} }\label{E42},\\[-.2cm]
\P(X^{(n)}_1\leq i^*+\eta)\sim &\, F(\eta):=\sum_{i=0}^\II f(\eta-i)\non.\\
f(\eta)& \mbox{ depends only on } p\non.
\end{align}
A simulation with $p=1/4$ and $50000$ simulated words for each $n\in\{10000,11547,13333,15396\}$ is given in Figure~\ref{F1}. The fit is excellent.
A corresponding table of observed and theoretical non-periodic mean
and variance in the equal pairs case (as well as another table for the unequal pairs case) is given below, all results rounded to 3 decimal places. We define $\bar X^{(n)}_j:=\frac1N\sum_{i=1}^NX_j^{(n),i}$
and $s_j^2(n):=\frac1{N-1}\sum_{i=1}^N(X_j^{(n),i}-\bar X^{(n)}_j)^2$
the sample mean and unbiased sample variance of a sample $(X_j^{(n),i})_{i=1}^N$.
See Theorem~\ref{EX1EX3} for asymptotics of
$\E X^{(n)}_1$ and $\E X^{(n)}_3$. Both simulations use $p=1/4$. The sample size $N$ for each row in the left table is $50000$, and in the right table it is $200000$, see also Figure~\ref{F3}.
\renewcommand{\arraystretch}{1.1}
$$\begin{tabular}{|r|c|c|c|c|}
\hline
$n$&$\E X^{(n)}_1$&$\bar X^{(n)}_1$&$\Var X^{(n)}_1$&$s_1^2(n)$\\
\hline\hline
$10000$&$12.692$&
$12.676$&$1.205$&$1.214$\\
\hline
$11547$&$12.942$&
$12.927$&$1.205$&$1.206$\\
\hline
$13333$&$13.192$&
$13.175$&$1.205$&$1.213$\\
\hline
$15396$&$13.442$&
$13.427$&$1.205$&$1.211$\\
\hline
\end{tabular}\quad
\begin{tabular}{|r|c|c|c|c|}
\hline
$n$&$\E X^{(n)}_3$&$\bar X^{(n)}_3$&$\Var X^{(n)}_3$&$s_3^2(n)$\\
\hline\hline
$500000$&$750.195$&
$750.198$&$129.889$&$130.053$\\
\hline
\end{tabular}$$
\begin{rem} Here we briefly sketch, how we obtained the graph of $f$ in Figure~\ref{F1}, where $p=1/4$. As before, we use random variables $\xi_i^{(n)}$ distributed Poisson$(np^2q^{2(i-1)})$, but now there is such a random variable for each $i\in\mathbb{Z}$ and each real $n>0$. For fixed such $n$ the
random variables $\big(\xi^{(n)}_i\big)_{i\in\mathbb{Z}}$ are assumed independent, and also the definition
$\xi^{(n)}:=\sum_{i\ge1}\xi^{(n)}_i$ is used for real $n>0$. We use
$i^*=i^*(n)=\ln(np^2/q^2)/(2L)$ again. For any $n$ satisfying $i^*+\eta\in\mathbb{Z}$, we have
\begin{align*}f(\eta)&=\lim_{k\to\infty}\P\Big(\xi^{(nq^{-2k})}-k=i^*+\eta\Big)
=\P\Big(\sum_{i\ge1}\kl\xi_i^{(n)}\ge1\kr+\sum_{j\ge0}(\kl\xi_{-j}^{(n)}\ge1\kr-1)=i^*+\eta\Big)\\
&=\P\Big(\sum_{i\ge1}\kl\xi_i^{(\nu)}\ge1\kr+\sum_{j\ge0}(\kl\xi_{-j}^{(\nu)}\ge1\kr-1)=0\Big),
\end{align*}\goodbreak
where for $n=\nu=\nu(\eta):=q^{2(1-\eta)}/p^2$ we have $i^*+\eta=0$,
and $\xi^{(\nu)}_i\sim$\,Poisson$\big(q^{2(i-\eta)}\big)$.
We want a good approximation of $f(\eta)$ only for $\eta\in[-3,5]$.
For such $\eta$ we have
$$\P\Big(\sum_{i>30}\kl\xi_i^{(\nu)}\ge1\kr>0\Big)=1-\prod_{i>30}e^{-q^{2i-2\eta}}\le1-\prod_{i>30}e^{-q^{2i-10}}
=1-\exp\left(-\frac{q^{52}}{1-q^2}\right)\approx7.28\cdot10^{-7}.
$$
and
$$\P\Big(\sum_{j>7}(\kl\xi_{-j}^{(\nu)}\ge1\kr-1)<0\Big)=
1-\prod_{j>7}\Big(1-e^{-q^{-2j-2\eta}}\Big)\le\sum_{j>7}e^{-q^{6-2j}}
\approx e^{-q^{-10}}\approx1.94\cdot10^{-8}.$$
So, up to an error smaller than $10^{-6}$, $f(\eta)$ is given by
\begin{align*}\P\Big(\sum_{i=1}^{30}\kl\xi_i^{(\nu)}\ge1\kr+\sum_{j=0}^7(\kl\xi_{-j}^{(\nu)}\ge1\kr-1)=0\Big)&=\P\Big(\sum_{i=-7}^{30}\kl\xi_i^{(\nu)}\ge1\kr\!=\!8\Big)\\
&=[z^8]\!\prod_{i=-7}^{30}\!\bigg(1+(z-1)\big(1-e^{-q^{2i-2\eta}}\big)\bigg),
\end{align*}
where, for each fixed $\eta$, the latter coefficient can easily be computed using Maple.
\end{rem}
\bth \label{EX1EX3} \textup{(see~\cite[Thm.~2,\,Thm.~3]{ABCKWW21})} Let
$L:=\ln(1/q)$ and $\chi:=2\pi\ii/L$. Then, as $n\to\infty$, the expectations of $X^{(n)}_i$, $i\in\{1,3\}$, satisfy
\begin{align}
\E X^{(n)}_1&\sim\frac{\ln (np^2)}{2L}+\frac12+\frac{\gamma}{2L}-\frac1{2L}\sum_{\ell\ne0}\Gamma\Big(\frac{\ell\chi}2\Big)
(np^2)^{-\ell\chi/2},\label{EX1}\\
\E X^{(n)}_3&\sim\frac{\ln^2(np^2)}{2L^2}+\Big[ \frac{\gamma}{L^2}+\frac1{2L} \Big]\ln(np^2)+\frac{\pi^2+6\gamma^2}{12L^2}+\frac{\gamma}{2L}-\frac1{12}\notag\\
&~~~~~~~~~~~~~~~-\frac{\ln (np^2)}{L^2}\sum_{\ell\ne0}\Gamma(\ell\chi)(np^2)^{-\ell\chi}
\label{EX2}\\
&~~~~~~~~~~~~~~~+\frac{1}{L^2}\sum_{\ell\ne0}\Gamma'(\ell\chi)(np^2)^{-\ell\chi}-\frac1{2L}\sum_{\ell\ne0}(-1)^\ell\Gamma\Big(\frac{\ell\chi}2\Big)
(np^2)^{-\ell\chi/2}.\notag
\end{align}
\ethGL

\section{The probability of avoiding certain pairs via Markov chains.}
\subsection {Two pairs $(i,i)$ and $(r,r)$ of identical letters}
The proofs of the theorems rest upon calculation of probabilities of avoiding certain pairs, which we will be doing by employing Markov chains. To illustrate that approach, we consider in greater detail the case of avoiding
two fixed pairs $(i,i)$ and $(r,r)$, where $i\ne r$, in a sequence of length $n$.
No distinction of letters different from $i,r$ is necessary, so for our Markov chain we can use a finite state space $S:=\{e,i,r,\Delta\}$, where $e:=\mathbb{N}\setminus \{i,r\}$ stands for ``everything else'', i.e., the set
$\mathbb{N}\setminus \{i,r\}$ is lumped together, and $\Delta$ denotes an additional cemetery state. The corresponding state diagram is
$$\includegraphics[width=5.5cm]{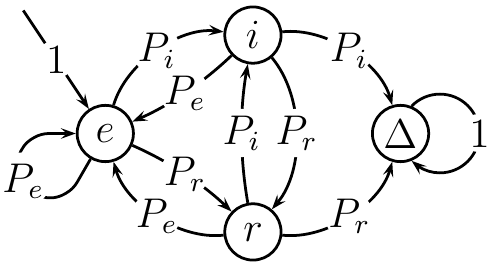}$$
From any realization $(z_k)_{k\ge1}$ of the i.i.d.\ sequence $(Z_k)_{k\ge1}$ we obtain a trajectory $(y_k)_{k\ge0}$ of this finite state Markov chain via
$$(y_0,y_1,y_2,\ldots)=(e,\phi(z_1),\phi(z_2),\phi(z_3),\ldots),$$
where $\phi(z_k):=\Delta$ if for some $j<k$ we have $(z_j,z_{j+1})\in\{(i,i),(r,r)\}$, and otherwise $$\phi(z_k)=\begin{cases}z_k,&z_k\in\{i,r\}\\e,&\textup{else}.\end{cases}$$
Example: If $n=8,i=1,r=2$ then the sequences $(1,2,3,1,2,3,4,5)$ and $(3,2,1,1,4,3,2,2)$ yield trajectories $(e,1,2,e,1,2,e,e,e)$ and $(e,e,2,1,\Delta,\Delta,\Delta,\Delta,\Delta)$.

Those trajectories $(y_k)_{k=0}^{n}$ satisfying $y_n\ne\Delta$ are in correspondence to sequences $(z_k)_{k=1}^n$ that avoid the pairs $(i,i)$ and $(r,r)$. Using the transition matrix
$$\Pi:=\begin{bmatrix}P_e&P_i&P_r&0\\P_e&0&P_r&P_i\\P_e&P_i&0&P_r\\0&0&0&1
\end{bmatrix},$$
where $P_e:=1-P_i-P_r$, the sought probability is $[1,0,0,0]\,\Pi^n\,[1,1,1,0]^t$, respectively, using the restriction $\bar\Pi$ of $\Pi$ to $\{e,i,r\}$, i.e.,
$$\bar\Pi:=\begin{bmatrix}P_e&P_i&P_r\\P_e&0&P_r\\P_e&P_i&0
\end{bmatrix},$$
and
initial probability $\pi(\cdot):=[1,0,0]$ and column vector of all ones $\mathds{1}$, that probability is
$$\P(X_{i,i}^{(n)}=X_{r,r}^{(n)}=0)=\pi\,\bar\Pi^n\,\mathds{1}.$$
A bound on such probability will now be derived in the following more general context.

We fix a finite non-empty set of forbidden pairs $$\mathcal{I}:=\{(k_i,m_i):i\in I\}$$ of size $|I|$, and
let $$J:=\bigcup_{i\in I}\{k_i,m_i\}=\{j_1,\ldots,j_{|J|}\},$$ where
$j_1<\ldots<j_{|J|}$. Moreover we fix $0<\delta\le 1/2$ and let
$$\mathcal{D}^{J}_\delta:=\{\bx\in\R^{|J|}:x_j\ge0\mbox{ for }j\in J,\,\sum_{j\in J}x_j\le 1-\delta\}.$$
\begin{lem}\label{exp3ClaPhi}
Let $\eps:=\sum_{i\in I}P_{k_i}P_{m_i}$. Then
\beq \label{exp3}\P(X_{k_i,m_i}^{(n)}=0,i\in I)\le
\delta^{-1/2}e^{-\eps n/2}\eeq
holds for $(P_j)_{j\in J}\in\mathcal{D}^{J}_\delta$.
Furthermore, there are functions $\la_1,C_1$ and $\Phi_n, n\ge1$, depending on $P_j,j\in J$,
that are $C^\infty$ and positive on an open set $\mathcal{F}$ satisfying
$\mathcal{D}^{J}_\delta\subseteq\mathcal{F}$, such that
\beq\label{ClaPhi}\P(X_{k_i,m_i}^{(n)}=0,i\in I)=C_1\la_1^n\Phi_n.\eeq
\end{lem}
\begin{rem}At several places we take the liberty to regard $(P_j)_{j\ge1}$ as variables (which is slight abuse of notation), to the effect, that several results in this section hold more generally also for strings of random variables with a distribution different from the geometric.
The reader must be prepared to see expressions involving $\lim_{P_j\to0}$, $\frac{\partial}{\partial P_j}$, and functions of $(P_j)_{j\in J}$ being $C^\infty$ in some domain, etc.\ all the time. In particular, we allow $(P_j)_{j\in J}$ to vary within the set
$\mathcal{D}^{J}_\delta$ above, which is a proper subset of the unit simplex of dimension $|J|$, because some of our results require $P_e$ to be bounded away from zero.
\end{rem}\goodbreak
{\bf Proof of Lemma~\ref{exp3ClaPhi}.}
Assume $P_j>0$ for $j\in J$, as well as $P_e:=1-\sum_{j\in J}P_j\ge\delta$.
Note that $\eps\le\sum_{j,\ell\in J}P_jP_{\ell}=(1-P_e)^2\le(1-\delta)^2<1-\delta$.
Define the matrix $\bar\Pi$ with rows and columns indexed by the set $J\cup\{e\}$ (which we assume ordered, starting with $e$ and followed by the elements of $J$ in ascending order) via
$$\bar\Pi_{k,m}:=\begin{cases}0,&(k,m)\in\mathcal{I},\\P_m,&\textup{else.}\end{cases}$$
We define a row vector
$\bw:=[\ \sqrt{P_e},\sqrt{P_{j_1}},\ldots,\sqrt{P_{j_{|J|}}}\ ]$, satisfying $\|\bw\|_2=1$, and a diagonal matrix $S:=\textup{Diag}(\bw)$, and the matrix
$$\hat\Pi:=S\bar\Pi S^{-1}=S\lb\mathds{1}\mathds{1}^t-\sum_{i\in I}\bbe_{k_i}\bbe_{m_i}^t\rb S,$$
where the column vectors $\bbe_j,j\in J\cup\{e\}$, denote the standard unit vectors in $\R^{|J|+1}$,
and observe, using the Frobenius norm $\|\hat\Pi\|_F=\sqrt{\sum_{k,m\in J\cup\{e\}}\hat\Pi_{k,m}^2}=\sqrt{1-\sum_{i\in I}P_{k_i}P_{m_i}}$, and $\pi=(\bbe_e)^t$,
\begin{align*}\mathbb{P}(X_{k_i,m_i}^{(n)}=0,i\in I)&=\pi\bar\Pi^n\mathds{1}
=\bw\hat\Pi^{n-1}\bw^t\\&\le\|\bw\|_2^2\|\hat\Pi\|_2^{n-1}\le\|\hat\Pi\|_F^{n-1}=(1-\eps)^{(n-1)/2}\le
\delta^{-1/2}(1-\eps)^{n/2}\le \delta^{-1/2}e^{-\eps n/2}.\end{align*}
Observe that $\bar\Pi$ is non-negative and primitive, therefore, by the Perron-Frobenius Theorem (see~\cite{SE81}), there is a unique positive eigenvalue $\la_1$, that is strictly larger in modulus than any other eigenvalue, and corresponding strictly positive left and right eigenvectors
$\bu$ and $\bv$, such that $\bar\Pi^n=\frac{\la_1^n}{\bu\bv}\bv\bu+\BO(n^{|J|}|\la_2|^n)$ element-wise,
where $\la_2$ is an eigenvalue of second largest modulus. This leads to
$$\mathbb{P}(X_{k_i,m_i}^{(n)}=0,i\in I)=\frac{(\pi\bv)(\bu\mathds{1})}{\bu\bv}\la_1^n+\BO(n^{|J|}|\la_2|^n).$$
By setting one or more of $(P_j)_{j\in J}$ to zero, one or more of the non-dominant eigenvalues $(\la_k)_{k\ge2}$ become zero, but there is a non-negative primitive submatrix constructed from the non-zero columns (and corresponding rows) of $\bar\Pi$, guaranteeing a unique positive eigenvalue larger in modulus than all other eigenvalues. As the row and column corresponding to state $e$ will always be part of that submatrix, the first components $u_e$ and $v_e$ of $\bu$ and of $\bv$ will be positive. By continuity, these properties also hold in a neighbourhood of such $(P_j)_{j\in J}$, which yields $\la_1$ being $C^\infty$ in some open superset $\bar{\mathcal{F}}$ of $\mathcal{D}^J_\delta$, by the implicit function theorem,
using the facts that the characteristic polynomial $p(\la)$ of $\bar\Pi$, considered as a function of $(\la,(P_j)_{j\in J})$, is $C^\infty$,
and the derivative of $p(\la)$ evaluated in a simple zero $\la_1$ is non-zero. On the set $\bar{\mathcal{F}}$, the components of $\frac1{u_e}\bu$ and $\frac1{v_e}\bv$ are $C^\infty$ functions of $(P_j)_{j\in J}$ as well.

We let $C_1:=\frac{(\pi\bv)(\bu\mathds{1})}{\bu\bv}$ and
$\Phi_n:=\frac1{C_1}\la_1^{-n}\mathbb{P}(X_{k_i,m_i}^{(n)}=0,i\in I)$.
Those are positive $C^\infty$ functions of $(P_j)_{j\in J}$
on an open set $\mathcal{F}$, satisfying $\mathcal{D}^J_\delta\subseteq\mathcal{F}\subseteq\bar{\mathcal{F}}$, the further restriction made necessary by the need to avoid $\bu\bv\le0$, which may occur for $(P_j)_{j\in J}$ outside $\mathcal{D}^J_\delta$.
Note that primitivity of $\bar\Pi$ may cease to hold when $P_e=0$. Moreover note that $|\la_2|$ is continuous on $\mathcal{D}_\delta^J$, but need not be differentiable on that set.
\epr
The bound \eqref{exp3} fits our needs when $\eps$ is large. Equation \eqref{ClaPhi} is useful in the case of small $\eps$, if asymptotics of $\la_1,C_1$ and $\Phi_n$ are known. In order to derive such asymptotics,
we let $\bar{\bar\Pi}$ be the matrix obtained from $\bar\Pi$ by deleting row and column corresponding to state $e$. Left and right eigenvectors $\bu=[1,\beta]$ and $\bv=[1/P_e,\mu^t]^t$, with row vector $\beta=(\beta_j)_{j\in J}$ and column vector $\mu=(\mu_j)_{j\in J}$, corresponding to the dominant eigenvalue $\lambda_1$ of $\bar\Pi$, lead to equations
\begin{align}
\lambda_1&=P_e(1+\sum_{j\in J}\beta_j)=P_e(1+\sum_{j\in J}P_j\mu_j),\label{sysla}\\
\beta&=\frac1{\lambda_1}\left[\beta\bar{\bar\Pi}+\bar\bp\right],\label{sysbe}\\
\mu&=\frac1{\lambda_1}\left[\bar{\bar\Pi}\mu+\mathds{1}\right],\label{sysmu}
\end{align}
with row vector $\bar \bp=(P_j)_{j\in J}$, and with ascending order of indices in $\beta, \mu,\bar \bp$. We keep denoting the column vector of all ones of appropriate dimension by $\mathds{1}$, and express $C_1$ in terms of $\beta$ and $\mu$ as follows:
\beq C_1=\frac{(\pi \bv)(\bu\mathds{1})}{\bu\bv}=
\frac{\frac1{P_e}(1+\sum_{j\in J}\beta_j)}
{\frac1{P_e}+\sum_{j\in J}\beta_j\mu_j}=
\frac{1+\beta\mathds{1}}
{1+P_e\beta\mu}.\label{C1gen}\eeq
Asymptotics up to any fixed order $K$ of
$\lambda_1,\beta,\mu$ are conveniently computed via fixed point iteration as described by the following algorithm:

\begin{algorithm}[H]\caption{Calculate asymptotics of
$\lambda_1,\beta,\mu$ up to fixed order.}
\label{alg1}
\begin{algorithmic}\REQUIRE $K \geq 0, k=0, \bar{\bar\Pi}, \bar \bp, \lambda=1,\bar\beta=[0,\ldots,0], \bar\mu=[1,\ldots,1]^t$
\WHILE{$k< K$}
\STATE $\bar\beta \leftarrow \frac1\lambda\big[\bar\beta\bar{\bar\Pi}+\bar\bp\big]$
\STATE $\lambda \leftarrow P_e[1+\bar\beta\mathds{1}]$
\STATE $\bar\mu \leftarrow \frac1\lambda\big[\bar{\bar\Pi}\bar\mu+\mathds{1}\big]$
\STATE $k \leftarrow k+1$
\ENDWHILE
\RETURN $\lambda,\bar\beta, \bar\mu$
\end{algorithmic}\end{algorithm}
The output $\la,\bar\beta, \bar\mu$ of the algorithm then satisfies $\la_1=\la+\BO^*_{K+1}$, $\beta=\bar\beta+\BO^*_{K+1}$, $\mu=\bar\mu+\BO^*_{K+1}$.
Here and in the following the notation $\BO_k^*$ always refers to the variables $(P_j)_{j\in J}$, but not to $P_e$. So, for instance,
$\BO_4^*$ is the same as $\BO(\gamma^4)$, where $\gamma=\sum_{j\in J}P_j$.

A few words on justification of the algorithm: First note, that nothing changes if the line
$\lambda \leftarrow P_e[1+\bar\beta\mathds{1}]$ is replaced by
$\lambda \leftarrow P_e[1+\bar \bp\bar\mu]$. This is seen to hold for
$k=0$, where $\bar\beta=\bar \bp$ has already been updated, but $\bar\mu=\mathds{1}$ has not, and for $k>0$ by a simple induction step. We can thus see Algorithm~1 as a combination of two algorithms, one of them only updating the pair $(\bar\beta,\la)$, the other only updating the pair
$(\la,\bar\mu)$, with those algorithms having identical updates of $\la$.
Let us concentrate on the latter algorithm. Denote $x=(\la,\bar\mu)$ and let $\bold{0}$ be the zero vector of appropriate dimension. Observe that the function
$F(x,\bar\bp)=\lb\begin{smallmatrix}
\la-P_e[1+\bar \bp\bar\mu]\\
\bar\mu-\frac1\lambda\big[\bar{\bar\Pi}\bar\mu+\mathds{1}\big]
\end{smallmatrix}\rb$
is $C^\infty$ in a neighbourhood of $(x_0,\bar\bp_0):=(1,\mathds{1},\bold{0})$, with $F(x_0,\bar\bp_0)=\bold{0}$. Now the Jacobian $J F(x_0,\bar\bp_0)$ is nonsingular, so there is a unique $C^\infty$ function
$x(\bar\bp)=(\la_1(\bar\bp),\mu(\bar\bp))$ defined in some neighbourhood $\mathcal{V}$ of $\bar\bp=\bold{0}$, satisfying
$x(\bold{0})=x_0$ and $F(x(\bar\bp),\bar\bp)=\bold{0}$ for $\bar\bp\in\mathcal{V}$, by the implicit function theorem. Denoting iterates by\begin{center}
$\la^{k+1}=f(\bar\mu^k)=P_e[1+\bar \bp\bar\mu^k]$\quad and \quad
$\bar\mu^{k+1}=\frac1{\la^{k+1}}g(\bar\mu^k)
=\frac1{\la^{k+1}}[\bar{\bar\Pi}\bar\mu^k+\mathds{1}]$,\end{center}
with $\la^0=1$ and $\bar\mu^0=\mathds{1}$, we can easily check $|\la_1(\bar\bp)-\la^0|=\BO_1^*$ and
$\|\mu_1(\bar\bp)-\bar\mu^0\|=\BO_1^*$, for $\bar\bp\in\mathcal{V}$.\\
Assume now that we have already shown $|\la_1(\bar\bp)-\la^{k-1}|=\BO_k^*$ and
$\|\mu(\bar\bp)-\bar\mu^{k-1}\|=\BO_k^*$. Then we have
$|\la_1(\bar\bp)-\la^k|=|f(\mu(\bar\bp))-f(\mu^{k-1})|=|P_e\bar\bp(\mu(\bar\bp)-\mu^{k-1})|
=\BO_1^*\|\mu(\bar\bp)-\bar\mu^{k-1}\|=\BO_{k+1}^*$, and
$\|\mu(\bar\bp)-\bar\mu^k\|
=\|\frac1{\la_1}g(\mu(\bar\bp))-\frac1{\la^{k}}g(\bar\mu^{k-1})\|
\le\|(\frac1{\la_1}-\frac1{\la^{k}})g(\mu(\bar\bp))\|
+\|\frac1{\la^{k}}(g(\mu(\bar\bp))-g(\bar\mu^{k-1}))\|=\BO_{k+1}^*$,
because $(\frac1{\la_1}\!-\!\frac1{\la^{k}})g(\mu(\bar\bp))
\!=\!\frac{\la^{k}-\la_1}{\la^{k}}\mu(\bar\bp)\!=\!\BO_{k+1}^*$,
and $\frac1{\la^{k}}(g(\mu(\bar\bp))\!-\!g(\bar\mu^{k-1}))
\!=\!\frac1{\la^{k}}\bar{\bar\Pi}(\mu(\bar\bp)\!-\!\bar\mu^{k-1})\!=\!\BO_{k+1}^*$.

The next lemma provides asymptotics of probabilities in the case of a single avoided pair.
\begin{lem}\label{Piiir}
The probabilities of avoiding the pair $(i,i)$, resp.\ $(i,r)$ for $i\ne r$,
in a sequence of length~$n$ satisfy
\begin{align}
\P(X_{i,i}^{(n)}=0)&=e^{-nP_i^2}+\BO\big(\sqrt{n}P_i^2e^{-\frac{2-\delta}4 nP_i^2}\big),\label{Pii}\\
\P(X_{i,r}^{(n)}=0)&=e^{-nP_iP_r}+\BO\big(P_iP_r e^{-\frac {2-\delta}4 nP_iP_r}\big),\label{Pir}
\end{align}
as $n\to \infty$, uniformly for $P_i\in\mathcal{D}^{\{i\}}_\delta$, resp.\ for $(P_i,P_r)\in\mathcal{D}^{\{i,r\}}_\delta$.
\end{lem}
\bpr
We first consider the forbidden pair $(i,i)$.
The matrix $\bar\Pi$, its characteristic polynomial $p$, and asymptotics of $\la_1$ and $C_1$ are given by
$$
\bar\Pi=\begin{bmatrix}P_e&P_i\\
P_e&0\\
\end{bmatrix},\qquad
\begin{aligned}
p(\la)&=\la^2-(1-P_i)\la-P_iP_e,\\
\la_1&=1-P_i^2+P_i^3-2P_i^4+\BO_5^*,\\
C_1&=1+P_i^2-2P_i^3+6P_i^4+\BO_5^*,
\end{aligned}$$
where we used Algorithm~1 (with $K=4$) and \eqref{C1gen}.

Following a suggestion by \cite{SA20}, we can easily derive $\la_1$ from  $p(\la)$, after replacing $P_e$ by $1-P_i$. We add  an extra variable $v$, carrying the weight of the $P_.:\Pt_.:=vP_.  $. We have the local expansion of the solution at $0$ by using the Maple package gfun (see \cite{SZ94}):
\[sol:=\mbox{gfun}[\mbox{algeqtoseries}](p(\la),v,\la,pr),\]
where $pr$ denotes the precision of the expansion into $v$. We obtain the solutions as  $sol[1],sol[2]$ and we keep the solution close to $1$.

Denoting $\pi\bar\Pi^n\mathds{1}=C_1\la_1^n+C_2\la_2^n$, with $\la_2$ the non-dominant eigenvalue of $\bar\Pi$, we have $C_1\la_1^0+C_2\la_2^0=1$, and
therefore $C_2=-P_i^2+2P_i^3-6P_i^4+\BO_5^*$, which leads to
$\Phi_n=1+\frac{C_2}{C_1}\big(\frac{\la_2}{\la_1}\big)^n=1+\BO(P_i^2)$, uniformly in $n$.
This is used in \eqref{ClaPhi}, together with $C_1=1+\BO(P_i^2)$ and
$$\la_1^n=e^{n\ln \la_1}=e^{n(-P_i^2+P_i^3+\BO(P_i^4))}=e^{-nP_i^2}(1+\BO(nP_i^3))$$
leading to $\P(X_{i,i}^{(n)}=0)=\lb1+\BO(P_i^2)+\BO(nP_i^3)\rb e^{-nP_i^2}$,
 for $nP_i^3=\BO(1)$, resp.\ for $nP_i^2=\BO(n^{1/3})$. Note that for fixed $\alpha,\beta\!>\!0$
the function $x^\alpha e^{-\beta x}$ is bounded for $x>0$, implying
$$nP_i^3e^{-nP_i^2}=\sqrt{n}P_i^2(nP_i^2)^{1/2}e^{-\frac {2+\delta}4 nP_i^2}e^{-\frac {2-\delta}4 nP_i^2}=\BO\Big(\sqrt{n}P_i^2e^{-\frac {2-\delta}4 nP_i^2}\Big).$$ Moreover also $P_i^2e^{-nP_i^2}=\BO\Big(\sqrt{n}P_i^2e^{-\frac {2-\delta}4 nP_i^2}\Big)$
holds, and \eqref{exp3} can be built in by observing that $nP_i^2=\Omega(n^{1/3})$ implies
$\delta^{-1/2}e^{-\frac n2P_i^2}=\BO\Big(\sqrt{n}P_i^2e^{-\frac {2-\delta}4 nP_i^2}\Big)$.
We have thus obtained \eqref{Pii}.

We now consider the forbidden pair $(i,r)$ with $i\ne r$. The matrix $\bar\Pi$, its characteristic polynomial $p$, and asymptotics of $\la_1$ and $C_1$ are given by
$$
\bar\Pi=\begin{bmatrix}P_e&P_i&P_r\\
P_e&P_i&0\\P_e&P_i&P_r\\
\end{bmatrix},\qquad
\begin{aligned}
p(\la)&=\la^3-\la^2+P_iP_r\la,\\
\la_1&=1-P_iP_r-P_i^2P_r^2-2P_i^3P_r^3+\BO_8^*,\\
C_1&=1+P_iP_r+3P_i^2P_r^2+10P_i^3P_r^3+\BO_8^*.
\end{aligned}$$

Clearly, $\la_1$, and therefore also $C_1$ and $\Phi_n$, are $C^\infty$ functions of the coefficient $P_iP_r$ of the characteristic polynomial $p$,
meaning that the error term $\BO_8^*$ is in fact $\BO(P_i^4P_r^4)$.
Sufficiently accurate for our purposes are the asymptotics
$\la_1=1-P_iP_r+\BO(P_i^2P_r^2)$ and $C_1=1+\BO(P_iP_r)$.

One of the eigenvalues is $0$, therefore a representation
$\pi\bar\Pi^n\mathds{1}=C_1\la_1^n+C_2\la_2^n$ as before also holds in this case, with $C_2=\BO(P_iP_r)$, and $\Phi_n=1+\frac{C_2}{C_1}\big(\frac{\la_2}{\la_1}\big)^n=1+\BO(P_iP_r)$, uniformly in $n$.
All this, together with $\la_1^n=e^{-nP_iP_r}(1+\BO(nP_i^2P_r^2))$, leads to \eqref{Pir} via \eqref{ClaPhi}, taking care of error terms as above.
\epr
The next corollary follows easily from equations \eqref{exp3}, \eqref{Pii} and \eqref{Pir}.
\bcor
The variances of $X_{i,i}^{(n)}$ and $X_{i,r}^{(n)}$ for $i\ne r$ satisfy
\begin{align}
\Var X_{i,i}^{(n)}&=e^{-nP_i^2}-e^{-2nP_i^2}
+\BO\big(\sqrt{n}P_i^2e^{-\frac {2-\delta}4 nP_i^2}\big),\label{Vii}\\
\Var X_{i,r}^{(n)}&=e^{-nP_iP_r}-e^{-2nP_iP_r}+\BO\big(P_iP_r e^{-\frac {2-\delta}4 nP_iP_r}\big),\label{Vir}
\end{align}
as $n\to \infty$, uniformly for $P_i\in\mathcal{D}^{\{i\}}_\delta$, resp.\ for $(P_i,P_r)\in\mathcal{D}^{\{i,r\}}_\delta$.
\ecor
In order to obtain asymptotics for the covariance
\begin{align*}\Cov(X_{i,i}^{(n)},X_{r,r}^{(n)})&=\P(X_{i,i}^{(n)}\!=\!X_{r,r}^{(n)}\!=\!1)-\P(X_{i,i}^{(n)}\!=\!1)\P(X_{r,r}^{(n)}\!=\!1)\\
&=\P(X_{i,i}^{(n)}\!=\!X_{r,r}^{(n)}\!=\!0)-\P(X_{i,i}^{(n)}\!=\!0)\P(X_{r,r}^{(n)}\!=\!0),\end{align*}
we need the following result.
\begin{lem}\label{Frobineq}
Let $A\in\R^{k\times k}$, with $k\ge2$, have spectral radius $\rho(A)\le1$ and
Frobenius norm $\|A\|_F=C'$. Then, with $C:=\max(C',k)$ and $C'':=2C^{k-1}$, we have
$$\|A^n\|_F\le C''n^{k-1}.$$
\end{lem}
\bpr
We use Schur decomposition, according to which there is a unitary matrix $Q$ such that $\bar A:=QAQ^{-1}$ is upper triangular and satisfies $\rho(\bar A)=\rho(A)$ and $\|\bar A\|_F=\|A\|_F$.
Then also
$$\|A^n\|_F=\|Q^{-1}\bar A^nQ\|_F=\|\bar A^n\|_F.$$
Moreover $\rho(\bar A^n)=\rho(A^n)\le1$, and $\bar A^n$ being triangular, we deduce
$|(\bar A^n)_{i,i}|\le1$. Regarding off diagonal elements of $\bar A^n$, we have
\begin{equation}\label{e:frob}|(\bar A^n)_{i,i+\ell}|\le \sum_{j=1}^\ell\binom{\ell-1}{j-1}\binom{n}{j}
\lp\frac C{\sqrt{j}}\rp^j,\end{equation}
as we now show. Note that $(\bar A^n)_{i,i+\ell}$ is a sum of products
$\bar a_{i_0,i_1}\cdot\bar a_{i_1,i_2}\cdots\bar a_{i_{n-1},i_n}$,
where the sum extends over all sequences $(i_k)_{k=0}^n$ that are increasing with $i_0=i$ and $i_n=i+\ell$. Such a sequence has at least one and at most $\ell$ jumps. For $j$ satisfying $1\le j\le\ell$, there
are $\binom{\ell-1}{j-1}$ ways to accommodate $j$ jump heights $(h_m)_{m=1}^j$, and for each of those there are $\binom{n}{j}$ ways to position those $j$ jumps.
In terms of cumulated jump heights $H_m:=i+\sum_{\mu=1}^mh_\mu,\,0\le m\le j$, we can rewrite above product as
$$\bar a_{i_0,i_1}\cdot\bar a_{i_1,i_2}\cdots\bar a_{i_{n-1},i_n}
=\bar{\bar a}\cdot\bar a_{H_0,H_1}\cdot\bar a_{H_1,H_2}\cdots\bar a_{H_{j-1},H_j},$$
where $\bar{\bar a}$ is a product of $n-j$ diagonal elements of $\bar A$, and therefore satisfies $|\bar{\bar a}|\le1$.
Furthermore, $\sum_{m=1}^j|\bar a_{H_{m-1},H_m}|^2\le\|\bar A\|_F^2\le C^2$, so by observing that the product $\prod_{m=1}^j|\bar a_{H_{m-1},H_m}|^2$ is maximized, if its terms are all equal to $\frac{C^2}j$, we obtain
$|\bar a_{i_0,i_1}\cdot\bar a_{i_1,i_2}\cdots\bar a_{i_{n-1},i_n}|\le
\lp\frac C{\sqrt{j}}\rp^j$, so \eqref{e:frob} is proven.\\
Since $C\ge k$ ensures that $\big(\big(C/\sqrt{\ell}\!\ \big)^\ell\big)_{1\le\ell<k}$ is increasing, we can extend the estimate \eqref{e:frob},
$$|(\bar A^n)_{i,i+\ell}|\le
\sum_{j=1}^\ell\binom{\ell-1}{j-1}\binom{n}{j}\lp\frac C{\sqrt{\ell}}\rp^\ell
=\lp\frac C{\sqrt{\ell}}\rp^\ell\binom{n+\ell-1}{\ell}\le
\lp\frac C{\sqrt{k-1}}\rp^{k-1}\binom{n+k-2}{k-1}
,$$
for $1\le i< i+\ell\le k$. We obtain
$$\|\bar A^n\|_F\le k\binom{n+k-2}{k-1}\lp\frac {C}{\sqrt{k-1}}\rp^{k-1}\le2C^{k-1}n^{k-1},
$$
because of $\binom{n+k-2}{k-1}\le n^{k-1}$ for $k\ge2$ and $n\ge1$,
and because of $\max\limits_{k\ge2}\frac{k}{(k-1)^{(k-1)/2}}=2$,
which completes the proof.
\epr
We now turn to asymptotics of covariances.
\begin{lem}\label{lemCovir} For $i\ne r$ and $(P_i,P_r)\in\mathcal{D}^{\{i,r\}}_\delta$ we have, for $nP_iP_r(P_i+P_r)^2=\BO(1)$,
\beq\label{CovirBO}
\Cov(X_{i,i}^{(n)},X_{r,r}^{(n)})=\BO\Big(P_iP_r+nP_iP_r(P_i+P_r)^2\Big)\P(X_{i,i}^{(n)}=0)\P(X_{r,r}^{(n)}=0).
\eeq
\end{lem}
\bpr
We first find asymptotics of $\la_1$ and $C_1$ from $\P(X_{i,i}^{(n)}=X_{r,r}^{(n)}=0)=C_1\la_1^n\Phi_n$, proceeding as in the previous lemma. The matrix $\bar\Pi$, its characteristic polynomial $p$, and asymptotics of $\la_1$ and $C_1$ are given by
$$
\bar\Pi=\begin{bmatrix}P_e&P_i&P_r\\
P_e&0&P_r\\P_e&P_i&0\\
\end{bmatrix},\qquad
\begin{aligned}
p(\la)&=\la^3-P_e\la^2-[P_e(P_i+P_r)+P_iP_r]\la-P_iP_rP_e,\\
\la_1&=1-P_i^2-P_r^2+P_i^3+P_r^3+\BO_4^*,\\
C_1&=1-P_i^2-P_r^2+2P_i^3+2P_r^3+\BO_4^*.
\end{aligned}$$
Again, we can also replace $P_e$ by $1-P_i-P_r$ and use gfun.
From Lemma~\ref{exp3ClaPhi} we know that $\la_1,C_1$ and $\Phi_n$ are $C^\infty$ functions of $P_i,P_r$ in some open superset $\mathcal{F}$ of $\mathcal{D}^{\{i,r\}}_\delta$, such that
\beq\label{e:pir}
\mathbb{P}(X_{i,i}^{(n)}=X_{r,r}^{(n)}=0)
=C_1(P_i,P_r)[\la_1(P_i,P_r)]^n\Phi_n(P_i,P_r)
\eeq
holds for $(P_i,P_r)\in\mathcal{D}^{\{i,r\}}_\delta$.
In fact, we will only need that those functions are $C^2$ in the following.

Note that $\P(X_{i,i}^{(n)}=0)$ can be obtained from \eqref{e:pir} as the limiting case $P_r\to0$.
Observe that we have \footnotesize{$$\lim_{P_i\to0}\frac{C_1(P_i,P_r)}{C_1(0,P_r)C_1(P_i,0)}=\lim_{P_r\to0}\frac{C_1(P_i,P_r)}{C_1(0,P_r)C_1(P_i,0)}=\lim_{P_i\to0}\frac{\Phi_n(P_i,P_r)}{\Phi_n(0,P_r)\Phi_n(P_i,0)}=\lim_{P_r\to0}\frac{\Phi_n(P_i,P_r)}{\Phi_n(0,P_r)\Phi_n(P_i,0)}=1,$$}\normalsize and therefore
$\frac{C_1(P_i,P_r)}{C_1(0,P_r)C_1(P_i,0)}=1+\BO(P_iP_r)$ and
$\frac{\Phi_n(P_i,P_r)}{\Phi_n(0,P_r)\Phi_n(P_i,0)}=1+\BO(P_iP_r)$.

To see that the latter holds uniformly in $n$ and $(P_i,P_r)\in\mathcal{D}^{\{i,r\}}_\delta$,
we start defining $\check\Pi:=\bar\Pi-\frac{\la_1}{\bu\bv}\bv\bu$,
so that $\bar\Pi=\frac{\la_1}{\bu\bv}\bv\bu+\check\Pi$, and
\beq\label{barlacheck}\pi\bar\Pi^n\mathds{1}=C_1\la_1^n+\pi\check\Pi^n\mathds{1},\eeq where we used that $\bu$ and $\bv$ are in the left resp.\ right kernel of the matrix $\check\Pi$.

Denoting the spectral radius of a
square matrix $A$ by $\rho(A)$, we clearly have $\rho(\check\Pi)=|\la_2|$, and since $\mathcal{D}^{\{i,r\}}_\delta$ is compact,
we have $\max_{(P_i,P_r)\in\mathcal{D}^{\{i,r\}}_\delta}\frac{|\la_2|}{\la_1}=:\kappa<1$.
All components of $\frac1{\la_1\kappa}\check\Pi$ are continuous, so there is a constant $C'$ such that $\|\frac1{\la_1\kappa}\check\Pi\|_F\le C'$ on
$\mathcal{D}^{\{i,r\}}_\delta$. By applying Lemma~\ref{Frobineq} below to the matrix
$\frac1{\la_1\kappa}\check\Pi$, we obtain
$$\Phi_n-1=\frac{\pi\check\Pi^n\mathds{1}}{C_1\la_1^n}=\BO(n^2\kappa^n)=
\BO(\bar\kappa^n),$$
for some $\kappa<\bar\kappa<1$, uniformly on $\mathcal{D}^{\{i,r\}}_\delta$.
Similarly, we obtain
$\frac{\partial\Phi_n}{\partial P_i}=\BO(\bar\kappa^n)$, $\frac{\partial\Phi_n}{\partial P_r}=\BO(\bar\kappa^n)$, and $\frac{\partial^2\Phi_n}{\partial P_i\partial P_r}\!=\!\BO(\bar\kappa^n)$, uniformly on $\mathcal{D}^{\{i,r\}}_\delta$, using, e.~g.,
$\frac{\partial\pi\check\Pi^n\mathds{1}}{\partial P_i}\!=\!
\sum_{0\le j< n}\pi\check\Pi^j\frac{\partial\check\Pi}{\partial P_i}\check\Pi^{n-1-j}\mathds{1}$, and again Lemma~\ref{Frobineq}.

Define
$\Psi_n(P_i,P_r):=\frac{\Phi_n(P_i,P_r)}{\Phi_n(0,P_r)\Phi_n(P_i,0)}-1$ and observe that $\lim_{n\to\infty}\frac{\partial^2}{\partial P_i\partial P_r}\Psi_n=0$ holds uniformly on $\mathcal{D}^{\{i,r\}}_\delta$. Note that we have $\Psi_n(P_i,0)=\Psi_n(0,P_r)=0$ for $0\le P_i,P_r\le 1-\delta$, yielding
$$\Psi_n(P_i,P_r)=\Psi_n(P_i,P_r)-\Psi_n(P_i,0)-\Psi_n(0,P_r)+\Psi_n(0,0)=P_iP_r\frac{\partial^2\Psi_n}{\partial P_i\partial P_r}(p_i,p_r)$$
by the (bivariate) Mean Value Theorem, where $0\le p_i\le P_i$ and $0\le p_r\le P_r$, see \cite[Thm.~9.40]{RU66}.
Defining $\bar C:=\max_{n\ge1}\max_{(p_i,p_r)\in\mathcal{D}^{\{i,r\}}_\delta}\left|\frac{\partial^2\Psi_n}{\partial P_i\partial P_r}(p_i,p_r)\right|$, we finally conclude $|\Psi_n(P_i,P_r)|\le\bar CP_iP_r$ for all $n\ge1$ and $(P_i,P_r)\in\mathcal{D}^{\{i,r\}}_\delta$, establishing the uniformity claim.
By our asymptotics for $\la_1$, we similarly obtain
$$\ln \la_1(P_i,P_r)-\ln \la_1(P_i,0)-\ln \la_1(0,P_r)+\ln \la_1(0,0)
=P_iP_r\frac{\partial^2\ln \la_1}{\partial P_i\partial P_r}(p_i,p_r)=
\BO(P_iP_r(P_i+P_r)^2),
$$
leading to
$$\frac{\la_1(P_i,P_r)}{\la_1(P_i,0)\la_1(0,P_r)}=1+\BO\big(P_iP_r(P_i+P_r)^2\big).$$
We summarize
$$
\frac{\mathbb{P}(X_{i,i}^{(n)}=X_{r,r}^{(n)}=0)}{\mathbb{P}(X_{i,i}^{(n)}=0)\mathbb{P}(X_{r,r}^{(n)}=0)}=1+\BO\big(P_iP_r+nP_iP_r(P_i+P_r)^2\big),
$$ finally arriving at \eqref{CovirBO}.
\epr
From \eqref{exp3} we derive $\Cov(X_{i,i}^{(n)},X_{r,r}^{(n)})= \BO\Big(e^{-\frac{n}2(P_i^2+P_r^2)}\Big)$, that together with \eqref{CovirBO}, where we use
\small{$$nP_iP_r(P_i+P_r)^2\P(X_{i,i}^{(n)}=0)\P(X_{r,r}^{(n)}=0)= \BO\Big(P_iP_r(nP_i^2+nP_r^2)e^{-\frac n2(P_i^2+P_r^2)}\Big)=\BO\Big(P_iP_re^{-\frac{2-\delta}4n(P_i^2+P_r^2)}\Big),$$}
\normalsize implies the next corollary, since $e^{-\frac{n}2(P_i^2+P_r^2)}=\BO\Big(P_iP_re^{-\frac{2-\delta}4n(P_i^2+P_r^2)}\Big)$, for $nP_iP_r(P_i+P_r)^2=\Omega(1)$.
\begin{cor}
For $i\ne r$, the covariance of $X_{i,i}^{(n)}$ and $X_{r,r}^{(n)}$ satisfies
\beq\label{Covirexp3}\Cov(X_{i,i}^{(n)},X_{r,r}^{(n)})=\BO\Big(P_iP_re^{-\frac {2-\delta}4 n(P_i^2+P_r^2)}\Big),\eeq
as $n\to \infty$, uniformly for $(P_i,P_r)\in\mathcal{D}^{\{i,r\}}_\delta$.
\end{cor}
\subsection {The variance of $X^{(n)}_1$}
In this subsection we use the results on variances and covariances in the case of avoided pairs of identical letters, that we have derived so far, to furnish a proof of equation \eqref{var1} of Theorem~\ref{thm:var123}.
\begin{lem}\label{varsumeqsumvar}
The variance of $X^{(n)}_1$ is asymptotically given by
$$\textup{Var}\,X^{(n)}_1=\sum_{i\ge1}\textup{Var}X_{i,i}^{(n)}+\BO\left(\frac1n\right)=S^{(n)}_1+\BO\lp\frac1{\sqrt{n}}\rp,$$
with $S^{(n)}_1$ given in \eqref{S1}. In particular the contribution of covariances is negligible.
\end{lem}
\bpr
Dealing with covariances first, note that \eqref{Covirexp3} guarantees that the double sum of covariances
$\sum_{i\ne r}\Cov(X_{i,i}^{(n)},X_{r,r}^{(n)})$
makes a negligible contribution to the variance of $X^{(n)}_1$:
We will use that
\beq\label{iialpha}\sum_{k\ge1}\lb nP_k^\beta\rb^\alpha e^{-nP_k^\beta}=\mathcal{O}(1)\ \mbox{ holds for }\alpha,\beta>0.\eeq
This follows from the following general result: If for some $c<1$ a set $\mathcal{P}=\{x_i:i\in\mathbb{N}\}$ satisfies $x_i>0$ and $\frac{x_{i+1}}{x_i}\le c$ for $i\in\mathbb{N}$, then $\sum_{x\in\mathcal{P}}x^\alpha e^{-x}<\infty$. For a proof observe that there is a constant $C_\alpha>0$ such that $x^\alpha e^{-x}\le\min(x^\alpha,C_\alpha x^{-\alpha})$ for $x>0$. Let $\bar x:=(C_\alpha)^{1/(2\alpha)}$. Then
$$\sum_{x\in\mathcal{P}}x^\alpha e^{-x}\le \sum_{x\in\mathcal{P}\cap\,]0,\bar x]}x^\alpha+\sum_{x\in\mathcal{P}\cap[\bar x,\infty[}C_\alpha x^{-\alpha}\le \bar x^\alpha\sum_{i\ge0}c^i+C_\alpha\bar x^{-\alpha}\sum_{i\ge0}c^i=2\frac{\sqrt{C_\alpha}}{1-c}.$$
With the help of \eqref{iialpha} we find
$$\sum_{i\ge1}\sum_{r\ge1}P_iP_r e^{-\frac {2-\delta}4n(P_i^2+P_r^2)}
=\frac1n\sum_{i\ge1}(nP_i^2)^{1/2}e^{-\frac {2-\delta}4nP_i^2}\sum_{r\ge1}(nP_r^2)^{1/2}e^{-\frac {2-\delta}4nP_r^2}=\BO\left(\frac1n\right).
$$
This leads to $\sum_{i\ne r}\Cov(X_{i,i}^{(n)},X_{r,r}^{(n)})=\BO(\frac1n)$.

We now turn to $\sum_{i\ge1}\textup{Var}X_{i,i}^{(n)}$.
Observe that the sum of error terms from \eqref{Vii} satisfies
$$\sum_{i\ge1}\sqrt{n}P_i^2 e^{-\frac {2-\delta}4 nP_i^2}=\BO\left(\frac1{\sqrt{n}}\right),$$ by \eqref{iialpha}. Therefore, up to an error term
$\BO\left(\frac1{\sqrt{n}}\right)$,
the variance $\Var X^{(n)}_1$ equals
$$\sum_{i\ge1}\lb e^{-nP_i^2}-e^{-2nP_i^2}\rb=\sum_{i\ge1}\lb 1-e^{-2nP_i^2}\rb-\sum_{i\ge1}\lb 1-e^{-nP_i^2}\rb=G(2np^2)-G(np^2),$$
which can be evaluated using $G$ from Appendix~\ref{A1}, directly leading to $S^{(n)}_1$ from \eqref{S1}.
\epr
\subsection {Contribution of covariances to the variance of $X^{(n)}_2$}
In this subsection we will prove the following lemma, which will also imply equations \eqref{var2} and \eqref{var3} of Theorem~\ref{thm:var123}.
\begin{lem}\label{var2sumeqsumvar2}The variance of $X^{(n)}_2$ is asymptotically given by
\beq\label{varX2}\textup{Var}\,X^{(n)}_2=\sum_{i,j\ge1}\textup{Var}X_{i,j}^{(n)}
+2\!\!\sum_{i,j,k\ge1}H(i,j,k)
+\BO\left(\frac{\ln n}{\sqrt{n}}\right)
=S^{(n)}_2+T^{(n)}_2+\BO\left(\frac{\ln n}{\sqrt{n}}\right),\eeq
where $H(i,j,k)=(e^{nP_iP_jP_k}-1)e^{-nP_iP_j-nP_jP_k}$, and $S^{(n)}_2,T^{(n)}_2$ are given in \eqref{S2} and \eqref{T2}.
Only covariances
$\Cov(X_{i,j}^{(n)},X_{j,k}^{(n)})$, resp.~$\Cov(X_{j,i}^{(n)},X_{k,j}^{(n)})$, with $i,j,k$ all different, and
$\Cov(X_{i,j}^{(n)},X_{j,i}^{(n)})$ with $i,j$ different, contribute significantly to $\Var X^{(n)}_2$.
\end{lem}
\bpr We start
considering distinct forbidden pairs $(i_1,j_1),\,(i_2,j_2)$, where we allow
$i_1\ne j_1$ or $i_2\ne j_2$ or both, and are again interested in  negligibility of covariance contributions.

Let $J:=\{i_1,j_1,i_2,j_2\}$, and assume $P_i>0$ for $i\in J$, as well as $P_e:=1-\sum_{i\in J}P_i\ge\delta$. Define the matrix $\bar\Pi$ with rows and columns indexed by the set $J\cup\{e\}$ (which we assume ordered, starting with $e$ and followed by the elements of $J$ in ascending order) via
$$\bar\Pi_{i,j}:=\begin{cases}0,&(i,j)\in\{(i_1,j_1),(i_2,j_2)\},\\P_j,&\textup{else,}\end{cases}$$

We will have to distinguish several cases, which however share some common features: The sought probability can be expressed as
$$\mathbb{P}(X_{i_1,j_1}^{(n)}=X_{i_2,j_2}^{(n)}=0)=\pi\bar\Pi^n\mathds{1}=
C_1\la_1^n\Phi_n,$$
where, as previously observed, $\la_1$, $C_1$ and $\Phi_n$ for $n\ge1$ are $C^\infty$ functions on an open superset of $\mathcal{D}^J_\delta$.
Limits $\lim_{n\to\infty}\Phi_n=1$, $\lim_{n\to\infty}\frac{\partial \Phi_n}{\partial P_{i_1}}=0$, etc., will again be uniform for
$(P_i)_{i\in J}\in\mathcal{D}^J_\delta$.
Denoting
$$\mathbb{P}(X_{i_1,j_1}^{(n)}=0)=
C_*\la_*^n\Phi_n^*,\qquad \mathbb{P}(X_{i_2,j_2}^{(n)}=0)=
C_\circ\la_\circ^n\Phi_n^\circ,$$
we observe
$$\lim_{P_i\to0}\frac{\la_1}{\la_*\la_\circ}=\lim_{P_i\to0}\frac{C_1}{C_*C_\circ}=\lim_{P_i\to0}\frac{\Phi_n}{\Phi_n^*\Phi_n^\circ}=1,\quad
\textup{for}\ i\in J,$$
leading to $\frac{\la_1}{\la_*\la_\circ}=1+\BO\left(\prod_{i\in J}P_i\right)$, $\frac{C_1}{C_*C_\circ}=1+\BO\left(\prod_{i\in J}P_i\right)$,
and $\frac{\Phi_n}{\Phi_n^*\Phi_n^\circ}=1+\BO\left(\prod_{i\in J}P_i\right)$, with implied constant independent of $n$.
(This independence can be shown as in the proof of Lemma~\ref{lemCovir}.)
As we will see,
more accurate representations for
$\la_1$, complementing those obtained by Algorithm~1, can always be found in the form
$$\la_1=1-P_{i_1}P_{j_1}-P_{i_2}P_{j_2}+Q+\BO_4^*,$$
where $Q=\BO_3^*$ and $Q\ge0$. We will observe, that in each of the cases
\beq\label{Qexpr}Q=\sum\nolimits_{i,r,t:(i,r),(r,t)\in\{(i_1,j_1),(i_2,j_2)\}}P_iP_rP_t\eeq holds.
Using $\la_*=1-P_{i_1}P_{j_1}+Q_*+\BO(P_{i_1}^2P_{j_1}^2)$ and $\la_\circ=1-P_{i_2}P_{j_2}+Q_\circ+\BO(P_{i_2}^2P_{j_2}^2)$ (depending on whether
$(i_1,j_1)=(i,i)$ or $(i,r)$, we have $Q_*=P_i^3$ or
$Q_*=0$, and similarly for $Q_\circ$, see the proof of Lemma~\ref{Piiir}),
we will obtain in most of the cases
\beq\label{barQerror}\ln\frac{\la_1}{\la_*\la_\circ}=Q-Q_*-Q_\circ+\BO(P_{i_1}P_{j_1}P_{i_2}P_{j_2}),\eeq
where the error term needs justification in each of these cases. In some cases this is done by employing the MVT, as in the proof of Lemma~\ref{lemCovir}.
This results in the following expression for a quotient of probabilities, that directly leads to an expression for the covariance, where we denote $\bar Q:=Q-Q_*-Q_\circ$,
$$\frac{\mathbb{P}(X_{i_1,j_1}^{(n)}=X_{i_2,j_2}^{(n)}=0)}{\mathbb{P}(X_{i_1,j_1}^{(n)}=0)\mathbb{P}(X_{i_2,j_2}^{(n)}=0)}
=\left(\frac{\la_1}{\la_*\la_\circ}\right)^n\frac{C_1}{C_*C_\circ}\frac{\Phi_n}{\Phi_n^*\Phi_n^\circ}\\
=e^{n\big(\bar Q+\BO(P_{i_1}P_{j_1}P_{i_2}P_{j_2})\big)}\bigg(1+\BO\Big(\prod_{i\in J}P_i\Big)\bigg),$$
$$\Cov(X_{i_1,j_1}^{(n)},X_{i_2,j_2}^{(n)})=
\Big[(e^{n\bar Q}-1)+\BO\Big(nP_{i_1}P_{j_1}P_{i_2}P_{j_2}+\prod_{i\in J}P_i\Big)e^{n\bar Q}\Big]\mathbb{P}(X_{i_1,j_1}^{(n)}=0)\mathbb{P}(X_{i_2,j_2}^{(n)}=0),$$
valid for $nP_{i_1}P_{j_1}P_{i_2}P_{j_2}=\BO(1)$.
It will turn out that in some of the cases we have $\bar Q=0$.
In cases where $\bar Q>0$ we always have $\bar Q=\BO\big(\prod_{i\in J}P_i\big)$ and
$\bar Q\le\frac{1-\delta}2\eps$, with $\eps:=P_{i_1}P_{j_1}+P_{i_2}P_{j_2}$.
Using the latter, and \eqref{exp3}, as well as $e^{n\bar Q}-1\le n\bar Qe^{n\bar Q}$, we obtain
$$e^{n\bar Q}\mathbb{P}(X_{i_1,j_1}^{(n)}\!=\!0)\mathbb{P}(X_{i_2,j_2}^{(n)}\!=\!0)
=\BO\Big(e^{-\frac{\delta}2n\eps}\Big),$$
$$(e^{n\bar Q}-1)\lb\mathbb{P}(X_{i_1,j_1}^{(n)}\!=\!0)\mathbb{P}(X_{i_2,j_2}^{(n)}\!=\!0)-e^{-n\eps}\rb
=
\BO\Big(n\bar Q\sqrt{n}\eps e^{-\frac{\delta}2 n\eps}\Big)
=\BO\Big(\sqrt{n}\bar Qe^{-\frac{\delta}4 n\eps}\Big).
$$
In case of $nP_{i_1}P_{j_1}P_{i_2}P_{j_2}=\Omega(1)$ we use
\eqref{exp3} to obtain
$$\Cov(X_{i_1,j_1}^{(n)},X_{i_2,j_2}^{(n)})= \BO\Big(e^{-\frac{n}2\eps}\Big)
=\BO\Big(nP_{i_1}P_{j_1}P_{i_2}P_{j_2}e^{-\frac{\delta}4n\eps}\Big),$$
and all this results in
\small{$$
\Cov(X_{i_1,j_1}^{(n)},X_{i_2,j_2}^{(n)})
=(e^{n\bar Q}-1)e^{-n(P_{i_1}\!P_{j_1}+P_{i_2}\!P_{j_2})}
+\BO\Big((nP_{i_1}P_{j_1}P_{i_2}P_{j_2}+\sqrt{n}\prod_{i\in J}P_i)e^{-\frac{\delta}4 n(P_{i_1}\!P_{j_1}+P_{i_2}\!P_{j_2})}\Big).
$$}\normalsize
We distinguish the following cases, only Cases 1, 5 and 6 involving $\bar Q\ne0$, and Case 6 slightly deviating from the general pattern outlined above.\\[.2cm]
\textbf{Case 1:} Pairs $(i,r),(r,t)$ with $i,r,t$ all different.

The matrix $\bar\Pi$ and its characteristic polynomial $p$ are given by
$$
\bar\Pi=\begin{bmatrix}P_e&P_i&P_r&P_t\\
P_e&P_i&0&P_t\\
P_e&P_i&P_r&0\\
P_e&P_i&P_r&P_t
\end{bmatrix},\qquad
p(\la)=\la^4-\la^3+P_r(P_i+P_t)\la^2-P_iP_rP_t\la.
$$
Using Algorithm~1 and \eqref{C1gen}, we obtain
\begin{align*}
\la_1&=1-P_iP_r-P_rP_t+P_iP_rP_t+\BO_4^*,\\
C_1&=1+P_iP_r+P_rP_t-2P_iP_rP_t+\BO_4^*.
\end{align*}
We can see that $\la_1=1-P_iP_r-P_rP_t+P_rP_iP_t+P_r^2\BO_2^*$ holds, by noting that $\la_1$ is a $C^\infty$ function of the coefficients
$P_r(P_i+P_t)$ and $-P_iP_rP_t$ of the polynomial $p$, and terms of order 2 or higher contribute $P_r^2\BO_2^*$.
Thus, by the MVT, for some $0<p_i<P_i,0<p_t<P_t$, $$\ln\frac{\la_1}{\la_*\la_\circ}=P_iP_t\dfrac{\partial^2\ln \la_1}{\partial P_i\partial P_t}(p_i,p_t)=P_iP_tP_r(1+\BO(P_r)).$$
So \eqref{barQerror} is established with $\bar Q=P_iP_rP_t$, which indeed
satisfies $\bar Q\le \frac14P_r(P_i+P_t)\le\frac{1-\delta}2\eps$,
since $\delta\le\frac12$.
\goodbreak
\ignore{Assume now $nP_iP_r=\BO(n^{1/3})$ and $nP_rP_t=\BO(n^{1/3})$, ensuring validity of
\eqref{Pir} for pairs $(i,r)$ and $(r,t)$.
In terms of the general discussion above we have here
$Q_*=Q_\circ=0$, $Q=P_iP_rP_t\le \frac14P_r(P_i+P_t)$, therefore
$e^{nQ}\mathbb{P}(X_{i,r}^{(n)}=0)\mathbb{P}(X_{r,t}^{(n)}=0)=\BO\Big(
e^{\frac n4P_r(P_i+P_t)}e^{-n(P_iP_r+P_rP_t)}\Big)=\BO\Big(e^{-\frac {3n}4(P_iP_r+P_rP_t)}\Big)$, using \eqref{Pir}, leading to
$$\Cov(X_{i,r}^{(n)},X_{r,t}^{(n)})=(e^{nP_iP_rP_t}-1)\mathbb{P}(X_{i,r}^{(n)}=0)\mathbb{P}(X_{r,t}^{(n)}=0)+ \BO(P_iP_rP_t+nP_iP_r^2P_t)e^{-\frac {3n}4(P_iP_r+P_rP_t)}.$$
Using $e^x-1\le xe^x$ for $x\ge0$, we obtain 
\begin{align*}(e^{nP_iP_rP_t}-1)&\Big[\mathbb{P}(X_{i,r}^{(n)}=0)\mathbb{P}(X_{r,t}^{(n)}=0)-e^{-n(P_iP_r+P_rP_t)}\Big]\\
&=\BO\lp nP_iP_rP_t\,e^{nP_iP_rP_t}
e^{-n(P_iP_r+P_rP_t)}(P_iP_r+P_rP_t+nP^2_iP^2_r+nP^2_rP^2_t)\rp\\
&=\BO\lp P_iP_rP_t\,e^{-\frac{n}{2}(P_iP_r+P_rP_t)}\rp, \end{align*}
using $P_iP_rP_t\le \frac14P_r(P_i+P_t)$ again, as well as $(nP_iP_r)^\alpha e^{-\beta nP_iP_r}=\BO(1)$ for fixed $\alpha,\beta>0$ (we used $\beta=\frac14$ and $\alpha\in\{1,2\}$).
By the latter we also derive
$nP_iP_r^2P_t\,e^{-\frac{n}{4}(P_iP_r+P_rP_t)}=\BO\lp n^{-1/2}P_i^{1/4}P_r^{1/2}P_t^{1/4}\rp$ and
$P_iP_rP_t\,e^{-\frac{n}{4}(P_iP_r+P_rP_t)}=\BO\lp n^{-1/2}P_i^{\frac34}P_r^{1/2}P_t^{\frac34}\rp=\BO\lp n^{-1/2}P_i^{1/4}P_r^{1/2}P_t^{1/4}\rp$, finally leading to
$$\Cov(X_{i,r}^{(n)},X_{r,t}^{(n)})=
\big(e^{nP_iP_rP_t}-1\big)e^{-nP_iP_r-nP_rP_t}
+\BO\lp \frac1n\,(nP_iP_r)^{1/4}\,(nP_rP_t)^{1/4}\,e^{-\frac{n}{4}(P_iP_r+P_rP_t)}\rp.$$}

\textbf{Case 2a:} Pairs $(i,r),(i,t)$ with $i,r,t$ all different.

The matrix $\bar\Pi$, its characteristic polynomial $p$, and asymptotics of $\la_1$ and $C_1$ are given by
$$
\bar\Pi=\begin{bmatrix}P_e&P_i&P_r&P_t\\
P_e&P_i&0&0\\
P_e&P_i&P_r&P_t\\
P_e&P_i&P_r&P_t
\end{bmatrix},\qquad
\begin{aligned}
p(\la)&=\la^4-\la^3+P_i(P_r+P_t)\la^2,\\
\la_1&=1-P_iP_r-P_iP_t+\BO_4^*,\\
C_1&=1+P_iP_r+P_iP_t+\BO_4^*.
\end{aligned}
$$
Again, $\la_1$ is a $C^\infty$ function of the coefficient
$P_i(P_r+P_t)$, leading to $\la_1=1-P_iP_r-P_iP_t+P_i^2\BO_2^*$,
which we use to derive
$\ln(\frac{\la_1}{\la_*\la_\circ})=P_rP_t\tfrac{\partial^2\ln \la_1}{\partial P_r\partial P_t}(p_r,p_t)=\BO(P_rP_tP_i^2)$, yielding \eqref{barQerror} with $\bar Q=0$. \\[.2cm]
\textbf{Case 2b:} Pairs $(r,i),(t,i)$ with $i,r,t$ all different.

Here the matrix (call it $\bar\Pi_b$) can be seen to be a similarity transformation involving diagonal matrices of the transposed matrix (call it $\bar\Pi_a$) in Case 2a, more precisely, with $\bp:=\pi\bar\Pi=[P_e,(P_i)_{i\in I}]$, we have
$\bar\Pi_b=\textup{Diag}(\bp)^{-1}\bar\Pi_a^t\textup{Diag}(\bp)$, leading
to $\bp\bar\Pi_b^{n-1}\mathds{1}=\bp\bar\Pi_a^{n-1}\mathds{1}$, and
implying that $p(\lambda),\la_1,C_1$, and also the covariance, are the same as in Case 2a. \\[.2cm]
\textbf{Case 3:} Pairs $(i,i),(r,t)$ with $i,r,t$ all different.

The matrix $\bar\Pi$, its characteristic polynomial $p$, and asymptotics of $\la_1$ and $C_1$ are given by
$$
\bar\Pi=\begin{bmatrix}P_e&P_i&P_r&P_t\\
P_e&0&P_r&P_t\\
P_e&P_i&P_r&0\\
P_e&P_i&P_r&P_t
\end{bmatrix},\qquad
\begin{aligned}
p(\la)&=\la^4-(1-P_i)\la^3-(P_i-P_i^2-P_rP_t)\la^2+P_iP_rP_t\la,\\
\la_1&=1-P_i^2-P_rP_t+P_i^3+\BO_4^*,\\
C_1&=1+P_i^2+P_rP_t-2P_i^3+\BO_4^*.
\end{aligned}
$$
Denoting by $\la_\circ=\lim_{P_i\to0}\la_1$ the largest zero of $\la^2-\la+P_rP_t$, and $r(\la)=\frac{p(\la)}\la$, we compute
$$r(\la_\circ+P_i^2\mu)=
P_i^2\la_\circ+P_i^2(P_i^2+2P_i\la_\circ+2{\la_\circ}^2-P_i-\la_\circ)\mu+P_i^4(P_i+3\la_\circ-1)\mu^2+P_i^6\mu^3=0,$$
and conclude by the implicit function theorem, using $\la_\circ=1+\BO_2^*$, that there is a unique $C^\infty$
function $\mu$ of $P_i,P_r,P_t$ near the origin, satisfying $\mu(0,0,0)=-1$, such that
$\la_1=\la_\circ+P_i^2\mu$. This leads to
$\frac{\partial^2 \la_1}{\partial P_i\partial P_r}=\BO(P_i)$, and similarly $\frac{\partial^2 \la_1}{\partial P_i\partial P_t}=\BO(P_i)$,
resulting in $\ln(\frac{\la_1}{\la_*\la_\circ})=\BO(P_rP_tP_i^2)$, yielding \eqref{barQerror}.\\[.2cm]
\textbf{Case 4:} Pairs $(i,j),(r,t)$ with $i,j,r,t$ all different.

The matrix $\bar\Pi$, its characteristic polynomial $p$, and asymptotics of $\la_1$ and $C_1$ are given by
$$
\bar\Pi=\begin{bmatrix}P_e&P_i&P_j&P_r&P_t\\
P_e&P_i&0&P_r&P_t\\
P_e&P_i&P_j&P_r&P_t\\
P_e&P_i&P_j&P_r&0\\
P_e&P_i&P_j&P_r&P_t
\end{bmatrix},\qquad
\begin{aligned}
p(\la)&=\la^5-\la^4+(P_iP_j+P_rP_t)\la^3,\\
\la_1&=1-P_iP_j-P_rP_t+\BO_4^*,\\
C_1&=1+P_iP_j+P_rP_t+\BO_4^*.
\end{aligned}
$$
Observe that
$\frac{\partial^2\ln \la_1}{\partial P_i\partial P_r}=\mathcal{O}_2^*$ and
$\frac{\partial^2\ln \la_1}{\partial P_j\partial P_t}=\mathcal{O}_2^*$
lead to $\ln(\frac{\la_1}{\la_*\la_\circ})=\BO(P_iP_jP_rP_t)$, yielding \eqref{barQerror}.\\[.2cm]
\goodbreak
\textbf{Case 5:} Pairs $(i,r),(r,i)$ with $i,r$ different.

The matrix $\bar\Pi$, its characteristic polynomial $p$, and asymptotics of $\la_1$ and $C_1$ are given by
$$
\bar\Pi=\begin{bmatrix}P_e&P_i&P_r\\
P_e&P_i&0\\
P_e&0&P_r\\
\end{bmatrix},\qquad
\begin{aligned}
p(\la)&=\la^3-\la^2+P_iP_r\la+P_eP_iP_r,\\
\la_1&=1-2P_iP_r+P_i^2P_r+P_iP_r^2+\BO_4^*,\\
C_1&=1+2P_iP_r-2P_i^2P_r-2P_iP_r^2+\BO_4^*.
\end{aligned}
$$
Note that $\la_1$ is a $C^\infty$ function of the coefficients
$P_iP_r$ and $P_iP_r(1-P_i-P_r)$, leading to $$\la_1=1-2P_iP_r+P_iP_r(P_i+P_r)+\BO(P_i^2P_r^2),$$
which, together with $\la_*=\la_\circ=1-P_iP_r+\BO(P_i^2P_r^2)$, we use to derive
$$\frac{\la_1}{\la_*\la_\circ}=1+P_i^2P_r+P_iP_r^2+\BO(P_i^2P_r^2).$$
This is in accordance with \eqref{barQerror}, with $\bar Q=P_i^2P_r+P_iP_r^2=P_iP_r(P_i+P_r)\le P_iP_r(1-\delta)
=\frac{1-\delta}2\eps$.

\ignore{In terms of the general discussion above we have here
$Q_*=Q_\circ=0$, $Q=P^2_iP_r+P_iP_r^2\le P_iP_r$, therefore
$e^{nQ}\mathbb{P}(X_{i,r}^{(n)}=0)\mathbb{P}(X_{r,i}^{(n)}=0)=\BO\lp
e^{nP_iP_r}e^{-2nP_iP_r}\rp=\BO\lp e^{-nP_iP_r}\rp$, using \eqref{Pir}, leading to
$$\Cov(X_{i,r}^{(n)},X_{r,i}^{(n)})=(e^{nP_iP_r(P_i+P_r)}-1)\mathbb{P}(X_{i,r}^{(n)}=0)\mathbb{P}(X_{r,i}^{(n)}=0)+ \BO(P_iP_r+nP_i^2P_r^2)e^{-nP_iP_r}.$$
Using $e^x-1\le xe^x$ for $x\ge0$, and assuming $P_iP_r=\BO(\frac{\ln n}n)$, we obtain by \eqref{Pir}
\begin{align*}(e^{nP_iP_r(P_i+P_r)}-1)&\Big[\mathbb{P}(X_{i,r}^{(n)}=0)\mathbb{P}(X_{r,i}^{(n)}=0)-e^{-2nP_iP_r}\Big]\\
&=\BO\lp nP_iP_r(P_i+P_r)\,e^{nP_iP_r(P_i+P_r)}
e^{-2nP_iP_r}(P_iP_r+nP^2_iP^2_r)\rp\\
&=\BO\lp \ln^2 (n)\,P_iP_r\,e^{-nP_iP_r}\rp, \end{align*}
and finally
$$\Cov(X_{i,r}^{(n)},X_{r,i}^{(n)})=
\big(e^{nP_iP_r(P_i+P_r)}-1\big)e^{-2nP_iP_r}+\BO\lp \ln^2 (n)\,P_iP_r\,e^{-nP_iP_r}\rp.$$}
\textbf{Case 6a:} Pairs $(i,i),(i,r)$ with $i,r$ different.

The matrix $\bar\Pi$, its characteristic polynomial $p$, and asymptotics of $\la_1$ and $C_1$ are given by
$$
\bar\Pi=\begin{bmatrix}P_e&P_i&P_r\\
P_e&0&0\\
P_e&P_i&P_r\\
\end{bmatrix},\qquad
\begin{aligned}
p(\la)&=\la^3-(1-P_i)\la^2-P_iP_e\la,\\
\la_1&=1-P_i^2-P_iP_r+P_i^2P_r+P_i^3+\BO_4^*,\\
C_1&=1+P_i^2+P_iP_r-2P_i^2P_r-2P_i^3+\BO_4^*.
\end{aligned}
$$
We start deriving the more precise estimate $\la_1=1-P_i^2-P_iP_r+P_i^2P_r+P_i^3+P_i^2\BO_2^*$:

Abbreviating $\sigma=P_i+P_r$, $\kappa=P_i-P_i^2$, we use $p(\la_1)=0$
to infer the existence of a function $\mu$ that satisfies
$\la_1=1-\kappa\sigma+P_i^2\mu$. Indeed, from
\begin{align*}
0&=\la_1^2-(1-P_i)\la_1-P_i(1-\sigma)\\
&=(1-\kappa\sigma+P_i^2\mu)^2-(1-P_i)(1-\kappa\sigma+P_i^2\mu)-P_i(1-\sigma)\\
&=\kappa^2\sigma^2+\sigma(P_i-\kappa-P_i\kappa)+P_i^2\mu(1+P_i-2\kappa\sigma)+P_i^4\mu^2\\
&=P_i^2\left[(1-P_i)^2\sigma^2+\sigma P_i+(1+P_i-2\kappa\sigma)\mu+P_i^2\mu^2\right]
\end{align*}
we conclude by the implicit function theorem that there is a unique $C^\infty$
function $\mu$ of $P_i,P_r$ near the origin, satisfying $\mu=\BO_2^*$.

Since $\lim_{P_r\to0}\la_1=\la_*$ and $\lim_{P_r\to0}\la_\circ=1$,
we have $\frac{\la_1}{\la_*\la_\circ}=1+\BO(P_r)$.
This estimate will now be refined. From
$\la_*=1-P_i^2+P_i^3+\BO(P_i^4)$ and $\la_\circ=1-P_iP_r+\BO(P_i^2P_r^2)$ we deduce $\la_*\la_\circ=1-P_i^2-P_iP_r+P_i^3+P_i^2\BO_2^*$ and
$$\frac{\la_1}{\la_*\la_\circ}=\frac{\la_*\la_\circ+P_i^2P_r+P_i^2\BO_2^*}{\la_*\la_\circ}=1+P_i^2P_r+P_i^2\BO_2^*=1+P_i^2P_r+P_i^2P_r\BO_1^*=1+P_i^2P_r+\BO(P_i^2P_r).$$
This is not quite \eqref{barQerror}, but
$\bar Q=P_i^2P_r=\frac{P_r}2P_i^2+\frac{P_i}2P_iP_r\le\frac{1-\delta}2\eps$ is satisfied, and $\BO(P_i^2P_r)$ turns out to be a sufficiently good substitute for $\BO(P_{i_1}P_{j_1}P_{i_2}P_{j_2})$.
\\[.2cm]
\textbf{Case 6b:} Pairs $(i,i),(r,i)$ with $i,r$ different.

Here the matrix $\bar\Pi$ can be seen to be a similarity transformation of the transposed matrix in Case 6a, implying that $p(\lambda),\la_1,C_1$, and also the covariance, are the same as in Case 6a.

We summarize the covariances $\Cov(X_{i_1,j_1}^{(n)},X_{i_2,j_2}^{(n)})$,
asymptotics valid for $(P_j)_{j\in J}\in\mathcal{D}^J_\delta$,
\begin{align*}
\Cov(X_{i,r}^{(n)},X_{r,t}^{(n)})&=
\big(e^{nP_iP_rP_t}\!-\!1\big)e^{-n(P_iP_r+P_rP_t)}+\BO\big(nP_iP_r^2P_t\!+\!\sqrt{n}P_iP_rP_t\big) e^{-\frac {\delta}4n(P_iP_r+P_rP_t)}
\tag{Case 1}\\
\Cov(X_{i,r}^{(n)},X_{i,t}^{(n)})&=\Cov(X_{r,i}^{(n)},X_{t,i}^{(n)})=
\BO\big(nP_i^2P_rP_t+P_iP_rP_t\big) e^{-\frac {\delta}4n(P_iP_r+P_iP_t)}\tag{Cases 2}\\
\Cov(X_{i,i}^{(n)},X_{r,t}^{(n)})&=
\BO\big(nP_i^2P_rP_t+P_iP_rP_t\big) e^{-\frac {\delta}4n(P_i^2+P_rP_t)}\tag{Case 3}\\
\Cov(X_{i,j}^{(n)},X_{r,t}^{(n)})&=
\BO\big(nP_iP_jP_rP_t\big) e^{-\frac {\delta}4n(P_iP_j+P_rP_t)}\tag{Case 4}\\
\Cov(X_{i,r}^{(n)},X_{r,i}^{(n)})&=
\big(e^{nP_iP_r(P_i+P_r)}\!-\!1\big)e^{-2nP_iP_r}+\BO\big(nP_i^2P_r^2\!+\!\sqrt{n}P_iP_r\big) e^{-\frac {\delta}2nP_iP_r}\tag{Case 5}\\
\Cov(X_{i,i}^{(n)},X_{i,r}^{(n)})&=\Cov(X_{i,i}^{(n)},X_{r,i}^{(n)})=
\BO\big(nP_i^2P_r+P_iP_r\big)e^{-\frac {\delta}4n(P_i^2+P_iP_r)}\tag{Cases 6}
\end{align*}
We continue showing that the multiple sums of error terms arising in \eqref{Vir} and Cases 1--6 are negligible. In addition to  \eqref{iialpha} we will
also use that
\beq\label{iralpha}\sum_{i,k\ge1}(nP_iP_k)^\alpha
e^{-nP_iP_k}=\mathcal{O}(\ln n)\ \mbox{ holds for }\alpha>0.\eeq
This can be deduced from \eqref{iialpha}, using
$\beta=1$, observing
$$\sum_{i,k\ge1}(nP_iP_k)^\alpha
e^{-nP_iP_k}=\sum_{\ell\ge2}(\ell-1)(n\tfrac pqP_\ell)^\alpha
e^{-n\tfrac pqP_\ell},$$
and furthermore
$$\sum_{\ell\ge2}\ell(nP_\ell)^\alpha
e^{-nP_\ell}=\sum_{\ell\ge2}\BO(\ln n-\ln(nP_\ell))(nP_\ell)^\alpha e^{-nP_\ell}=\BO(\ln n),$$
because of $x^\alpha\ln x=\BO(x^{\alpha/2})$.
Note that \eqref{iralpha} yields
$\sum_{i,r\ge1}\!P_iP_r e^{-\frac {2-\delta}4 nP_iP_r}\!=\!\BO\big(\frac{\ln n}{n}\big)$, which settles \eqref{Vir}, and also Case~5, where the double sum is $\BO\Big(\frac{\ln n}{\sqrt{n}}\Big)$,
and Case~4, with quadruple sum of order $\BO\Big(\frac{\ln^2 n}{n}\Big)$.
Using $P_iP_t\le\sqrt{P_iP_t}$, Case 1 can be reduced to bounding the sum
$$\frac1{\sqrt{n}}\sum_{i,r,t\ge1}\sqrt{nP_iP_r}\sqrt{nP_rP_t}e^{-\frac {\delta}4n(P_iP_r+P_rP_t)}=
\BO\bigg(\frac1{\sqrt{n}}\sum_{i,r\ge1}\sqrt{nP_iP_r}\,e^{-\frac {\delta}4nP_iP_r}\bigg)
=\BO\left(\frac{\ln n}{\sqrt{n}}\right),
$$
where for the inner sum (w.r.t.~$t$) we used \eqref{iialpha}. Similarly Cases 2
give rise to triple sums of order $\BO\left(\frac{\ln n}{n}\right)$. The same is true for Case~3, which is seen by upper bounding the triple sums by
\small{$$\frac1n\sum_{i,r,t\ge1}(nP_i^2)^\alpha(nP_rP_t)^\alpha e^{-\frac {\delta}4n(P_i^2+P_rP_t)}=
\frac1n\sum_{i\ge1}(nP_i^2)^{\alpha} e^{-\frac {\delta}4nP_i^2}\sum_{r,t\ge1}(nP_rP_t)^{\alpha}e^{-\frac {\delta}4nP_rP_t}=
\BO\left(\frac{\ln n}n\right),$$}\normalsize
where $\alpha\in\{1/2,1\}$. Finally, the following estimates
\beq\label{case61}\sum_{i,r\ge1}nP_i^2P_re^{-\frac {\delta}4n(P_i^2+P_iP_r)}=
\frac1{\sqrt{n}}\sum_{i\ge1}(nP_i^2)^{1/2}e^{-\frac {\delta}4nP_i^2}\sum_{r\ge1}nP_iP_re^{-\frac {\delta}4nP_iP_r}=\BO\left(\frac{\ln n}{\sqrt{n}}\right),\eeq
$$\sum_{i,r\ge1}P_iP_re^{-\frac {\delta}4n(P_i^2+P_iP_r)}\le
\frac1{n^{\frac34}}\sum_{i\ge1}(nP_i^2)^{1/4}e^{-\frac {\delta}4nP_i^2}\sum_{r\ge1}(nP_iP_r)^{1/2}e^{-\frac {\delta}4nP_iP_r}=\BO\left(\frac{\ln n}{n^{\frac34}}\right),$$
deal with Cases~6. The total contribution of error terms is therefore of order $\BO\left(\frac{\ln n}{\sqrt{n}}\right)$.

\ignore{Next we bound the multiple sums of the main terms in \eqref{Vir} and Cases 1 and 5 over the range of indices where $\max(P_{i_1}P_{j_1},P_{i_2}P_{j_2})=\BO(\frac{\ln n}n)$ does not hold. We will actually sum over larger sets described by inequalities
$P_{i_1}P_{j_1}+P_{i_2}P_{j_2}\ge\mbox{const\,}\frac{\ln n}n$:
$$\sum_{P_iP_r\ge 2\frac{\ln n}n}\lb e^{-nP_iP_r}-e^{-2nP_iP_r}\rb
\le \sum_{P_iP_r\ge 2\frac{\ln n}n}nP_iP_r e^{-nP_iP_r}
\le \frac1n\sum_{i,r\ge1}nP_iP_r e^{-\frac n2P_iP_r}
=\BO\left(\frac{\ln n}{n}\right),$$
$$\sum_{P_iP_r+P_rP_t\ge 2\frac{\ln n}{n}}\underbrace{\lb e^{nP_iP_rP_t}-1\rb e^{-nP_iP_r-nP_rP_t}}_{\le nP_iP_rP_t\, e^{-\frac{3n}4P_iP_r-\frac{3n}4P_rP_t}}
\le \frac1n\sum_{i,r,t\ge1}nP_iP_rP_t\, e^{-\frac{n}4P_iP_r-\frac{n}4P_rP_t}
=\BO\left(\frac{\ln n}{n}\right),$$
$$\sum_{P_iP_r\ge 2\frac{\ln n}n}\lb e^{nP_iP_r(P_i+P_r)}-1\rb e^{-2nP_iP_r}
\le \sum_{P_iP_r\ge 2\frac{\ln n}n}nP_iP_r\, e^{-nP_iP_r}
\le\frac1n\sum_{i,r\ge1}nP_iP_r\, e^{-\frac n2P_iP_r}
=\BO\left(\frac{\ln n}{n}\right),$$
where we used previous estimates from Cases 1 and 5.
Next we estimate the multiple sums of the upper bounds given in \eqref{Virle} and \eqref{Covle} over ranges of indices like
those above:
$$\sum_{P_iP_r\ge 6\frac{\ln n}n}\lb e^{-\frac n3nP_iP_r}-e^{-2\frac n3P_iP_r}\rb=\BO\left(\frac{\ln n}{n}\right),\qquad
\sum_{i,j,r,t\ge1} nP_iP_jP_rP_t\,e^{-\frac n3(P_iP_j+P_rP_t)}=\BO\left(\frac{\ln^2 n}{n}\right),$$
$$\sum_{P_iP_r+P_rP_t\ge 6\frac{\ln n}{n}}\!\!\!\!\!\! nP_iP_rP_t\,e^{-\frac n3(P_iP_r+P_rP_t)}
\le\frac1n\sum_{i,r,t\ge1}nP_iP_rP_t\,e^{-\frac n6(P_iP_r+P_rP_t)}
=\BO\left(\frac{\ln n}{n}\right),$$
$$\sum_{P_i^2+P_rP_t\ge 6\frac{\ln n}{n}}\!\!\!\!\!\!\!\!\!\!\!\!\! nP_iP_rP_t\,e^{-\frac n3(P_i^2+P_rP_t)}=\BO\left(\frac{\ln n}{n}\right),\qquad
\sum_{P_i^2+P_iP_r\ge 3\frac{\ln n}{n}}\!\!\!\!\!\!\!\! nP_i^2P_r\,e^{-\frac n3(P_i^2+P_iP_r)}=\BO\left(\frac{\ln n}n\right),$$
so the error term in \eqref{varX2} is confirmed.}
We are left with dealing with the sums of the main terms of Cases~1 and 5, and \eqref{Vir}. Note that Case~1 has a twin case,
$\Cov(X_{i,r}^{(n)},X_{r,t}^{(n)})=\Cov(X_{r,i}^{(n)},X_{t,r}^{(n)})$.

Denote $H(i,j,k)=(e^{nP_iP_jP_k}-1)e^{-nP_iP_j-nP_jP_k}$ and
$H^\circ(i,j)=(e^{nP_iP_j(P_i+P_j)}-1)e^{-2nP_iP_j}$.
Observe that
\begin{align*}\sum_{i\ne j}(e^{nP^2_iP_j}-1)(e^{nP_iP_j^2}-1)e^{-2nP_iP_j}
&\le\sum_{i,j}nP_i^2P_je^{nP^2_iP_j}nP_iP_j^2e^{nP_iP_j^2}e^{-2nP_iP_j}\\
&\le\sum_{i,j}n^2P_i^3P_j^3e^{-2\delta nP_iP_j}=\BO\lp\frac{\ln n}{n}\rp\end{align*}
and $(e^{a+b}-1)=(e^{a}-1)+(e^{b}-1)+(e^{a}-1)(e^{b}-1)$ imply
$\sum_{i\ne j}H^\circ(i,j)=2\sum_{i\ne j}H(i,j,i)+\BO\lp\frac{\ln n}{n}\rp$.
Therefore we have
\begin{align*}2\!\!\!\!\!\!\sum_{{i,j,k\ge1}\atop{|\{i,j,k\}|=3}}\!\!\!H(i,j,k)
+\sum_{{i,j\ge1}\atop{|\{i,j\}|=2}}H^\circ(i,j)
&\sim2\!\!\!\sum_{{i,j,k\ge1}\atop{j\not\in\{i,k\}}}\!\!\!H(i,j,k)
\\[-.75cm]
&=2\!\!\sum_{i,j,k\ge1}\!\!\!H(i,j,k)-4\overbrace{\sum_{i,j\ge1}\!\!H(i,i,j)}^{\BO\left(\frac{\ln n}{\sqrt{n}}\right)}+2\overbrace{\sum_{i\ge1}H(i,i,i)}^{\BO\left(\frac{1}{\sqrt{n}}\right)},\end{align*}
where we have estimated two of the sums using \eqref{iialpha} and \eqref{case61}. Asymptotics of the sum $\sum_{i,j,k\ge1}H(i,j,k)$
are computed in Appendix~\ref{A3},
confirming $T^{(n)}_2$ as given in \eqref{T2}. The sum
$$\sum_{i,r\ge1}\lb e^{-nP_iP_r}-e^{-2nP_iP_r}\rb
=\sum_{i,r\ge1}\lb 1-e^{-2nP_iP_r}\rb-\sum_{i,r\ge1}\lb 1-e^{-nP_iP_r}\rb=\tilde G(2np^2)-\tilde G(np^2)
,$$
which, as we have seen, is an asymptotic equivalent of $\sum_{i,r\ge1}\Var X_{i,r}^{(n)}$, is evaluated in Appendix~\ref{A2}, confirming $S^{(n)}_2$ as given in \eqref{S2}.
This completes the proof of the lemma, and also proves \eqref{var3}, as we have seen, that multiple sums of covariances
$\Cov(X_{i_1,j_1}^{(n)},X_{i_2,j_2}^{(n)})$ with $i_1=j_1$, but $i_2\ne j_2$, are negligible.
\epr
\begin{rem}
Along the lines of the two preceding proofs an
 independent proof of Theorem~\ref{EX1EX3} could easily be furnished. We would use \eqref{exp3}, \eqref{Pii}, \eqref{Pir} to identify
$\sum_{i\ge1}(1-e^{-nP_i^2})$ and $\sum_{i\ne j}(1-e^{-nP_iP_j})$
as asymptotic equivalents of $\E X^{(n)}_1$ and $\E X^{(n)}_3$, leading to
$\E X^{(n)}_1\sim G(np^2)$ and $\E X^{(n)}_3\sim \tilde G(np^2)-G(np^2)$, with
$G, \tilde G$ from Appendices \ref{A1} and \ref{A2}.
\end{rem}

\subsection {More than two pairs of identical letters}
We now turn to the case of $k$ pairs $(i_1,i_1),\ldots,(i_k,i_k)$, allowing for $k>2$.
\begin{lem} 
Fix a set $I:=\{i_1,i_2,\ldots,i_k\}$ of size $k$, assuming $i_k<\ldots<i_1$, and thus $P_{i_1}<\ldots<P_{i_k}$.
Let $\eps:=\sum_{i\in I}P_i^2$.
Then we have
\beq\label{Pbound}\mathbb{P}(X_{i,i}^{(n)}=0,i\in I)=\sum_{j=1}^{k+1}C_j\lambda_j^n=\begin{cases}\BO(\frac1n),&\textup{for }\eps\ge\frac{3\ln n}n,\\
C_1\la_1^n+\BO(\frac1n),&\textup{for }\eps\le 1/4, \end{cases}\eeq
with all $\lambda_i$ different, and error terms holding uniformly in $k$. More precisely, we have $\la_1>|\la_j|>0$ for $2\le j\le k+1$, and $-P_{i_k}<\la_{k+1}<
-P_{i_{k-1}}<\la_{k}<\ldots<-P_{i_1}<\la_2<0$. Moreover,
\begin{align}\label{la1k}\la_1&=1-\sum_{i\in I}P_i^2+\sum_{i\in I}P_i^3+\BO(\eps^2),\\
\label{C1k}C_1&=1+\BO(\eps),\end{align}
again with error terms holding uniformly in $k$.
\end{lem}
\bpr As before, we let $e:=\mathbb{N}\setminus I$ and
$P_e:=1-\sum_{i\in I} P_i$, and introduce the matrix
$$\bar\Pi=\begin{bmatrix}P_e&P_{i_1}&P_{i_2}&\cdots&P_{i_k}\\
P_e&0&P_{i_2}&\cdots&P_{i_k}\\
P_e&P_{i_1}&0&\cdots&P_{i_k}\\
\vdots&\vdots&\vdots&\ddots&\vdots\\
P_e&P_{i_1}&P_{i_2}&\cdots&0
\end{bmatrix}.$$
In order to find eigenvalues and corresponding left and right eigenvectors of $\bar\Pi$,
we have to solve the following systems,
\beq\label{evk}
\begin{aligned}[c]
\lambda&=P_e\Big(1+\sum_{j\in I} \beta_j\Big)\\
\lambda\beta_i&=P_i\Big(1+\sum_{j\in I\setminus\{i\}} \beta_j\Big),\quad i\in I\\
\end{aligned}
\qquad\qquad
\begin{aligned}[c]
\lambda&=P_e\Big(1+\sum_{j\in I} P_j\mu_j\Big)\\
\lambda\mu_i&=1+\sum_{j\in I\setminus\{i\}} P_j\mu_j,\quad i\in I\\
\end{aligned}
\eeq
Note that $(\mu_i)_{i\in I}$ solves the right system if and only if
$(\beta_i)_{i\in I}=(P_i\mu_i)_{i\in I}$ solves the left system.
From the left system we easily obtain
\beq\label{betai}\beta_i=\frac{\la P_i}{P_e(\la+P_i)},\textup{ for }i\in I,\eeq
and, upon inserting into the first equation of the left system,
\beq\label{la1impl}\la=P_e+\sum_{i\in I}\frac{\la P_i}{\la +P_i}
=P_e+\sum_{i\in I}P_i-\sum_{i\in I}\frac{P_i^2}{\la +P_i}
=1-\sum_{i\in I}\frac{P_i^2}{\la +P_i}.\eeq
There are at most $k+1$ different solutions to \eqref{la1impl}, those being exactly the eigenvalues of $\bar\Pi$. Defining
$f(\la):=\la-1+\sum_{i\in I}\frac{P_i^2}{\la +P_i}$, we observe the following $k+1$ sign changes on the interval $[-P_{i_k},1]$,
$$\sum\nolimits_{i\in I}\tfrac{P_i^2}{1 +P_i}=f(1)>0>f(0)=-P_e,\quad\lim_{\la\nearrow -P_{i}}f(\la)=-\infty,\
\lim_{\la\searrow -P_{i}}f(\la)=\infty,\mbox{ for }i\in I,$$
from which we obtain the result regarding the locations of the eigenvalues.

We continue with the proof of \eqref{Pbound}. The first estimate, $\BO(\frac1n)$,  directly follows from \eqref{exp3}.
For the second, note that $\eps\le 1/4$ implies $P_{i_k}\le 1/2$. We then use \eqref{barlacheck} and $S$ and $\bw$ as defined in the proof of Lemma~\ref{exp3ClaPhi}. Then for some
orthogonal matrix $Q$ the matrix $\tilde\Pi:=QS\check\Pi S^{-1}Q^{-1}$
is diagonal and satisfies $\rho(\tilde\Pi)=|\la_{k+1}|<P_{i_k}\le 1/2$, and $|\la_{k+1-j}|<P_{i_{k-j}}\le\frac{q^j}2$ for $j\ge1$,
implying
$\|\tilde\Pi^n\|_F\le \frac1{1-q}2^{-n}$. This leads to
$$\pi\check\Pi^n\mathds{1}=\bw(S\check\Pi S^{-1})^{n-1}\bw^t
=\bw Q^{-1}\tilde\Pi^{n-1}Q\bw^t
\le\|\tilde\Pi^{n-1}\|_2\le
\|\tilde\Pi^{n-1}\|_F\le\tfrac2{1-q}2^{-n}=\BO\lp\tfrac1n\rp.$$

Turning now to asymptotic expansions of $\la_1$ and $C_1$, we first provide a convenient representation of the latter in the spirit of
\eqref{la1impl}, starting from \eqref{C1gen},
\beq\label{Ca1impl}C_1=\dfrac{1+\sum_{i\in I}\beta_i}{1+P_e\sum\limits_{i\in I}\frac{\beta_i^2}{P_i}}=\dfrac{\la_1}{P_e+\sum\limits_{i\in I}\frac{\la_1^2P_i}{(\la_1+P_i)^2}}=\dfrac{\la_1}{\la_1-\sum\limits_{i\in I}P_i\lb\frac{\la_1}{\la_1+P_i}-\frac{\la_1^2}{(\la_1+P_i)^2}\rb}
=\dfrac{1}{1-\sum\limits_{i\in I}\frac{P_i^2}{(\la_1+P_i)^2}}.
\eeq
Note that asymptotic estimates of higher order than those given in \eqref{la1k} and \eqref{C1k} could easily be obtained by Algorithm~\ref{alg1},
but as we need error terms uniformly in $k$, we choose another route.
We assume $\eps\le 1/9$ and observe
$f(1-\frac32\eps)=-\frac32\eps+\sum_{i\in I}\frac{P_i^2}{1-\frac32\eps +P_i}
<-\frac32\eps+\sum_{i\in I}\frac{P_i^2}{1-\frac32\eps}
\le-\frac32\eps+\frac65\eps \le0$, which implies $\la_1>1-\frac32\eps$.
Using $\la_1+P_i\le 1+1/3=4/3$ in equation \eqref{la1impl},
we obtain
$$\la_1=1-\sum_{i\in I}\frac{P_i^2}{\la_1 +P_i}\le
1-\tfrac34\sum_{i\in I}P_i^2=1-\tfrac34\eps.$$
Next we employ $1-x\le\frac1{1+x}\le1-x+2x^2$, holding for $x\in[-1/2,1]$, in
\begin{align*}\la_1&\le1-\sum_{i\in I}\frac{P_i^2}{1-\frac34\eps +P_i}\le
1-\sum_{i\in I}P_i^2+\sum_{i\in I}P_i^3-\tfrac34\eps^2,\\
\la_1&\ge1-\sum_{i\in I}\frac{P_i^2}{1-\tfrac32\eps +P_i}\ge
1-\sum_{i\in I}P_i^2+\sum_{i\in I}P_i^3-\tfrac32\eps^2-2\sum_{i\in I}P_i^4+6\eps\sum_{i\in I}P_i^3-\tfrac92\eps^3\\
&\ge 1-\sum_{i\in I}P_i^2+\sum_{i\in I}P_i^3-4\eps^2,
\end{align*}
proving \eqref{la1k}.
Similarly, \eqref{C1k} follows from \eqref{Ca1impl}, using
$\la_1+P_i\ge1-\tfrac32\eps\ge 5/6$:
$$1\le C_1=\Big[{1-\sum_{i\in I}\tfrac{P_i^2}{(\la_1+P_i)^2}}\Big]^{-1}\le\Big[1-\tfrac{36}{25}\sum_{i\in I}P_i^2\Big]^{-1}\le1+2\eps.$$
This completes the proof of the lemma.\epr
\goodbreak
{\bf Proof of Theorem~\ref{thm:asyindep}:}
 We first prove \eqref{eq:asyindep} in the case that $x_i=0$ for all $i\in I$. Letting $\eps:=\sum_{i\in I}P_i^2$ again, by the previous lemma we have
$$
\P(X_{i,i}^{(n)}=0,i\in I)
=C_1\la_1^n+\sum_{j=2}^{k+1}C_j\la_j^n=\prod_{i\in I}e^{-n(P_i^2-P_i^3)}\Big(1+\BO(\eps)+n\BO(\eps^2)\Big)+\BO\Big(\frac1n\Big).
$$
By letting $P_j\to0$ for $j\in I\setminus\{i\}$, we obtain
$$
\P(X_{i,i}^{(n)}=0)=e^{-n(P_i^2-P_i^3)}\Big(1+\BO(P_i^2)+n\BO(P_i^4)\Big)+\BO\Big(\frac1n\Big),
$$
and finally
\beq\label{e:asyindep}
\P(X_{i,i}^{(n)}=0,i\in I)-\prod_{i\in I}\P(X_{i,i}^{(n)}=0)=\left(\prod_{i\in I}e^{-n(P_i^2-P_i^3)}\right)\Big(\BO(\eps)+n\BO(\eps^2)\Big)+\BO\Big(\frac1n\Big)=\BO\Big(\frac1n\Big),
\eeq
using $\prod_{i\in I}e^{-n(P_i^2-P_i^3)}\le e^{-n\eps(1-P_1)}$, and the fact that $e^{-x(1-P_1)}(x+x^2)$ is bounded for $x\ge0$.

Clearly, equation \eqref{eq:asyindep} holds for $I=\{i\}$ and all
$x_i\in\{0,1\}$. Assume that equation \eqref{eq:asyindep} has been shown for all $I$ with $|I|=k$.
Consider $I'$ with $|I'|=k+1$. Then, as we have just shown, equation \eqref{eq:asyindep} holds for $I'$ when $\sum_{i\in I'}x_i=0$. It also holds when $\sum_{i\in I'}x_i=1$: If $x_j=1$, $x_i=0$ for $i\in I'\setminus\{j\}$, then
\begin{align*}
\P(X_{i,i}^{(n)}=x_i,i\in I')&=\P(X_{i,i}^{(n)}=0,i\in I'\setminus\{j\})-
\P(X_{i,i}^{(n)}=0,i\in I'),\\
\prod_{i\in I'}\P(X_{i,i}^{(n)}=x_i)&=\prod_{i\in I'\setminus\{j\}}\P(X_{i,i}^{(n)}=0)-\prod_{i\in I'}\P(X_{i,i}^{(n)}=0),
\end{align*}
so, by taking the difference of these equations, we have
$$\P(X_{i,i}^{(n)}=x_i,i\in I')-\prod_{i\in I'}\P(X_{i,i}^{(n)}=x_i)=\BO\Big(\frac1n\Big).$$
Similarly, by induction on
$\kappa:=\sum_{i\in I'}x_i$, we can prove that \eqref{e:asyindep} holds for all $I'$ with $|I'|=k+1$  and all $x\in\{0,1\}^{k+1}$.
Clearly the error terms $\BO\big(\frac1n\big)$ may now suffer from dependence on $|I|$, but not on $I$, as the values $\{P_i\}_{i\in I}$ did not enter the proof.
\epr
We conclude this subsection with the following conjecture.
\bconj \label  {conj:aipairw}
The same kind of asymptotic independence as in Theorem~\ref{thm:asyindep} holds for $\big(X_{k_i,m_i}^{(n)}\big)_{i\ge1}$, when the sets $\{k_i,m_i\}$ are pairwise disjoint.
\econj

\subsection{Some further results on the probability of avoiding a prescribed set of pairs
}
In this section we aim at a better understanding of $\la_1$ and $C_1$ given in \eqref{ClaPhi}, as examples like
\small{$$\begin{bmatrix}\la_1\\C_1\end{bmatrix}=
\begin{bmatrix}1-P_iP_r-P_i^2P_r^2-2P_i^3P_r^3+\BO_8^*\\
1+P_iP_r+3P_i^2P_r^2+10P_i^3P_r^3+\BO_8^*\end{bmatrix}
,\ \textup{ resp. }\
\begin{bmatrix}\la_1\\C_1\end{bmatrix}=
\begin{bmatrix}1-P_i^2-P_iP_r+P_i^2P_r+P_i^3+\BO_4^*\\
1+P_i^2+P_iP_r-2P_i^2P_r-2P_i^3+\BO_4^*\end{bmatrix}$$}\normalsize
from the proof of Lemma~\ref{Piiir}, resp. from Case~6a in the proof of Lemma~\ref{var2sumeqsumvar2}, suggest that there may be a simple relationship between $\la_1$ and $C_1$.
This turns out to be the case, see \eqref{hipreC} below, and our method of proof also allows for
a representation of the generating function of the probabilities in
\eqref{ClaPhi}. Besides shedding light on above mystery, we hope that the results of this section will turn out
useful when computing asymptotics of higher moments of $X^{(n)}_2$ and
$X^{(n)}_3$, a task however not further pursued in the present paper.

We start with a finite non-empty set of forbidden pairs $\mathcal{I}:=\{(k_i,m_i):i\in I\}$ and
let $J:=\bigcup_{i\in I}\{k_i,m_i\}$.  Using $\bar{\bar\Pi}$ and $\bar\bp$ introduced shortly before Algorithm~1, we define $\tsup[2]{\Pi}:=\mathds{1}\bar\bp-\bar{\bar\Pi}$, i.e.,
$$\tsup[2]{\Pi}_{k,m}:=\begin{cases}P_m,&(k,m)\in \mathcal{I},\\0,&\textup{else,}\end{cases}$$
and use it to define $\psi_1:=\bar\bp\mathds{1}=\sum_{j\in J}P_j=1-P_e$ and
$$\psi_{i+1}:=\bar\bp\tsup[2]{\Pi}{}^i\mathds{1}=
\sum\nolimits_{k_0,\ldots,k_i:(k_0,k_1),\ldots,(k_{i-1},k_i)\in\mathcal{I}}P_{k_0}P_{k_1}\cdots P_{k_i},$$
for $i\ge1$. Note that $\psi_2=\eps$, with $\eps$ introduced in Lemma~\ref{exp3ClaPhi}, and $\psi_3$ is a generalization of $Q$ introduced in \eqref{Qexpr}. Moreover $\psi_i\le(1-P_e)^i=\BO^*_i$ holds for $i\ge1$.
Denote the identity matrix of appropriate dimension by $\Id$ and define a meromorphic function in terms of a resolvent, $$\Psi(z):=\bar\bp\Big(\frac1z\Id+\tsup[2]{\Pi}\Big)^{-1}\mathds{1}=z\bar\bp(\Id+z\tsup[2]{\Pi})^{-1}\mathds{1}=-\sum_{i\ge1}\psi_{i}(-z)^i,$$ with the series converging for $|z|<\frac1{1-P_e}$.
The derivative $\Psi'(z)=\frac1{z^2}\bar\bp\big(\frac1z\Id+\tsup[2]{\Pi}\big)^{-2}\mathds{1}$ will be needed later on.
Denote $$p_{\mathcal{I}}^{(n)}=p_{\mathcal{I}}^{(n)}((P_j)_{j\in J}):=\P(X_{k_i,m_i}^{(n)}=0,i\in I),$$
and for $v\in[0,1]$ consider now the functions $C(v)$ and $\la(v)$ defined via \eqref{ClaPhi} by
$$p_{\mathcal{I}}^{(n)}((vP_j)_{j\in J})\sim C(v)\la(v)^n.$$
Arguing as in the proof of Lemma~\ref{exp3ClaPhi}, i.e., invoking the Perron-Frobenius theorem and the implicit function theorem, these functions are analytic in an open subset of $\mathbb{C}$ containing the interval $[0,1]$.
The following theorem shows how to express $\la(v)$, $C(v)$, and $\mathcal{P}_{\mathcal{I}}(z)\!:=\!\sum_{n\ge0}p_{\mathcal{I}}^{(n)}z^n$,
in terms of $\Psi$.
\begin{thm} The function $\la(v)$ is a solution to the following equation,\goodbreak
\beq\label{hiprela}\la(v)=\frac{1-v\psi_1}{1-\Psi(\frac v{\la(v)})}.\eeq
The function $C$ satisfies
\beq\label{hipreC}C(v)=\la(v)-v\la'(v),\eeq
which, in terms of coefficients, means $[v^n]C(v)=-(n-1)[v^n]\la(v)$.\\
Moreover, the generating function of the sequence $(p_{\mathcal{I}}^{(n)})_{n\ge0}$ satisfies
\beq\mathcal{P}_{\mathcal{I}}(z)\label{hipreP}=
\frac1{1-(1-\psi_1)z-\Psi(z)}.\eeq
\end{thm}
\bpr
We start with \eqref{sysla} -- \eqref{sysmu}, i.e., $\lambda=P_e(1+\beta\mathds{1}),
\beta=\frac1{\lambda}\big[\beta\bar{\bar\Pi}+\bar\bp\big],
\mu=\frac1{\lambda}\big[\bar{\bar\Pi}\mu+\mathds{1}\big]
$, and replace $P_j$ with $vP_j$ for $j\in J$, leading to
$$\lambda=(1-v\psi_1)(1+\beta\mathds{1}),\quad
\beta=\frac{v}{\lambda}\big[\beta\bar{\bar\Pi}+\bar\bp\big],\quad
\mu=\frac1{\lambda}\big[v\bar{\bar\Pi}\mu+\mathds{1}\big]
,$$ where here and in the following $\la$, $\beta$, $\mu$ are short for $\la(v)$, $\beta(v)$, and $\mu(v)$. Rewriting the equation for $\beta$ in terms of $\tsup[2]{\Pi}$, we obtain
$\beta\big(\Id+\frac{v}{\la}\tsup[2]{\Pi}\big)
=\frac{v}{\la}(1+\beta\mathds{1})\bar\bp=\frac{v}{1-v\psi_1}\bar\bp,$
furthermore
\beq\label{betav}\beta=\frac{\la}{1-v\psi_1}\bar\bp\left(\frac{\la}{v}\Id+\tsup[2]{\Pi}\right)^{-1},\eeq and finally
$\frac{\la}{1-v\psi_1}-1=\beta\mathds{1}=\frac{\la}{1-v\psi_1}\Psi\left(\frac{v}{\la}\right),$
from which \eqref{hiprela} immediately follows.

For the proof of \eqref{hipreC}, we rewrite \eqref{hiprela} as
$\Psi(\frac v\la)-\psi_1\frac v\la=1-\frac1\la$ and differentiate w.r.t.\
$\frac v\la$, yielding
\beq\label{psiprime}\Psi'\left(\frac v\la\right)-\psi_1=\frac1{\la^2}\frac{\partial \la}{\partial\frac v\la}=\frac1{\la^2}\frac{\partial \la}{\partial v}
\left(\frac{\partial \frac v\la}{\partial v}\right)^{-1}
=\frac1{\la^2}\la'\left(\frac1\la-\frac{v\la'}{\la^2}\right)^{-1}
=\frac{\la'}{\la-v\la'}.\eeq
Rewriting the equation for $\mu$ in terms of $\tsup[2]{\Pi}$, we obtain
$\big(\Id+\frac{v}{\la}\tsup[2]{\Pi}\big)\mu
=\frac{1}{\la}\mathds{1}(1+v\bar\bp\mu)=\frac{1}{1-v\psi_1}\mathds{1},$
hence
\beq\label{muv}\mu=\frac{\la}{v(1-v\psi_1)}\left(\frac{\la}{v}\Id+\tsup[2]{\Pi}\right)^{-1}\mathds{1}.\eeq
Combining \eqref{betav} and \eqref{muv}, we obtain
$$(1-v\psi_1)^2\beta\mu=\frac{\la^2}{v}\bar\bp\Big(\frac{\la}{v}\Id+\tsup[2]{\Pi}\Big)^{-2}\mathds{1}=v\Psi'\left(\frac v\la\right).$$ 
This, and  \eqref{psiprime}, we plug into \eqref{C1gen}, thus establishing \eqref{hipreC},
\begin{align*}C(v)&=\frac{1+\beta\mathds{1}}{1+(1\!-\!v\psi_1)\beta\mu}
=\frac\la{1\!-\!v\psi_1+(1\!-\!v\psi_1)^2\beta\mu}=\frac\la{1+v(\Psi'(\frac v\la)-\psi_1)}=\frac\la{1+v\frac{\la'}{\la-v\la'}}\\
&=\la(v)-v\la'(v).\end{align*}
For the proof of \eqref{hipreP} observe that
$p_{\mathcal{I}}^{(n)}=\bp\bar\Pi^{n-1}\mathds{1}$ holds for $n\ge1$, with $\bp=[P_e,\bar\bp]$, and $\bar\Pi$ from the proof of Lemma~\ref{exp3ClaPhi}, yielding
$$\mathcal{P}_{\mathcal{I}}(z)
=1+\bp\Big(\frac1z\Id-\bar\Pi\Big)^{-1}\mathds{1}.$$
Let $\tsup[1]{\Pi}:=\mathds{1}\bp-\bar\Pi$ and observe 
$\bp\tsup[1]{\Pi}{}^n\mathds{1}=\bar\bp\tsup[2]{\Pi}{}^n\mathds{1}$ for 
$n\ge1$, as well as $\bp\mathds{1}=1,\, \bar\bp\mathds{1}=1-P_e$, which leads to
$$\bp(\tfrac1z\Id+\tsup[1]{\Pi})^{-1}\mathds{1}=
z\bp(\Id+z\tsup[1]{\Pi})^{-1}\mathds{1}=z-z^2\bp\tsup[1]{\Pi}(\Id+z\tsup[1]{\Pi})^{-1}\mathds{1}
=z-z\bar\bp\tsup[2]{\Pi}(\tfrac1z\Id+\tsup[2]{\Pi})^{-1}\mathds{1}
=P_ez+\Psi(z).
$$
By a well known resolvent identity, we have
$$(\tfrac1z\Id-\bar\Pi)^{-1}-(\tfrac1z\Id+\tsup[1]{\Pi})^{-1}
=(\tfrac1z\Id-\bar\Pi)^{-1}(\bar\Pi+\tsup[1]{\Pi})(\tfrac1z\Id+\tsup[1]{\Pi})^{-1},$$
and thus
$$\bp(\tfrac1z\Id-\bar\Pi)^{-1}\mathds{1}-\bp(\tfrac1z\Id+\tsup[1]{\Pi})^{-1}\mathds{1}
=\bp(\tfrac1z\Id-\bar\Pi)^{-1}\mathds{1}\,
\bp(\tfrac1z\Id+\tsup[1]{\Pi})^{-1}\mathds{1},$$
i.e.,
$$\mathcal{P}_{\mathcal{I}}(z)-1-(P_ez+\Psi(z))=(\mathcal{P}_{\mathcal{I}}(z)-1)(P_ez+\Psi(z)),$$
from which \eqref{hipreP} immediately follows.
\epr
Using \eqref{hiprela},
we can express $\la_1=\la(1)$ in terms of $(\psi_i)_{i\ge2}$ as follows,
\begin{align*}\la_1=
1&-\psi_2+\psi_3-(\psi_2^2+\psi_4)+(3\psi_2\psi_3+\psi_5)-(2\psi_2^3+4\psi_2\psi_4+2\psi_3^2+\psi_6)\\&+(10\psi_2^2\psi_3+5\psi_2\psi_5+5\psi_3\psi_4+\psi_7)\\
&-(5\psi_2^4+15\psi_2^2\psi_4+15\psi_2\psi_3^2+6\psi_2\psi_6
+6\psi_3\psi_5+3\psi_3^2+\psi_8)\\
&+(35\psi_2^3\psi_3+21\psi_2^2\psi_5+42\psi_2\psi_3\psi_4+7\psi_3^3
+7\psi_2\psi_7+7\psi_3\psi_6+7\psi_4\psi_5+\psi_9)+
\BO_{10}^*.
\end{align*}
This is found by computing the ninth Taylor polynomial of $\la(v)$ at $v=0$ and evaluating it at $v=1$. Clearly, more terms of $\la_1$ can easily be extracted using gfun.
Furthermore, by \eqref{hipreC}, we have
$$
C_1=1+\psi_2-2\psi_3+3(\psi_2^2+\psi_4)-4(3\psi_2\psi_3+\psi_5)+5(2\psi_2^3+4\psi_2\psi_4+2\psi_3^2+\psi_6)+\BO_7^*.$$
The expansion obtained from \eqref{hipreP} also turns out to use only $(\psi_i)_{i\ge2}$, and starts
\begin{align*}\mathcal{P}_{\mathcal{I}}(z)&=
1+z+(1-\psi_2)z^2+(1-2\psi_2+\psi_3)z^3+(1-3\psi_2+2\psi_3+\psi_2^2-\psi_4)z^4\\
&\hphantom{=1+z\,\,}+(1-4\psi_2+3\psi_3-2\psi_4+3\psi_2^2+\psi_5-2\psi_2\psi_3)z^5+\BO(z^6).
\end{align*}
To give an example of \eqref{hipreP} in action, consider the set of forbidden pairs $\mathcal{I}=\{(k,k),(k,\ell),(\ell,k)\}$ with $k\ne\ell$. Then we have 
$\tsup[2]{\Pi}=\left[\begin{smallmatrix}P_k&P_\ell\\ P_k&0\\ \end{smallmatrix}\right]$, which leads to 
$\Psi(z)=\left[\begin{smallmatrix}P_k&P_\ell\\ \end{smallmatrix}\right]
\Big(\frac1z\Id+\tsup[2]{\Pi}\Big)^{-1}
\left[\begin{smallmatrix}\vphantom{P_k}1\\ \vphantom{P_k}1\\ \end{smallmatrix}\right]=\frac{z(P_k+P_\ell-P_kP_\ell z)}{1+P_k z-P_kP_\ell z^2}$ and finally to $\mathcal{P}_{\mathcal{I}}(z)
=\frac{1+P_k z-P_kP_\ell z^2}
{1-(1-P_k)z-P_k(1-P_k-P_\ell)z^2+P_kP_\ell(1-P_k-P_\ell)z^3}$.

Using the function $\widetilde\Psi(z):=-\sum_{i\ge2}\psi_{i}(-z)^i$, equations \eqref{hiprela} and \eqref{hipreP} can be recast in the following, somewhat simpler forms,
$$\la=\frac{1}{1-\widetilde\Psi(\frac v\la)},\qquad\qquad
\mathcal{P}_{\mathcal{I}}(z)=\frac1{1-z-\widetilde\Psi(z)}.$$
We will meet the latter generating function again in Section~5,
where, employing a combinatorial approach, we are able to show that in case of one or two forbidden pairs, the generating function is rational with a denominator of degree at most three, which allows for very explicit expressions for the coefficients.

\ignore{In case of one or two forbidden pairs, it turns out that the generating function $\mathcal{P}_{\mathcal{I}}(z)$ is rational with a denominator of degree at most three. Therefore exact expressions for probabilities can readily be found using partial fraction decompositions.
For full details in this line of research see
section~5 in the extended preprint \cite{LSW22}, where a combinatorial approach is used for the computation of $\mathcal{P}_{\mathcal{I}}(z)$,
utilizing the methodology of the paper by \cite{Bassino}.}

\subsection {Limiting distribution of $X^{(n)}_3$}

\begin{conj}\label{conj:limdist1}
The asymptotic distribution of $X^{(n)}_3$ is Gaussian
\end{conj}
\bpr Note that the following proof is non-rigorous, as it is based on heuristic assumptions.

We assume asymptotic independence of $X_{i,j}^{(n)}$, as  the covariance total contribution is $\BO(1)$.
We consider pairs $(i,j)$ such that $i\neq j$.
The probability $\P[X_{i,j}^{(n)}=1]$ of pair $(i,j)$ occurring depends on $(i,j)$ only via $u:=i+j$, and is a decreasing function of $u$, which we denote $p_{n,u}$, with known asymptotics from \eqref{exp3} and \eqref{Pir}.
The number of pairs such that $i+j=u$ is given by
$c(u)=u-1- \kl even(u) \kr$.
Assuming that only the pairs most likely to occur, i.e., exactly those with $i+j\le\tilde u$ for some threshold $\tilde u$, contribute to $X_{i,j}^{(n)}$ (which we know is close to its expectation),
we are led to
$$\sum_{v=1}^{\ut}c(v)=
\sum_{v=1}^{\ut} (v-1)-\lf \frac{\ut}{2} \rf=\frac{\ut^2}{2}-\frac{\ut}{2}-\lf \frac{\ut}{2} \rf  \sim \frac{\ln(n)^2}{2L^2}+\BO(\ln(n)),$$
so we define $\tilde u:=\lf\frac{\ln(n)}{L}\rf$ to have a good match.
Taking into account also pairs $(i,j)$ with $i+j=\tilde u+1$, we have to
add a binomially distributed random variable
$\operatorname{Bin}\big(  c(\ut+1), p_{n,\ut+1}\big)$, which is asymptotically Gaussian. Similar corrections have to be added for
pairs $(i,j)$ with $i+j=\tilde u+k$ with $k\ge2$, the contributions rapidly becoming small as $k$ increases because of $p_{n,\tilde u+k}=\BO(q^k)$ as $k\to\infty$.
As some of the pairs with $i+j=\tilde u$ may be missing, we have to
subtract $\operatorname{Bin}\big(  c(\ut), 1-p_{n,\ut}\big)$.
Similar corrections have to be subtracted for
pairs $(i,j)$ with $i+j=\tilde u-k$ with $1\le k\le\ut-2$, all of these
corrections being asymptotically Gaussian. Again contributions rapidly become small as $k$ increases, because of $1-p_{n,\tilde u-k}\le\exp(-cq^{-k})$ for some $c>0$.
So the asymptotic total random contribution is Gaussian.
\epr

\begin{figure}[htbp]
	\centering
		\includegraphics[width=0.75\textwidth,angle=0]{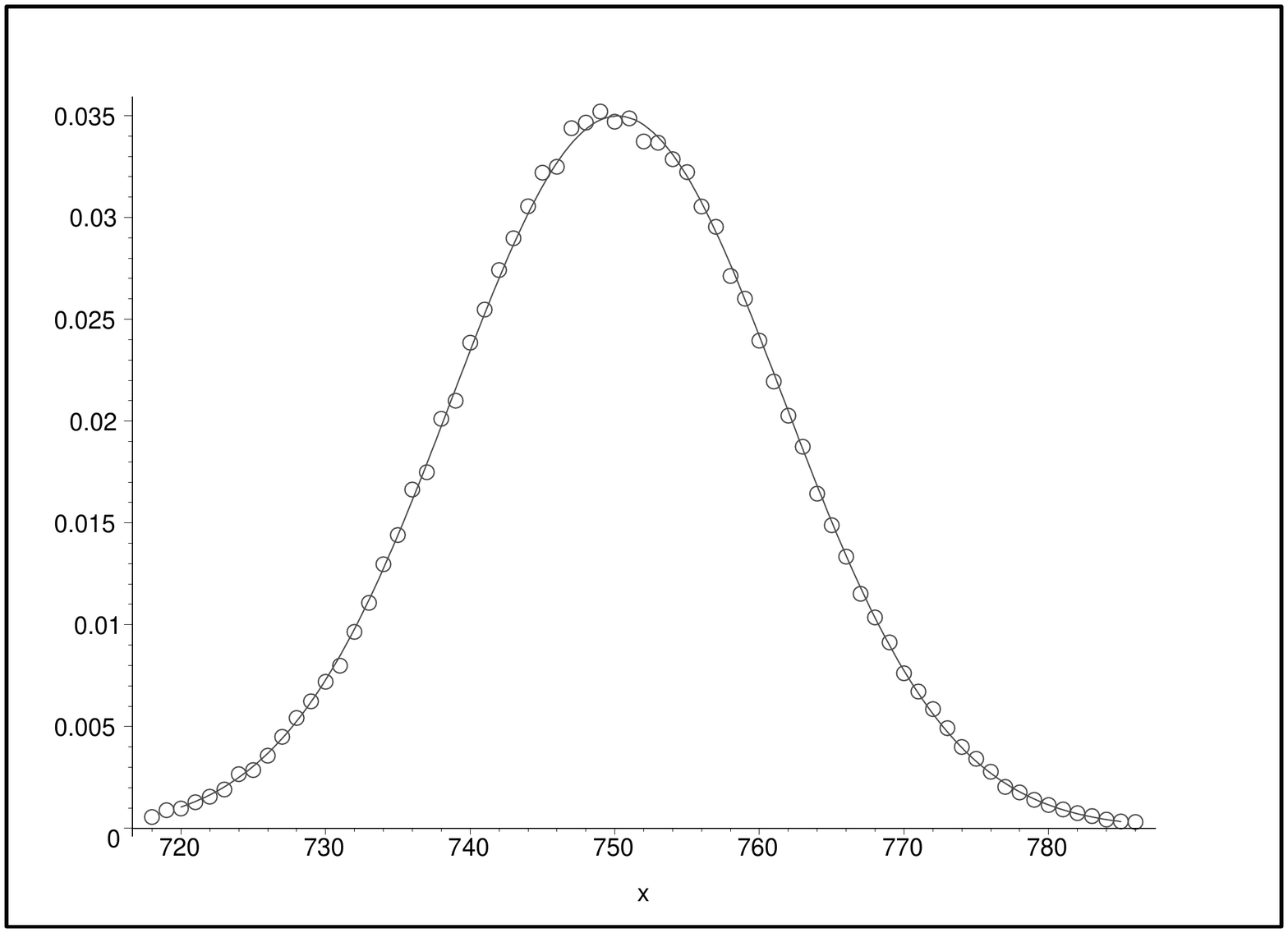}
	\caption{ Comparison between Gaussian density $f(x)$ (line) and the simulation of $X^{(n)}_3$ (circles), with $p=1/4,n=500000$, and number of simulated words $N=200000$.}
	\label{F3}
\end{figure}

The result of a simulation with $p=1/4$,  $n=500000$, and number of simulated words $N=200000$ can be seen in Figure \ref{F3}. The observed mean $\bar X^{(n)}_3\approx750.19$ and observed variance $s_3^2(n)\approx130.05$ are very close to
$\E X^{(n)}_3\approx750.19$ and $\Var X^{(n)}_3\approx129.88$. The density of a Gaussian with mean $\E X^{(n)}_3$ and variance $\Var X^{(n)}_3$ is also  shown in Figure \ref{F3}. The fit is excellent.

A rigorous proof of Conjecture~\ref{conj:limdist1} eludes us for now. What we have tried is the following. Define
random variables $\zeta^{(n)}_{i,j}$ with the same distribution as $X^{(n)}_{i,j}$, for $n,i,j\ge1$,
but such that for fixed $n$ the random variables $(\zeta^{(n)}_{i,j})_{i,j\ge1}$ are independent.
Furthermore define $\zeta^{(n)}:=\sum_{i\ne j}\zeta^{(n)}_{i,j}$, and let $\kappa^{(n)}_m$,
resp.\ $\bar\kappa^{(n)}_m$, be the $m$th cumulant of $\zeta^{(n)}$, resp.\ $X^{(n)}_3$.
Then show that
\begin{enumerate}
\item[i)] the sequence $(\zeta^{(n)})_{n\ge1}$ satisfies a CLT,
\item[ii)] the cumulants $\kappa^{(n)}_m$ and $\bar\kappa^{(n)}_m$ are close enough for the CLT proof to work also for $(X^{(n)}_3)_{n\ge1}$.
\end{enumerate}
Task i) is doable. We have $\kappa^{(n)}_2\sim S_2^{(n)}=\frac{\ln2}{\ln^2q}\ln n+\BO(1)$, by Theorem~\ref{thm:var123} and Lemma~\ref{var2sumeqsumvar2}, and can show $\kappa^{(n)}_m=\BO(\ln n)$ for
$m>2$. This gives $(\kappa^{(n)}_2)^{-\frac m2}\kappa^{(n)}_m\to0$ as $n\to\infty$, for each $m>2$, therefore, by the Frechet-Shohat theorem,
$(\Var \zeta^{(n)})^{-\frac12}(\zeta^{(n)}-\E\zeta^{(n)})$ converges in distribution to a standard normal random variable.
See section 4.7 in the extended preprint of \cite{LSW22} for first steps in the sketched direction.

For task ii), we know $|\kappa^{(n)}_2-\bar\kappa^{(n)}_2|\sim T^{(n)}_2=\BO(1)$.
Thus a bound like $|\kappa^{(n)}_m-\bar\kappa^{(n)}_m|=\BO(1)$ (or even $|\kappa^{(n)}_m-\bar\kappa^{(n)}_m|=\BO(\ln^{m/2-\eps}n)$ with some $\eps>0$) holding for
$m>2$ would guarantee the above CLT argument to carry over to the sequence $(X^{(n)}_3)_{n\ge1}$.
Now $\kappa^{(n)}_m-\bar\kappa^{(n)}_m$ involves infinite sums of mixed $m$th moments, and we are not quite sure, if our methods to deal with covariances would easily adapt to higher moments. Moreover the number of cases to distinguish (analogous to the 6 cases we had for $m=2$) grows rapidly with $m$.
So, unfortunately, we can not report progress here.%
\ignore{
\subsection{Another approach to the asymptotic distribution of $X^{(n)}_3$
}
Let us  turn to the couple $(i,j)$. Let us first assume  full independence. We first consider $X^{(n)}_2:=\sum_{i,j}  X_{i,j}^{(n)}$. The mean is simply given by
$M(X^{(n)}_2)=\sum_{i,j}  E(X_{i,j}^{(n)})$. The variance $VAR(X^{(n)}_2)$ is made of two parts: the first part is related to
$\sum_{i,j}  VAR(X_{i,j}^{(n)})$, the second part is related to $\sum_{i,j}\sum_{u,v} COV(X_{i,j}^{(n)},X_{u,v}^{(n)})$ which, as we have seen,
is only due to $\sum_{i,j,k} COV(X_{i,j}^{(n)},X_{j,k}^{(n)})$. So we first proceed as in the $(i,i)$ case. Next  we will deal with
$X^{(n)}_3:=\sum_{i\neq j}  X_{i,j}^{(n)}$ i.e.  subtract the contribution of the $(i,i)$ couples.  Finally, we will deal with the dependence property.
\subsubsection{Mean and $\sum_{i,j}  VAR(X_{i,j}^{(n)})$ contribution}
Set
\bals
&G_{n,2}(\tet):=E \lp  e^{\tet X^{(n)}_2}\rp,\quad S_{n,2}(\tet):=\ln \lp G_{n,2}(\tet)\rp.\\
\end{align*}
We proceed as in Sec. \ref{S31}.
\bals
\mbox{ Set }Q_2:=&\frac{1}{q},L_2:=\ln(Q_2)=-\ln(q),\ns:=n_2 k,\log:=\log_{Q_2},L_1=2 L_2,\\
\P(X_{i,j}(0)=1) \sim& \Pb(i,j)=e^{-n_2 q^{i+j}},\\
V_v=&\sum_i\sum_j\sum_{k=0}^v (-1)^k\bin{v}{k}\Pb(i,j)^k=\sum_i\sum_j\sum_{k=1}^v (-1)^{k+1}\bin{v}{k}\lb 1-\Pb(i,j)^k\rb,\\
&\sum_i\sum_j \lb 1-\Pb(i,j)^k\rb =\sum_{u=2} (u-1)\lb 1-e^{-\ns q^u} \rb =\sum_{u=1} (u-1)\lb 1-e^{-\ns q^u} \rb\\
&=\sum_{u=1} u\lb 1-e^{-\ns q^u} \rb-\sum_{u=1} \lb 1-e^{-\ns q^u} \rb:= V_{v,1}-V_{v,2},
\end{align*}
The second summation leads, using the previous results on $(i,i)$ to
$$\sum_{u=1} \lb 1-e^{-\ns q^u} \rb=U_1(Q_2)+\log_{Q_2}k+\be_k(Q_2),$$
where
$$
U_1(Q_2):=\frac{\ln(n_2)}{L_2}-\frac12+\frac{\gam}{L_2},\be_k(Q_2)=-\frac{1}{L_2}\sum_{\ell \in \Z\setminus \{0\}} \Gam\lp \frac{2\ii \ell\pi}{L_2}\rp e^{-2 \ii \ell\pi\log_{Q_2} (kn_2)},$$
This leads to
$$
V_{v,2}=-\sum_{k=1}^v (-1)^{k+1}\bin{v}{k}[U_1(Q_2)+\log_{Q_2}k+\be_k(Q_2)],$$
omitting the periodic contribution,
\bals
V_{v,2}=&-U_1(Q_2)+C_1(v),C_1(1)=0, C_1(2)=\log_{Q_2}2\\
V_{1,2}=&-U_1(Q_2)-\be_1(Q_2),\\
V_{2,2}=&-U_1(Q_2)+ \log_{Q_2}2  -\sum_{k=1}^2
(-1)^{k+1}\bin{2}{k}\be_k(Q_2).
\end{align*}
From Appendix \ref{A11},  the first summation gives
\bals
V_{v,1}=&\sum_{k=1}^v (-1)^{k+1}\bin{v}{k}\lb \frac{\ln
(n_2)^2+2\ln(n_2)\ln(k)+\ln(k)^2}{2L_2^2}+\frac{\gam(\ln(n_2)+\ln(k))}{L_2^2}\right.\\
&\left. -\lb -\frac{\pi^2/12+\gam^2/2}{L_2^2}+\frac{1}{12}\rb \rb+\bet_{v,1},\\
=&U_2(Q_2)+B_{v,1}+\bet_{v,1},
\end{align*}
where
\bals
U_2(Q_2)=& \frac{\ln(n_2)^2}{2L_2^2}+\frac{\gam\ln(n_2)}{L_2^2}
-\lb -\frac{\pi^2/12+\gam^2/2}{L_2^2}+\frac{1}{12}\rb \\
&=\frac{\ln(n)^2}{2L_2^2}+\ln(n)\lb {\frac {{ \gam}+2\,\ln  \left( p \right) +2\,{\it L_2}}{{{ L_2}}^{2
}}} \rb\\
&+1/12\,{\frac {{\pi}^{2}+6\,{{\gam}}^{2}+23\,{{L_2}}^{2}+24\,{
\gam}\,\ln  \left( p \right) +24\,{\gam}\,{L_2}+24\, \left(
\ln  \left( p \right)  \right) ^{2}+48\,\ln  \left( p \right) {L_2}
}{{{L_2}}^{2}}},
\end{align*}
Now define
$$
B_{v,1}:=\sum_{k=2}^v (-1)^{k+1}\bin{v}{k} \lb
\frac{2\ln(n_2)\ln(k)+\ln(k)^2}{2L_2^2}+\frac{\gam\ln(k)}{L_2^2}\rb=C_2(v)\ln(n)+C_3(v),$$
with, say,
$C_2(1)=C_3(1)=0$, and
\bals
B_{1,1}=&0,\\
B_{2,1}=&-\lb  \frac{2\ln(n_2)\ln(2)+\ln(2)^2}{2L_2^2}+\frac{\gam\ln(2)}{L_2^2}\rb,\\
\bet_{v,1}=&\sum_{k=1}^v (-1)^{k+1} \bin{v}{k} \be_{k,1},\bet_{1,1}=\be_{1,1},\bet_{2,1}=2 \be_{1,1}-\be_{2,1},\\
\be_{k,1}=&\frac{1}{L_2}\sum_{\ell \in \Z\setminus \{0\}}\lb -\frac{\ln(n_2 k)}{L_2} \Gam\lp \frac{2\ii \ell\pi}{L_2}\rp
+\frac{\Gam'\lp \frac{2\ii \ell\pi}{L_2}\rp}{L_2}\rb  e^{-2 \ii \ell\pi\log (k np^2)},\\
V_{1,1}=&U_2(Q_2)+\be_{1,1},\\
V_{2,1}=&U_2(Q_2)+B_{2,1}+\bet_{2,1},
\end{align*}
Finally, omitting the periodic contribution,
\bals
V_1=&V_{1,1}+V_{1,2}\sim U_2(Q_2)-U_1(Q_2),\\
V_2=&V_{2,1}+V_{2,2}\sim U_2(Q_2)-U_1(Q_2)+B_{2,1}+C_1(2),\\
V_v\sim& U_2(Q_2)-U_1(Q_2)+
C_2(v)\ln(n)+C_4(v),C_4(v)=C_1(v)+C_3(v),C_4(1)=0,
\end{align*}
We can now compute the mean $M_1(X^{(n)}_2)$ and $VAR_1(X^{(n)}_2)$
(see~(\ref{E40})).  We get
\bals
M_1(X^{(n)}_2)&=V_1=U_2(Q_2)-U_1(Q_2)+\be_{1,1}-\be_1(Q_2),\\
VAR_1(X^{(n)}_2)&= \K_2(X^{(n)}_2)=V_1-V_2=-B_{2,1}- \log_{Q_2}2+\be_{1,1}-\be_1(Q_2)-\bet_{2,1}  +\sum_{k=1}^2 (-1)^{k+1}\bin{2}{k} \be_k(Q_2),\\
&=\lb  \frac{2\ln(n_2)\ln(2)+\ln(2)^2}{2L_2^2}+\frac{\gam\ln(2)}{L_2^2}\rb- \log_{Q_2}2- \be_{1,1}+\be_{2,1}+\be_1(Q_2) -\be_2(Q_2)
\end{align*}
\subsubsection{Subtraction of  the contribution from all couples  $(i,i)$}
We obtain for the non-periodic contribution $M_{NP}(X^{(n)}_3)$
\bals
M_{NP}(X^{(n)}_3)&:=U_2(Q_2)-U_1(Q_2)-U_1(Q_1)=\frac{\ln(n)^2}{2L_2^2}
+\ln(n)\lb 1/2\,{\frac {{L_2}+2\,{\gam}+4\,\ln  \left( p \right) }{{{L_2}}^{2}}} \rb  \\
&-1/12{\frac {-{\pi}^{2}-6\,{{\gam}}^{2}+{{L_2}}^{2}-24\,{\it
\gam}\,\ln  \left( p \right) -6\,{\gam}\,{L_2}-24\, \left( \ln
 \left( p \right)  \right) ^{2}-12\,\ln  \left( p \right) {L_2}}{{{
L_2}}^{2}}},
\end{align*}
Next, for the periodic contribution $M_{P}(X^{(n)}_3)$, we obtain
\bals
M_{P}(X^{(n)}_3)&=\be_{1,1}-\be_1(Q_2) -\be_1(Q_1)=\frac{1}{L_2}\sum_{\ell \in \Z\setminus \{0\}}\lb -\frac{\ln(n_2 )}{L_2} \Gam\lp \frac{2\ii \ell\pi}{L_2}\rp
+\frac{\Gam'\lp \frac{2\ii \ell\pi}{L_2}\rp}{L_2}\rb e^{-2 \ii \ell\pi\log_{Q_2} ( np^2)}\\
&+\frac{1}{L_2}\sum_{\ell \in \Z\setminus \{0\}}\lb\Gam\lp \frac{2\ii \ell\pi}{L_2}\rp\rb  e^{-2 \ii \ell\pi\log_{Q_2} ( np^2)}\\
&+\frac{1}{L_1}\sum_{\ell \in \Z\setminus \{0\}}\lb\Gam\lp \frac{2\ii \ell\pi}{L_1}\rp\rb  e^{-2 \ii \ell\pi\log_{Q_1} ( np^2)}.\\
\end{align*}
Similarly, the variance of $X^{(n)}_3$ (without dependence property) is given by
\bals
VAR_2(X^{(n)}_3)&=\lb
\frac{2\ln(n_2)\ln(2)+\ln(2)^2}{2L_2^2}+\frac{\gam\ln(2)}{L_2^2}\rb-
\log_{Q_2}2- \be_{1,1}+\be_{2,1} +\be_1(Q_2)-\be_2(Q_2)\\
&\qquad{}-[\log_{Q_1} 2+\be_2(Q_1)-\be_1(Q_1)].
\end{align*}
After some algebra, we can check that this is identical to Thm \ref{thm:var123} (except for the term $T^{(n)}_2$, of course).
Also  the mean is identical  with the expression given in \cite{ABCKWW21}, Thm 3. Our  approach  is simple and
general. Note again that the mean  does not rely on the dependence property.
\subsubsection{The dominant terms, without the periodic components}
 We have, successively (we keep only dominant terms), and setting
$$\Xt_{i,j}^{(n)} :=X_{i,j}^{(n)} -\E(X_{i,j}^{(n)} ),
\qquad\qquad\mu_k(X):=\E[(X-\E(X))^k],$$
we have
$$
S_{n_2}(\tet)  =\sum_{v=1}^\II \frac{(-1)^{v+1}}{v}\lp e^\tet -1 \rp^v V_v,$$
We first deal without dependence property,
\bals
 \K_2(X^{(n)}_2)&=\mu_2(X^{(n)}_2)=VAR_1(X^{(n)}_2)=-C_2(2)\ln(n)-C_4(2),\\
 \K_3(X^{(n)}_2)&=\mu_3(X^{(n)}_2)\\
&=V_1-3V_2+2V_3\\
&=[-3C_2(2)+2C_2(3)]\ln(n)+[-3C_4(2)+2C_4(3)]\\
&=C_5\ln(n)+C_6, \mbox{ say},
\end{align*}
but
$$
 \K_2(X^{(n)}_2)=\E\lb \lp   \sum_{i\neq j}  \Xt_{i,j}^{(n)}  +\sum_{i}  \Xt_{i,i}^{(n)} \rp^2  \rb=\K_2(X^{(n)}_3)+\K_2(X^{(n)}_1)=\K_2(X^{(n)}_3)+C_7,\mbox{ say },$$
hence,
$$
 VAR_2(X^{(n)}_3)=\mu_2(X^{(n)}_3)=\K_2(X^{(n)}_3)=-C_2(2)\ln(n)-C_4(2)-C_7=-C_2(2)\ln(n)+C_8,\mbox{
   say },$$
Now we calculate
$$
 \K_3(X^{(n)}_2)=\E\lb \lp   \sum_{i\neq j}  \Xt_{i,j}^{(n)}  +\sum_{i}  \Xt_{i,i}^{(n)} \rp^3  \rb=\K_3(X^{(n)}_3)+\K_3(X^{(n)}_1)=\K_3(X^{(n)}_3)+C_9,\mbox{ say },
$$
hence,
$$
\K_3(X^{(n)}_3)=\mu_3(X^{(n)}_3)=C_5\ln(n)+C_6-C_9=C_5\ln(n)+C_{10},\mbox{ say }.$$
Now we deal  with dependence property,
$$
\sigd=VAR(X^{(n)}_3)=\K_2(X^{(n)}_3)+T^{(n)}_2=-C_2(2)\ln(n)+C_{11},\mbox{ say },$$
Then
$\Kb_3(X^{(n)}_3)$ is the new third cumulant, including the
correlations, and $T_3$ is the correction, similar to $T_2$, so that
$$
 \Kb_3(X^{(n)}_3)=\K_3(X^{(n)}_3)+T_3=C_5\ln(n)+C_{10}+T_3.$$
A bound such that  $|T_3|=\BO(1)$ would be enough.
We now have
$$
 \E\lb  \exp\lb  \tet\frac{X^{(n)}_3-M(X^{(n)}_3)}{\sig} \rb\rb \sim\exp\lb \frac{\tet^2}{2}+\frac{\tet^3}{3!}\frac{C_5\ln(n)+C_{10}+T_3}{(-C_2(2)\ln(n)+C_{11})^{3/2}}+\ldots\rb
$$
Let us, for instance, deal with the fourth cumulant. We have, without dependence
$$
\K_4(X^{(n)}_2)=V_1-7V_2+12V_3-6V_4=C_{12}\ln(n)+C_{13}, \mbox{ say },$$
and
$$\mu_4(X^{(n)}_2)=\K_4(X^{(n)}_2)+3\K_2^2(X^{(n)}_2),$$
but
$$\mu_4(X^{(n)}_2)=\mu_4(X^{(n)}_3)+\mu_4(X^{(n)}_1)+6\mu_2(X^{(n)}_3)\mu_2(X^{(n)}_1),$$
hence
$$\mu_4(X^{(n)}_3)=\mu_4(X^{(n)}_2)-[\mu_4(X^{(n)}_1)+6\mu_2(X^{(n)}_3)\mu_2(X^{(n)}_1)]$$
and
$$\mu_4(X^{(n)}_3)=\K_4(X^{(n)}_2)+3\K_2^2(X^{(n)}_2) -\mu_4(X^{(n)}_1)-6\mu_2(X^{(n)}_3)\mu_2(X^{(n)}_1);$$
Now, with dependence,
\bals
\mub_4(X^{(n)}_3)&=\mu_4(X^{(n)}_3)+T_4,\\
\Kb_4(X^{(n)}_3)&=\K_4(X^{(n)}_2)+3\mu_2^2(X^{(n)}_2)-\mu_4(X^{(n)}_1)-6\mu_2(X^{(n)}_3)\mu_2(X^{(n)}_1)+T_4-3\mu_2^2(X^{(n)}_3)-6\mu_2(X^{(n)}_3)T_2-3T_2^2\\
&=\K_4(X^{(n)}_2)+3\mu_2^2(X^{(n)}_1)-\mu_4(X^{(n)}_1)-6\mu_2(X^{(n)}_3)T_2+T_4-3T_2^2\\
&=C_{14}\ln(n)+C_{15}+T_4-3T_2^2, \mbox{ say. }
\end{align*}
Again a  bound such that $|T_4|=\BO(1)$ would be enough.
We proceed similarly for all cumulants.
\begin{conj}\label{conj:limdist2}
\[|T_.|=\BO(1).\]
\end{conj}
 This  leads to the asymptotic Gaussian distribution.
}
\ignore{{\bf An attempt to show finiteness of $T_3$}
We are more generally interested in
terms (again $n$ indicates the length of the word)
$$\E_n \biggl[\prod_{i\in I}\tilde X_{k_i,m_i}^{(n)}\biggr]=\biggl[\sum_{I'\subseteq I}(-1)^{|I'|}\rho^{(n)}_{I'}\biggr]\prod_{i\in I}p^{(n)}_{\{i\}},$$
where $\rho^{(n)}_{I'}:=\dfrac{p^{(n)}_{I'}}{\prod_{i\in I'}p^{(n)}_{\{i\}}}$, for $I'\subseteq I$. Stressing the dependence on $I$, we now rewrite \eqref{ClaPhi}
as $$\P(X_{k_i,m_i}^{(n)}=0,i\in I)=C_I\la_I^n\Phi_{I,n},$$  and define $\Lambda_I:=\frac{\la_I}{\prod_{i\in I}\la_{\{i\}}}$, and $\Gamma_{I,n}:=\frac{C_I\Phi_{I,n}}{\prod_{i\in I}C_{\{i\}}\Phi_{\{i\},n}}$, such that $\rho_I=\Gamma_{I,n}\Lambda_I^n$.
For fixed $\delta>0$ there is a $C^\infty$ function $R_n$ depending on $\bar\bp\in\mathcal{D}_\delta^J:=\{\bx\in\mathcal{D}^{J}:\bx^t\mathds{1}\le1-\delta\}$,
such that $\rho_I^{(n)}=R_n(\bar\bp)$, and $\Lambda_I$ and $\Gamma_{I,n}$ are $C^\infty$ functions of $\bar\bp\in\mathcal{D}_\delta^J$ as well. Using $\ln\la_I=-\psi_2+\psi_3-(\frac32\psi_2^2+\psi_4)+(4\psi_2\psi_3+\psi_5)+\BO_6^*$, we obtain
$$\Lambda_I=e^{\bar\psi_3-\bar\psi_4+\bar\psi_5+\BO_6^*},$$
with
\begin{align*}\bar\psi_3&:=\psi_3-\sum_{(k,k)\in\mathcal{I}}P_k^3,\\ \bar\psi_4&:=(\tfrac32\psi_2^2+\psi_4)-\tfrac32\sum_{(k,m)\in\mathcal{I}}P_k^2P_m^2
-\sum_{(k,k)\in\mathcal{I}}P_k^4,\\
\bar\psi_5&:=(4\psi_2\psi_3+\psi_5)-5\sum_{(k,k)\in\mathcal{I}}P_k^5,\end{align*}
and similarly $\Gamma_{I,n}=1+\BO(\bar\psi_3)+\BO(\bar\psi_4)$.
Some of the terms $\rho_{I'}^{(n)}$ can be expressed by setting some of the variables $P_\cdot$ to zero in $R_n(P_{j_1},\ldots,P_{j_{|J|}})$, and the MVT may be applied.
{\bf Example:} Consider a set of $3$ pairs $\mathcal{I}=\{(k_i,m_i):i\in I\}=\{(i,j),(j,k),(k,\ell)\}$, with $J:=\{i,j,k,\ell\}$ of size 4. We suppress $n$.
Then $$\rho_{I}=\dfrac{\P(X_{i,j}^{(n)}=X_{j,k}^{(n)}=X_{k,\ell}^{(n)}=0)}{\P(X_{i,j}^{(n)}=0)\P( X_{j,k}^{(n)}=0)\P(X_{k,\ell}^{(n)}=0)}=R(P_i,P_j,P_k,P_\ell),$$
and moreover
$$
\rho_{\{1,2\}}=\dfrac{\P(X_{i,j}^{(n)}=X_{j,k}^{(n)}=0)}{\P(X_{i,j}^{(n)}\!=\!0)\P( X_{j,k}^{(n)}\!=\!0)}=R(P_i,P_j,P_k,0),\quad
\rho_{\{2,3\}}=\dfrac{\P(X_{j,k}^{(n)}=X_{k,\ell}^{(n)}=0)}{\P(X_{j,k}^{(n)}\!=\!0)\P( X_{k,\ell}^{(n)}\!=\!0)}=R(0,P_j,P_k,P_\ell),
$$
$$\rho_{\{1,3\}}=\dfrac{\P(X_{i,j}^{(n)}=X_{k,\ell}^{(n)}=0)}{\P(X_{i,j}^{(n)}=0)\P( X_{k,\ell}^{(n)}=0)},\qquad \rho_{\{2\}}=R(0,P_j,P_k,0)=1,\qquad \rho_{\{1\}}=\rho_{\{3\}}=\rho_{\{\}}=1.$$
Therefore
$$
\sum_{I'\subseteq I}(-1)^{|I'|}\rho_{I'}=-\underbrace{\lb R(P_i,P_j,P_k,P_\ell)\!-\!R(0,P_j,P_k,P_\ell)\!-\!R(P_i,P_j,P_k,0)\!+\!R(0,P_j,P_k,0)\rb}_{P_iP_\ell\frac{\partial^2R}{\partial P_i\partial P_\ell}(p_i,P_jP_k,p_\ell)\ \textup{ for some }0<p_i<P_i,0<p_\ell<P_\ell}+
\underbrace{\lb\rho_{\{1,3\}}-1\rb}_{\textup{see Case~4}}.
$$
Observe that $\lim_{P_j\searrow0}\Lambda_I=\lim_{P_k\searrow0}\Lambda_I=
\lim_{P_i+P_\ell\searrow0}\Lambda_I=1$, therefore $\Lambda_I=1+\BO(P_iP_jP_k)+\BO(P_jP_kP_\ell)$, and similarly $\Gamma_{I,n}=1+\BO(P_iP_jP_k)+\BO(P_jP_kP_\ell)$, uniformly in $n$, this uniformity being shown as in the proof of Lemma~\ref{lemCovir}.
Clearly, using the expression preceding \eqref{hiprela}, a much more precise term for $\Lambda_I$ can be found via $\la_I$:
We have
$\psi_2=P_iP_j+P_jP_k+P_kP_\ell$, $\psi_3=P_iP_jP_k+P_jP_kP_\ell$, and
$\psi_4=P_iP_jP_kP_\ell$,
resulting in
$$\la_I=1-(P_iP_j+P_jP_k+P_kP_\ell)+(P_iP_jP_k+P_jP_kP_\ell)-(P_iP_j+P_jP_k+P_kP_\ell)^2-P_iP_jP_kP_\ell+\BO_5^*,$$
but $R(P_i,P_j,P_k,P_\ell)=1+\BO(nP_iP_jP_k)+\BO(nP_jP_kP_\ell)$ fits our purposes for now. Both underbraced terms above are $\BO(nP_iP_jP_kP_\ell)$, and the following contribution to $T_3$ can now be bounded:
\begin{align*}
\sum_{i,j,k,\ell\ge1:|\{i,j,k,\ell\}|=4}\biggl|\E \lb\tilde X_{i,j}^{(n)}\tilde X_{j,k}^{(n)}\tilde X_{k,\ell}^{(n)}\rb\biggr|&=\BO\biggl(
\sum_{i,j,k,\ell\ge1}nP_iP_jP_kP_\ell e^{-nP_iP_j-nP_jP_k-nP_kP_\ell}\biggr)\\
&=\BO\biggl(
\sum_{i,j,k,\ell\ge1}n^\frac12P_i(P_jP_k)^\frac12P_\ell e^{-nP_iP_j-nP_kP_\ell}\biggr)=\BO\lp\frac1{\sqrt{n}}\rp.
\end{align*}}

\section{Combinatorial Pattern Matching Approach}

For a combinatorial approach, we utilize the methodology of \cite{Bassino}.  The full strength of~\cite{Bassino}
is not needed, because (in the present analysis) we are only studying
``reduced'' sets of patterns.
In a reduced set of patterns, no word is a subword of another word.
Here, we are always analyzing patterns
of length~2, so our patterns are necessarily (already) reduced.  So we
\emph{only} need to understand Sections~4.1 and 4.2
of~\cite{Bassino}.

Since we follow the notation and overall approach of~\cite{Bassino},
the reader might want to review the first 10 pages of~\cite{Bassino},
through Section~4.2.  The basic methodology is to use an
inclusion-exclusion approach to enumerating patterns.  This approach
allows an \emph{exact} derivation of the probabilities of each set of
patterns.  For this approach, Section~4.1 of~\cite{Bassino} explains
how to utilize decorated texts, in which some occurrences of patterns
are ``distinguished'' (while others might not be distinguished).

Collections of overlapping distinguished texts are gathered together
into clusters.  With this methodology, ``the set of decorated
texts~$T$ decomposes as sequences of either arbitrary letters of the
alphabet~$\mathcal{A}$ or clusters: $T = (\mathcal{A} + C)^{*}$''.
Using $\xi(z,t) = \sum_{w\in C}\pi(w)z^{|w|}t^{\tau(w)}$, where
$\pi(w)$ is the probability of a text, and $\tau(w)$ is the number of
distinguished occurrences of subwords in $w$, the generating function
of all decorated texts is $T(z,t) = 1/(1-A(z)-\xi(z,t))$.  

Finally, using inclusion-exclusion, it follows that the probability generating
function $F_{\mathcal{U}}(z,x)$, in which powers of $z$ mark the
length of texts, and powers of $x$ mark the total number of
occurrences of patterns in $\mathcal{U}$, we obtain
$F_{\mathcal{U}}(z,x) = 1/(1-A(z)-\xi(z,x-1))$.  This is the set of
core ideas from~\cite{Bassino} that forms the foundation of the
analysis in the present section.

We define $X^{(n)}$ as the total number of distinct (adjacent) pairs
in a word $ Z_{1},\ldots,Z_{n}$, and we have
$$X^{(n)} = \sum_{i=1}^{\infty}\sum_{j=1}^{\infty}X_{i,j}^{(n)}$$

\begin{note}\label{partialfractiondecomposition1}
The roots of the polynomials in the denominators of the generating
functions in Table~\ref{bigtable1}
and in Table~\ref{bigtable2} exist and are unique (or
there is a removable singularity that can be defined by
using continuity).
\end{note}

\begin{table}[!htb]
\begin{center}
\renewcommand{\arraystretch}{2.1}
\begin{tabular}{l|l|l}
\hline\hline
A & \em gen.\ func. & $\frac{1}{1-z+P_i P_j z^{2}}$ \\
\hline
& \em par.\ frac. & $\frac{1 - \sqrt{1-4P_i P_j }}{2P_i P_j\sqrt{1-4P_i P_j }}
\lp1-\frac{2P_i P_j }{1 - \sqrt{1-4P_i P_j }}z\rp^{-1}-
\frac{1 + \sqrt{1-4P_i P_j }}{2P_i P_j\sqrt{1-4P_i P_j }}
\lp1-\frac{2P_i P_j }{1 + \sqrt{1-4P_i P_j }}z\rp^{-1}
$ \\
\hline
& \em coeff.\ of $z^{n}$ & $\frac{(2P_i P_j )^{n+1}}{\sqrt{1-4P_i P_j }}
\bigg(
 \frac{1}{
\big(1 - \sqrt{1-4P_i P_j }\big)^{n+1}}
-\frac{1}{
\big(1 + \sqrt{1-4P_i P_j }\big)^{n+1}}\bigg)$ \\
\hline\hline\hline
B & \em gen.\ func. & $\lp 1-z+\frac{P_i^2 z^{2}}{1 + P_i z}\rp^{-1}
= \frac{1 + P_i z}{1-(1-P_i )z-P_i (1-P_i )z^{2}}$ \\
\hline
& \em par.\ frac. & $\lp\frac12-\frac{1+P_i}{2\sqrt{(1-P_i )(1+3P_i )}}\rp\lp1-\frac{-2P_i (1-P_i )}{1-P_i +\sqrt{(1-P_i )(1+3P_i )}}z\rp^{-1}$\\
& & $\qquad\qquad{}-\lp\frac12+\frac{1+P_i}{2\sqrt{(1-P_i )(1+3P_i )}}\rp\lp1-\frac{-2P_i (1-P_i )}{1-P_i -\sqrt{(1-P_i )(1+3P_i )}}z\rp^{-1}$ \\
\hline
& \em coeff.\ of $z^{n}$ & $\lp\frac12-\frac{1+P_i}{2\sqrt{(1-P_i )(1+3P_i )}}\rp\lp\frac{-2P_i (1-P_i )}{1-P_i +\sqrt{(1-P_i )(1+3P_i )}}\rp^{n}$\\
& &
$\qquad\qquad{}-\lp\frac12+\frac{1+P_i}{2\sqrt{(1-P_i )(1+3P_i )}}\rp\lp\frac{-2P_i (1-P_i )}{1-P_i -\sqrt{(1-P_i )(1+3P_i )}}\rp^{n}$ \\
\hline\hline\hline
C & \em gen.\ func. & $\frac{1}{1-z+P_i P_j z^{2} + P_k P_\ell z^{2}}$ \\
\hline
& \em par.\ frac. & {\small $\frac{1 + \sqrt{1-4(P_iP_j\!+\!P_kP_\ell) }}{2\sqrt{1-4(P_iP_j\!+\!P_kP_\ell) }}
\Big(\!1\!-\!\frac{2(P_iP_j\!+\!P_kP_\ell) }{1 - \sqrt{1-4(P_iP_j\!+\!P_kP_\ell) }}z\!\Big)^{\!\!-1}\!\!\!-\!
\frac{1 - \sqrt{1-4(P_iP_j\!+\!P_kP_\ell) }}{2\sqrt{1-4(P_iP_j\!+\!P_kP_\ell) }}
\Big(\!1\!-\!\frac{2(P_iP_j\!+\!P_kP_\ell) }{1 + \sqrt{1-4(P_iP_j\!+\!P_kP_\ell) }}z\!\Big)^{\!\!-1}
$} \\
\hline
& \em coeff.\ of $z^{n}$ & $
\frac{(2(P_i P_j  + P_k P_\ell ))^{n+1}}
{\sqrt{1-4(P_i P_j  + P_k P_\ell )}}\bigg(
\frac{1}{
\big(1-\sqrt{1-4(P_i P_j  + P_k P_\ell )}\big)^{n+1}}- \frac{1}{
\big(1+\sqrt{1-4(P_i P_j  + P_k P_\ell )}\big)^{n+1}}\bigg)$ \\
\hline\hline\hline
D & \em gen.\ func. & $\lp 1-z+\frac{P_i z(P_i z + P_\ell z)}{1 + P_i z}\rp^{-1}
= \frac{1 + P_i z}
{1 - (1 - P_i )z
+(P_i^2  + P_\ell P_i  - P_i )z^{2}}$ \\
\hline
& \em par.\ frac. & $\lp\frac12-\frac{1+P_i}{2\sqrt{(1-P_i )(1+3P_i )-4P_iP_\ell}}\rp\lp1-\frac{-2P_i (1-P_i-P_\ell )}{1-P_i +\sqrt{(1-P_i )(1+3P_i )-4P_iP_\ell}}z\rp^{-1}$\\
& & $\qquad\qquad{}-\lp\frac12+\frac{1+P_i}{2\sqrt{(1-P_i )(1+3P_i )-4P_iP_\ell}}\rp\lp1-\frac{-2P_i (1-P_i-P_\ell )}{1-P_i -\sqrt{(1-P_i )(1+3P_i )-4P_iP_\ell}}z\rp^{-1}$ \\
\hline
& \em coeff.\ of $z^{n}$ & $\lp\frac12-\frac{1+P_i}{2\sqrt{(1-P_i )(1+3P_i )-4P_iP_\ell}}\rp\lp\frac{-2P_i (1-P_i-P_\ell )}{1-P_i +\sqrt{(1-P_i )(1+3P_i )-4P_iP_\ell}}\rp^{n}$\\
& & $\qquad\qquad{}-\lp\frac12+\frac{1+P_i}{2\sqrt{(1-P_i )(1+3P_i )-4P_iP_\ell}}\rp\lp\frac{-2P_i (1-P_i-P_\ell )}{1-P_i -\sqrt{(1-P_i )(1+3P_i )-4P_iP_\ell}}\rp^{n}$\\
\hline\hline
\end{tabular}
\end{center}
\begin{small}
\caption{\label{bigtable1}\small
Table of generating functions, partial fraction decompositions, and
coefficients of~$z^{n}$, $n\ge2$, in each.
}
\end{small}
\end{table}




\begin{table}[!htb]
\begin{center}
\renewcommand{\arraystretch}{2.1}
\begin{tabular}{l|l|l}
\hline\hline
E & \em gen.\ func. & $\lp 1-z+\frac{P_i^2 z^{2}}{1 + P_i z}
+ \frac{P_k^2 z^{2}}{1 + P_k z}\rp^{-1}= \frac{(1 + P_i z)(1 + P_k z)}{(1-z)(1 + P_i z)(1 + P_k z)+P_i^2 z^{2}(1 + P_k z) + P_k^2 z^{2}(1 + P_i z)}$\\
\hline
& \em $a,b,c,d$ &
$a = P_i P_k (P_i  + P_k  - 1)$,
$b = P_i^2 + P_iP_k + P_k^2 -P_i - P_k$, $c = P_k +P_i -1$,
$d = 1$\\
\hline
& \em par.\ frac. & $\frac{(1 + P_i z)(1 + P_k z)}{az^{3} + bz^{2} + cz + d}$\\
& & ${}= (1 + P_i z)(1 + P_k z)
\big(\frac{rs}{(r-t)(s-t)(1-z/t)}
+ \frac{rt}{(r-s)(t-s)(1-z/s)}
+ \frac{st}{(s-r)(t-r)(1-z/r)}\big)$ \\
\hline
& \em coeff.\ of $z^{n}$ & $\frac{(1 + P_i t)(1 + P_k t)rs}{(r-t)(s-t)t^{n}}
+ \frac{(1 + P_i s)(1 + P_k s)rt}{(r-s)(t-s)s^{n}}
+ \frac{(1 + P_i r)(1 + P_k r)st}{(s-r)(t-r)r^{n}}$ \\
\hline\hline\hline
F & \em gen.\ func. & $\lp 1-z+\frac{P_i^2 z^{2}}{1 + P_i z}
+ P_k P_\ell z^{2}\rp^{-1}= \frac{1 + P_i z}{(1-z)(1 + P_i z)+P_i^2 z^{2}
+ P_k P_\ell z^{2}(1 + P_i z)}$\\
\hline
& \em $a,b,c,d$ &
$a = P_i P_k P_\ell $,
$b = P_i P_i  + P_k P_\ell -P_i $,
$c = P_i -1$,
$d = 1$\\
\hline
& \em par.\ frac. & $\frac{1 + P_i z}{az^{3} + bz^{2} + cz + d} = (1 + P_i z)
\big(\frac{rs}{(r-t)(s-t)(1-z/t)}
+ \frac{rt}{(r-s)(t-s)(1-z/s)}
+ \frac{st}{(s-r)(t-r)(1-z/r)}\big)$ \\
\hline
& \em coeff.\ of $z^{n}$ & $\frac{(1 + P_i t)rs}{(r-t)(s-t)t^{n}}
+ \frac{(1 + P_i s)rt}{(r-s)(t-s)s^{n}}
+ \frac{(1 + P_i r)st}{(s-r)(t-r)r^{n}}$ \\
\hline\hline\hline
G & \em gen.\ func. & $\frac{1}{1-z+P_iP_jz^2 + P_j P_\ell z^{2}-P_i P_j P_\ell z^{3}}$ \\
\hline
& \em $a,b,c,d$ &
$a = -P_i P_j P_\ell $,
$b = P_iP_j+P_j P_\ell $,
$c = -1$,
$d = 1$\\
\hline
& \em par.\ frac. & $\frac{1}{az^{3} + bz^{2} + cz + d} =
\frac{rs}{(r-t)(s-t)(1-z/t)}
+ \frac{rt}{(r-s)(t-s)(1-z/s)}
+ \frac{st}{(s-r)(t-r)(1-z/r)}$ \\
\hline
& \em coeff.\ of $z^{n}$ & $\frac{rs}{(r-t)(s-t)t^{n}}
+ \frac{rt}{(r-s)(t-s)s^{n}}
+ \frac{st}{(s-r)(t-r)r^{n}}$ \\
\hline\hline\hline
H & \em gen.\ func. & $\lp 1-z+\frac{2P_i P_j z^{2}}{1 - P_i P_j z^{2}}
- \frac{P_i^2 P_j z^{3}}{1 - P_i P_j z^{2}}
- \frac{P_i P_j^2 z^{3}}{1 - P_i P_j z^{2}}\rp^{-1}= \frac{1 - P_i P_j z^{2}}{(1-z)(1 - P_i P_j z^{2})
+ 2P_i P_j z^{2}
- P_i^2 P_j z^{3}
- P_i P_j^2 z^{3}}$\\
\hline
& \em $a,b,c,d$ &
$a = P_i P_j (1 - P_i  - P_j )$,
$b = P_i P_j $,
$c = - 1$,
$d = 1$\\
\hline
& \em par.\ frac. & $\frac{1 - P_i P_j z^{2}}{az^{3} + bz^{2} + cz + d} = (1 - P_i P_j z^{2})
\big(\frac{rs}{(r-t)(s-t)(1-z/t)}
+ \frac{rt}{(r-s)(t-s)(1-z/s)}
+ \frac{st}{(s-r)(t-r)(1-z/r)}\big)$ \\
\hline
& \em coeff.\ of $z^{n}$ & $\frac{(1 - P_i P_j t^{2})rs}{(r-t)(s-t)t^{n}}
+ \frac{(1 - P_i P_j s^{2})rt}{(r-s)(t-s)s^{n}}
+ \frac{(1 - P_i P_j r^{2})st}{(s-r)(t-r)r^{n}}$ \\
\hline\hline
\end{tabular}
\end{center}
\begin{small}
\caption{\label{bigtable2}\small
Table of generating functions, partial fraction decompositions, and
coefficients of~$z^{n}$, $n\ge2$, in each.
}
\end{small}
\end{table}

\begin{lem}\label{meanijneqlemma}
For $n\geq 2$, and for $i\neq j$,
the probability that $ij$ occurs
(at least once) as an adjacent pattern in
$ Z_{1},\ldots,Z_{n}$ is exactly
\begin{equation*}
E[X_{i,j}^{(n)}]
= 1
- \frac{(2P_i P_j )^{n+1}}{\sqrt{1-4P_i P_j }}\Bigg(
\frac{1}{\big(1 - \sqrt{1-4P_i P_j }\!\ \big)^{n+1}}
- \frac{1}{\big(1 + \sqrt{1-4P_i P_j }\!\ \big)^{n+1}}\Bigg).
\end{equation*}
\end{lem}
\bpr
The proof of Lemma~\ref{meanijneqlemma} is in subsection~\ref{meanijneqsection}.
\epr
\begin{lem}\label{meaniilemma}
For $n\geq 2$,
the probability that $ii$ occurs
(at least once) as an adjacent pattern in
$ Z_{1},\ldots,Z_{n}$ is exactly
\begin{align*}
E[X_{i,i}^{(n)}]
&= 1-\lp\frac12-\frac{1+P_i}{2\sqrt{(1-P_i )(1+3P_i )}}\rp\lp\frac{-2P_i (1-P_i )}{1-P_i +\sqrt{(1-P_i )(1+3P_i )}}\rp^{n}\\
& \qquad\qquad{}+\lp\frac12+\frac{1+P_i}{2\sqrt{(1-P_i )(1+3P_i )}}\rp\lp\frac{-2P_i (1-P_i )}{1-P_i -\sqrt{(1-P_i )(1+3P_i )}}\rp^{n}.
\end{align*}
Again, for $n<2$, we have $E[X_{i,i}^{(n)}]=0$.
\end{lem}
\bpr
The proof of Lemma~\ref{meaniilemma} is in subsection~\ref{meaniisection}.
\epr

\subsection{Main results}
By adding the results from Lemmas~\ref{meanijneqlemma} and~\ref{meaniilemma},
we establish the following theorem:
\begin{thm}\label{theorem51}
For $n\geq 2$,
the mean number of distinct (adjacent) pairs in a word
$ Z_{1},\ldots,Z_{n}$ is exactly
\begin{align*}
E[X^{(n)}] &=
\sum_{i=1}^{\infty}\sum_{j\neq i}\Bigg[
1
- \frac{(2P_i P_j )^{n+1}}{\sqrt{1-4P_i P_j }}\bigg(
\frac{1}{\big(1 - \sqrt{1-4P_i P_j }\!\ \big)^{n+1}}
- \frac{1}{\big(1 + \sqrt{1-4P_i P_j }\!\ \big)^{n+1}}\bigg)\Bigg]\\
&\qquad{}+\sum_{i=1}^{\infty}
\Bigg[1-\lp\frac12-\frac{1+P_i}{2\sqrt{(1-P_i )(1+3P_i )}}\rp\lp\frac{-2P_i (1-P_i )}{1-P_i +\sqrt{(1-P_i )(1+3P_i )}}\rp^{n} \cr
&\qquad\qquad\qquad{}+\lp\frac12+\frac{1+P_i}{2\sqrt{(1-P_i )(1+3P_i )}}\rp\lp\frac{-2P_i (1-P_i )}{1-P_i -\sqrt{(1-P_i )(1+3P_i )}}\rp^{n}\Bigg]
\end{align*}
For $n<2$, we have $E[X^{(n)}]=0$.
\end{thm}

In Section~\ref{section53}, we give all of the analogous parts of the
analysis for $E[(X^{(n)})^{2}]$, but we do not wrap the results into
a statement in a theorem, because the second moment has many parts,
and the notation is cumbersome.

\subsection{Analysis of the average number of distinct (adjacent) pairs}

\subsubsection{Analysis of distinct (adjacent) two letter patterns
$ij$ with $i \neq j$}\label{meanijneqsection}
If we fix $i\neq j$ and we analyze the occurrences
of the pattern $ij$, then the only ``cluster'' (to use
Bassino et~al.'s terminology) is $ij$ itself.  So the generating
function $\xi(z,t)$ of the set of clusters $C = \{ij\}$
becomes only (compare with (6) in Bassino et~al.):
$$\xi(z,t) = P_i P_j tz^{2}.$$
The generating function of the \emph{decorated texts}
(with $z$ marking the length of the words, and $t$
marking the number of decorated occurrences of $ij$, and the
coefficients are the associated probabilities)
is
$$T(z,t) = \frac{1}{1-A(z) - \xi(z,t)} =
\frac{1}{1-z-P_i P_j tz^{2}},$$
where $A(z) = z$ is the probability generating function of the alphabet
$\mathcal{A}$.

Now we use $F(z,x)$ to denote
the bivariate probability generating function of occurrences of
$ij$ (with $z$ marking the length of the words, and $x$
marking the number of occurrences of $ij$, and the coefficients are the
associated probabilities), i.e., we define
$$F(z,x) := \sum_{n=0}^{\infty}\sum_{k=0}^{\infty}
P( Z_{1},\ldots,Z_{n}\textrm{ has exactly $k$ occurrences of $ij$ as a subword})x^{k}z^{n}.
$$
We know from inclusion-exclusion
(see~\cite[Chapter~3]{FS09}
or~\cite{Bassino}) that $F(z,x) = T(z,x-1)$, so we obtain
$$F(z,x) = T(z,x-1) = \frac{1}{1-z-P_i P_j (x-1)z^{2}}.$$
The probability generating function of words with
\emph{zero occurrences} of pattern $ij$ can be obtained by considering
the case $k=0$, corresponding to the coefficients of $x^{0}$.  To
extract those coefficients, we can evaluate $F(z,x)$ at $x=0$, and we obtain
$$[x^{0}]F(z,x) = F(z,0) =
\frac{1}{1-z+P_i P_j z^{2}},$$
so, finally, the probability generating function of the words with
\emph{at least one occurrence of} $ij$ is
\begin{equation*}
\sum_{n=0}^{\infty}E[X_{i,j}^{(n)}]z^{n}
= \frac{1}{1-z} - \frac{1}{1-z+P_i P_j z^{2}}
\end{equation*}
and it follows, using Table~\ref{bigtable1}A, that
\begin{align*}
\sum_{n=0}^{\infty}E[X_{i,j}^{(n)}]z^{n}
&= \sum_{n=0}^{\infty}z^{n}
+ \frac{2P_i P_j }{(1 + \sqrt{1-4P_i P_j })
\sqrt{1-4P_i P_j }}\sum_{n=0}^{\infty}
\bigg(\frac{2P_i P_j }{1 + \sqrt{1-4P_i P_j }}\bigg)^{n}z^{n}\\
&\qquad\qquad{}- \frac{2P_i P_j }{(1 - \sqrt{1-4P_i P_j })
\sqrt{1-4P_i P_j }}\sum_{n=0}^{\infty}
\bigg(\frac{2P_i P_j }{1 - \sqrt{1-4P_i P_j }}\bigg)^{n}z^{n}
\end{align*}
and we conclude with the \emph{exact} expression for
$E[X_{i,j}^{(n)}]$ in Lemma~\ref{meanijneqlemma}.

\subsubsection{Analysis of distinct (adjacent) two letter patterns $ij$ with $i = j$}\label{meaniisection}
Now we fix $i$ and we analyze the
occurrences of the pattern $ii$.  The clusters have the
form $ii\cdots i$, i.e., they are all words that consist
of 2 or more consecutive occurrences of~$i$.
So the generating function $\xi(z,t)$ of the set of clusters $C = \{ii,iii,iiii,iiiii,\ldots\}$ becomes
$$\xi(z,t) = \frac{P_i^2 tz^{2}}{1 -
  P_i tz}.$$
The analysis is similar to the reasoning in
subsection~\ref{meanijneqsection},
and we get
\begin{equation*}
\sum_{n=0}^{\infty}E[X_{i,i}^{(n)}]z^{n}
= \frac{1}{1-z} - \frac{1}{1-z+\frac{P_i^2 z^{2}}{1 +
    P_i z}}
\end{equation*}
and then, using
Table~\ref{bigtable1}B, we have
\begin{align*}
\sum_{n=0}^{\infty}E[X_{i,i}^{(n)}]z^{n} &= \frac{1}{1-z}-
\frac{(1+P_i z)(-1+P_i +\sqrt{(1-P_i )(1+3P_i )}\!\ )}
{2\sqrt{(1-P_i )(1+3P_i )}\Big(1-\frac{-2P_i (1-P_i )}{1-P_i +\sqrt{(1-P_i )(1+3P_i )}}z\Big)}\\
&\qquad\qquad{}+\frac{(1+P_i z)(-1+P_i -\sqrt{(1-P_i )(1+3P_i )}\!\ )}{2\sqrt{(1-P_i )(1+3P_i )}\Big(1-\frac{-2P_i (1-P_i )}{1-P_i -\sqrt{(1-P_i )(1+3P_i )}}z\Big)}
\end{align*}
and we conclude with the \emph{exact} expression for
$E[X_{i,i}^{(n)}]$ in Lemma~\eqref{meaniilemma}.

\subsection{Analysis of the second moment of the number of distinct (adjacent) pairs}\label{section53}

Now we study the second moment of $X^{(n)}$, namely,
$E[(X^{(n)})^{2}]$.  We have
$$(X^{(n)})^{2} = \sum_{i=1}^{\infty}\sum_{j=1}^{\infty}X_{i,j}^{(n)}
\sum_{k=1}^{\infty}\sum_{\ell=1}^{\infty}X_{k,\ell}^{(n)}$$
so the second moment is, by linearity of expectation,
$$E[(X^{(n)})^{2}] = \sum_{i=1}^{\infty}\sum_{j=1}^{\infty}
\sum_{k=1}^{\infty}\sum_{\ell=1}^{\infty}E[X_{i,j}^{(n)}X_{k,\ell}^{(n)}].$$

We break the analysis into 4 regimes, namely:
\begin{itemize}
\item $i=j$ and $k=\ell$
\item $i=j$ and $k\neq\ell$
\item $i\neq j$ and $k=\ell$
\item $i\neq j$ and $k\neq\ell$
\end{itemize}

\subsubsection{$i=j$ and $k=\ell$}
In the case $i=j$ and $k=\ell$, we have two possibilities, namely,
either $i=j=k=\ell$ or $i=j\neq k=\ell$.
\paragraph{$i=j=k=\ell$}
In the case $i=j=k=\ell$, we have
$X_{i,j}^{(n)}X_{k,\ell}^{(n)}
= X_{i,i}^{(n)}$,
so we get
$E[X_{i,j}^{(n)}X_{k,\ell}^{(n)}]
= E[X_{i,i}^{(n)}]$,
which we already handled in Lemma~\ref{meaniilemma}.

\paragraph{$i=j\neq k=\ell$}

In the case $i=j\neq k=\ell$, we need to
analyze the occurrences of the patterns $ii$ and $kk$.
The clusters each have the
form $ii\cdots i$ or $kk\cdots k$, i.e., they are all words that consist
of 2 or more consecutive occurrences of~$i$,
or consist of 2 or more consecutive occurrences of~$k$.
So the generating function $\xi(z,t,u)$ of the set of clusters $C = \{ii,iii,iiii,iiiii,\ldots,kk,kkk,kkkk,kkkkk,\ldots\}$ becomes
$$\xi(z,t,u) = \frac{P_i^2 tz^{2}}{1 - P_i tz}
+ \frac{P_k^2 uz^{2}}{1 - P_k uz}$$
(with $z$ marking the length of the words,
and $t$ marking the number of decorated occurrences of $ii$,
and $u$ marking the number of decorated occurrences of $kk$,
and the coefficients are the associated probabilities).
The methodology now proceeds in a very similar way to the method from
Section~\ref{meanijneqsection}, but $\xi$, $T$, and $F$ all have an
additional variable, as compared to that earlier (more simple) analysis.
We have
$$
T(z,t,u)
= \frac{1}{1-A(z)-\xi(z,t,u)}
= \frac{1}{1-z-\frac{P_i^2 tz^{2}}{1 - P_i tz}
- \frac{P_k^2 uz^{2}}{1 - P_k uz}},
$$
and it follows that the probability generating function of occurrences of
$ii$ and $kk$ (with $z$ marking the length of the words,
and $x$ marking the number of occurrences of $ii$,
and $y$ marking the number of occurrences of $kk$,
and the coefficients are the associated probabilities) is
$$F(z,x,y) = T(z,x-1,y-1)
= \frac{1}{1-z-\frac{P_i^2 (x-1)z^{2}}{1 - P_i (x-1)z}
- \frac{P_k^2 (y-1)z^{2}}{1 - P_k (y-1)z}}.$$
It follows that the probability generating function of the words with
\emph{at least one occurrence of} $ii$ and
\emph{at least one occurrence of} $kk$
is
\begin{align*}
\sum_{n=0}^{\infty}E[X_{i,j}^{(n)}X_{k,\ell}^{(n)}]z^{n}
&= \sum_{n=0}^{\infty}E[X_{i,i}^{(n)}X_{k,k}^{(n)}]z^{n}\\
&= \frac{1}{1-z} - F(z,0,1)
- F(z,1,0) + F(z,0,0)\\
&= \frac{1}{1-z}
- \frac{1}{1-z+\frac{P_i^2 z^{2}}{1 + P_i z}}
- \frac{1}{1-z+\frac{P_k^2 z^{2}}{1 + P_k z}}
+ \frac{1}{1-z+\frac{P_i^2 z^{2}}{1 + P_i z}
+ \frac{P_k^2 z^{2}}{1 + P_k z}}
\end{align*}

The partial fraction decomposition for the second term is given in
Table~\ref{bigtable1}B.

The third term is the same as the second term, using $k$ instead of $i$.

The partial fraction decomposition for the
fourth term is given in
Table~\ref{bigtable2}E.

\subsubsection{$i=j$ and $k\neq\ell$}
\paragraph{$i=j$ and $k$ and $\ell$ are distinct}\label{subsubsection221}

The clusters each have the
form $ii\cdots i$ or $k\ell$, i.e., they are all words that consist
of either 2 or more consecutive occurrences of~$i$,
or simply the word $k\ell$.
So $\xi(z,t,u)$ of the set of clusters $C = \{ii,iii,iiii,iiiii,\ldots,k\ell\}$ becomes
$$\xi(z,t,u) = \frac{P_i^2 tz^{2}}{1 - P_i tz}
+ P_k P_\ell uz^{2}$$
(with $z$ marking the length of the words,
and $t$ marking the number of decorated occurrences of $ij$,
and $u$ marking the number of decorated occurrences of $k\ell$,
and the coefficients are the associated probabilities).  It follows that
\begin{align*}
T(z,t,u)
&= \frac{1}{1-z-\frac{P_i^2 tz^{2}}{1 - P_i tz}
- P_k P_\ell uz^{2}},
\end{align*}
and
$$F(z,x,y) = T(z,x-1,y-1)
= \frac{1}{1-z-\frac{P_i^2 (x-1)z^{2}}{1 - P_i (x-1)z}
- P_k P_\ell (y-1)z^{2}}.$$
It follows that the probability generating function of the words with
\emph{at least one occurrence of} $ij$ and
\emph{at least one occurrence of} $k\ell$
is
\begin{align*}
\sum_{n=0}^{\infty}E[X_{i,j}^{(n)}X_{k,\ell}^{(n)}]z^{n}
&= \frac{1}{1-z} - F(z,0,1)
- F(z,1,0) + F(z,0,0)\\
&= \frac{1}{1-z}- \frac{1}{1-z+\frac{P_i^2 z^{2}}{1 + P_i z}}
- \frac{1}{1-z+P_k P_\ell z^{2}}
+ \frac{1}{1-z+\frac{P_i^2 z^{2}}{1 + P_i z}
+ P_k P_\ell z^{2}}
\end{align*}

The partial fraction decomposition for the second term is given in
Table~\ref{bigtable1}B.

The partial fraction decomposition for the third term is given in
Table~\ref{bigtable1}A,
using~$k$ and~$\ell$ instead of~$i$ and~$j$.

The partial fraction decomposition for the fourth term is given in
Table~\ref{bigtable2}F.

\paragraph{$i=j=k\neq\ell$}\label{subsubsection222}

The clusters each have the
form $ii\cdots i$ or $ii\cdots i\ell$, i.e., they are all words that consist
of 2 or more consecutive occurrences of~$i$,
or of 1 or more consecutive occurrences of~$i$ followed by $\ell$.
So $\xi(z,t,u)$ of the set of clusters $C = \{ii,iii,iiii,iiiii,\ldots,i\ell,ii\ell,iii\ell,iiii\ell,\ldots\}$ becomes
$$\xi(z,t,u)
= \frac{P_i^2 tz^{2}}{1 - P_i tz}
+ \frac{P_i P_\ell uz^{2}}{1 - P_i tz}
= \frac{P_i z(P_i tz + P_\ell uz)}{1 - P_i tz},
$$
and
$$F(z,x,y) = \frac{1}{1-z-\frac{P_i z(P_i (x-1)z + P_\ell (y-1)z)}{1 - P_i (x-1)z}}.$$
It follows that the probability generating function of the words with
\emph{at least one occurrence of} $ij$ and
\emph{at least one occurrence of} $k\ell$
is
\begin{align*}
\sum_{n=0}^{\infty}E[X_{i,j}^{(n)}X_{k,\ell}^{(n)}]z^{n}
&= \frac{1}{1-z} - F(z,0,1)
- F(z,1,0) + F(z,0,0)\\
&= \frac{1}{1-z}- \frac{1}{1-z+\frac{P_i^2 z^{2}}{1 + P_i z}}
- \frac{1}{1-z+P_i P_\ell z^{2}}
+ \frac{1}{1-z+\frac{P_i z(P_i z + P_\ell z)}{1 + P_i z}}
\end{align*}

The partial fraction decomposition for the second term is given in
Table~\ref{bigtable1}B.

The partial fraction decomposition for the third term is given in
Table~\ref{bigtable1}A,
using $\ell$ instead of~$j$.

The partial fraction decomposition for the fourth term is given in
Table~\ref{bigtable1}D.

\paragraph{$i=j=\ell\neq k$}\label{subsubsection223}

The cluster have the
form $ki\cdots i$ or $ii\cdots i$, i.e., they are all words that consist
of $k$ followed by 1 or more consecutive occurrences of~$i$,
or of 2 or more consecutive occurrences of~$i$.
So $\xi(z,t,u)$ of the set of clusters $C = \{ki,kii,kiii,kiiii,\ldots,ii,iii,iiii,iiiii,\ldots\}$ becomes
$$\xi(z,t,u)
= \frac{P_i P_k uz^{2}}{1 - P_i tz}
+ \frac{P_i^2 tz^{2}}{1 - P_i tz}
= \frac{P_i z(P_k uz + P_i tz)}{1 - P_i tz},
$$
and
$$F(z,x,y)
= \frac{1}{1-z-\frac{P_i z(P_k (y-1)z + P_i (x-1)z)}{1 - P_i (x-1)z}}.$$
It follows that the probability generating function of the words with
\emph{at least one occurrence of} $ij$ and
\emph{at least one occurrence of} $k\ell$
is
\begin{align*}
\sum_{n=0}^{\infty}E[X_{i,j}^{(n)}X_{k,\ell}^{(n)}]z^{n}
&= \frac{1}{1-z} - F(z,0,1)
- F(z,1,0) + F(z,0,0)\\
&= \frac{1}{1-z}- \frac{1}{1-z+\frac{P_i^2 z^{2}}{1 + P_i z}}
- \frac{1}{1-z+P_i P_k z^{2}}
+ \frac{1}{1-z+\frac{P_i z(P_k z + P_i z)}{1 + P_i z}}
\end{align*}

The partial fraction decomposition for the second term is given in
Table~\ref{bigtable1}B.

The partial fraction decomposition for the third term is given in
Table~\ref{bigtable1}A,
using $k$ instead of~$j$.

The partial fraction decomposition for the fourth term is given in
Table~\ref{bigtable1}D,
using $k$ instead of~$\ell$.

\subsubsection{$i\neq j$ and $k=\ell$}

\paragraph{$k=\ell$ and $i$ and $j$ are distinct}

Same as section~\ref{subsubsection221} but with $i$ and $k$ exchanged,
and with $j$ and $\ell$ exchanged.

\paragraph{$k=\ell=i\neq j$}

Same as section~\ref{subsubsection222} but with $i$ and $k$ exchanged,
and with $j$ and $\ell$ exchanged.

\paragraph{$k=\ell=j\neq i$}

Same as section~\ref{subsubsection223} but with $i$ and $k$ exchanged,
and with $j$ and $\ell$ exchanged.

\subsubsection{$i\neq j$ and $k\neq\ell$}
\paragraph{$i$ and $j$ and $k$ and $\ell$ are distinct}\label{subsubsection241}

The clusters are $ij$ and $k\ell$.
So $\xi(z,t,u)$ of the set of clusters $C = \{ij,k\ell\}$ becomes
$$\xi(z,t,u)
= P_i P_j tz^{2} + P_k P_\ell uz^{2},
$$
and
$$F(z,x,y)
= \frac{1}{1-z-P_i P_j (x-1)z^{2} - P_k P_\ell (y-1)z^{2}}.$$
It follows that the probability generating function of the words with
\emph{at least one occurrence of} $ij$ and
\emph{at least one occurrence of} $k\ell$
is
\begin{align*}
\sum_{n=0}^{\infty}E[X_{i,j}^{(n)}X_{k,\ell}^{(n)}]z^{n}
&= \frac{1}{1-z} - F(z,0,1)
- F(z,1,0) + F(z,0,0)\\
&= \frac{1}{1-z}- \frac{1}{1-z+P_i P_j z^{2}}
                - \frac{1}{1-z+P_k P_\ell z^{2}}\\
&\qquad{} + \frac{1}{1-z+P_i P_j z^{2} + P_k P_\ell z^{2}}
\end{align*}

The partial fraction decomposition for the second term is given in
Table~\ref{bigtable1}A.

The partial fraction decomposition for the third term is given in
Table~\ref{bigtable1}A,
using~$k$ and~$\ell$ instead of~$i$ and~$j$.

The partial fraction decomposition for the fourth term is given in
Table~\ref{bigtable1}C.

\paragraph{$k=i$ and $j$ and $\ell$ are distinct}\label{subsubsection242}

The clusters
are $ij$ and $i\ell$.  So, by the same analysis from
section~\ref{subsubsection241}, we get
\begin{align*}
\sum_{n=0}^{\infty}E[X_{i,j}^{(n)}X_{i,\ell}^{(n)}]z^{n}
&= \frac{1}{1-z} - F(z,0,1)
- F(z,1,0) + F(z,0,0)\\
&= \frac{1}{1-z}- \frac{1}{1-z+P_i P_j z^{2}}
                - \frac{1}{1-z+P_i P_\ell z^{2}}\\
&\qquad{} + \frac{1}{1-z+P_i P_j z^{2} + P_i P_\ell z^{2}}
\end{align*}

Exactly as in section~\ref{subsubsection241} above:

The partial fraction decomposition for the second term is given in
Table~\ref{bigtable1}A.

The partial fraction decomposition for the third term is given in
Table~\ref{bigtable1}A,
using~$\ell$ instead of~$j$.

The partial fraction decomposition for the fourth term is given in
Table~\ref{bigtable1}C,
using~$i$ instead of~$k$.

\paragraph{$k=j$ and $i$ and $\ell$ are distinct}\label{subsubsection243}

The clusters are $ij$, $ij\ell$ and $j\ell$.
So $\xi(z,t,u)$ of the set of clusters $C = \{ij,ij\ell,j\ell\}$ becomes
$$\xi(z,t,u)
=  P_i P_j tz^{2}+ P_j P_\ell uz^{2} + P_i P_j P_\ell tuz^{3},
$$
and
$$F(z,x,y)
= \frac{1}{1-z - P_iP_j(x-1)z^{2}-P_j P_\ell (y-1)z^{2}-P_i P_j P_\ell (x-1)(y-1)z^{3}}.$$
It follows that the probability generating function of the words with
\emph{at least one occurrence of} $ij$ and
\emph{at least one occurrence of} $j\ell$
is
\begin{align*}
\sum_{n=0}^{\infty}E[X_{i,j}^{(n)}X_{j,\ell}^{(n)}]z^{n}
&= \frac{1}{1-z} - F(z,0,1)
- F(z,1,0) + F(z,0,0)\\
&= \frac{1}{1-z} - \frac{1}{1\!-\!z\!+\! P_i P_j z^{2}}
- \frac{1}{1\!-\!z\!+\! P_j P_\ell z^{2}}
+ \frac{1}{1-z + P_i P_j z^{2}+ P_j P_\ell z^{2}-P_i P_j P_\ell z^{3}}
\end{align*}

The partial fraction decompositions for the second and third terms are given in Table~\ref{bigtable1}A, once
using~$j$ and~$\ell$ instead of~$i$ and~$j$.

The partial fraction decomposition for the fourth term is given in
Table~\ref{bigtable2}G.

\paragraph{$i=\ell$ and $k$ and $j$ are distinct}\label{subsubsection244}

Same as section~\ref{subsubsection243} but with $i$ and $k$ exchanged,
and with $j$ and $\ell$ exchanged.

\paragraph{$\ell=j$ and $i$ and $k$ are distinct}\label{subsubsection245}

The clusters are $ij$ and $kj$.  So, by the same analysis from
section~\ref{subsubsection241}, we get
\begin{align*}
\sum_{n=0}^{\infty}E[X_{i,j}^{(n)}X_{k,j}^{(n)}]z^{n}
&= \frac{1}{1-z} - F(z,0,1)
- F(z,1,0) + F(z,0,0)\\
&= \frac{1}{1-z}- \frac{1}{1-z+P_i P_j z^{2}}
                - \frac{1}{1-z+P_k P_j z^{2}}\\
&\qquad{} + \frac{1}{1-z+P_i P_j z^{2} + P_k P_j z^{2}}
\end{align*}

Exactly as in section~\ref{subsubsection241} above:

The partial fraction decomposition for the second term is given in
Table~\ref{bigtable1}A.

The partial fraction decomposition for the third term is given in
Table~\ref{bigtable1}A,
using~$k$ instead of~$i$.

The partial fraction decomposition for the fourth term is given in
Table~\ref{bigtable1}C,
using~$j$ instead of~$\ell$.

\paragraph{$i=k$ and $j=\ell$ are distinct}\label{subsubsection246}

In this case we have
$X_{i,j}^{(n)}X_{k,\ell}^{(n)}
= X_{i,j}^{(n)}$,
so we get
$E[X_{i,j}^{(n)}X_{k,\ell}^{(n)}]
= E[X_{i,j}^{(n)}]$,
which we already handled in Lemma~\ref{meanijneqlemma}.

\paragraph{$i=\ell$ and $j=k$ are distinct}\label{subsubsection247}

The clusters each have the
form $ijiji\ldots$ or $jijij\ldots$.
So $\xi(z,t,u)$ of the set of clusters $C = \{ij,iji,ijij,ijiji\ldots,ji,jij,jiji,jijij\ldots\}$ becomes
$$\xi(z,t,u)
= \frac{P_i P_j tz^{2}}        {1 - P_i P_j tuz^{2}}
+ \frac{P_i^2 P_j tuz^{3}}{1 - P_i P_j tuz^{2}}
+ \frac{P_j P_i uz^{2}}        {1 - P_i P_j tuz^{2}}
+ \frac{P_i P_j^2 tuz^{3}}{1 - P_i P_j tuz^{2}},
$$
and $F(z,x,y) = 1/(1-z-\xi(z,x-1,y-1))$.
It follows that the probability generating function of the words with
\emph{at least one occurrence of} $ij$ and
\emph{at least one occurrence of} $k\ell$
is
\begin{align*}
\sum_{n=0}^{\infty}E[X_{i,j}^{(n)}X_{k,\ell}^{(n)}]z^{n}
&= \frac{1}{1-z} - F(z,0,1)
- F(z,1,0) + F(z,0,0)\\
&= \frac{1}{1-z}
- \frac{1}{1-z+P_i P_j z^{2} }- \frac{1}{1-z + P_j P_i z^{2}}\\
&\qquad{}+ \frac{1}{1-z+\frac{P_i P_j z^{2}}{1 - P_i P_j z^{2}}
- \frac{P_i^2 P_j z^{3}}{1 - P_i P_j z^{2}}
+ \frac{P_j P_i z^{2}}        {1 - P_i P_j z^{2}}
- \frac{P_i P_j^2 z^{3}}{1 - P_iP_jz^{2}}}
\end{align*}

The partial fraction decomposition for the second and for the third
term is given in
Table~\ref{bigtable1}A.

The partial fraction decomposition for the fourth term is given in
Table~\ref{bigtable2}H.

As mentioned immediately after Theorem~\ref{theorem51}, we do not wrap
all of the analysis from Section~\ref{section53} into a theorem
(because it would be very lengthy), but we 
have precisely analyzed every aspect that is needed for exactly
characterizing the second moment $E[(X^{(n)})^{2}]$.

\acknowledgements
\label{sec:ack}

We would like to thank  B. Pittel for providing the Poisson distribution of $X_{i,j}(m)$, and  B. Salvy for suggesting the use of
the Maple package gfun.

Furthermore we would like to thank an anonymous referee, whose suggestions led to substantial improvements of the paper. In particular we want to thank for pointing out to us
a connection to Stirling numbers, that is stated in Remark~\ref{rem:OEIS}.

M.D.~Ward's research is supported by National Science Foundation (NSF)
grants 0939370, 1246818, 2005632, 2123321, 2118329, and 2235473, by the
Foundation for Food and Agriculture Research (FFAR) grant 534662, by
the National 
Institute of Food and Agriculture (NIFA) grants 2019-67032-29077,
2020-70003-32299, 2021-38420-34943, and 2022-67021-37022,
by the Society Of Actuaries grant 19111857, by
Cummins Inc., by Gro Master, by Lilly Endowment, and by Sandia National Laboratories.

\nocite{*}
\bibliographystyle{abbrvnat}
\bibliography{louchard-schachinger-ward-V4}
\label{sec:biblio}

\addcontentsline{toc}{section}{Appendices}
\begin{appendices}
\section{Some Mellin transforms}
%
%
To keep the paper self contained we give here a short outline on how to use Mellin transforms to obtain asymptotic expansions. The reader seeking more detail is referred to \cite{FGD95} for a nice exposition.
Subsections A.1, A.2, and A.3 are devoted to asymptotic equivalents of three sums that play a crucial role in our paper.

The Mellin transform $f^*(s)$ of $f(x)$, also denoted
$\BMC\lb f(x);s \rb$, is given by
\[f^*(s)=\int_0^\II f(x) x^{s-1} dx.\]
The interior of the set of $s$ for which the integral converges is an open strip $\langle a,b\rangle:=\{s\in\mathbb{C}:a<\Re s<b\}$, called the \emph{fundamental strip}, with $a, b$ depending on how $f$ behaves at $0$ and $\infty$. For example, we have $\BMC\lb e^{-x};s \rb=\Gamma(s)$,
with fundamental strip $\langle 0,\infty\rangle$, and
$\BMC\lb 1-e^{-x};s \rb=-\Gamma(s)$,
with fundamental strip $\langle -1,0\rangle$.
When computing the Mellin transform of so called \emph{harmonic sums}, the rescaling rule turns out to be very useful:
\[\BMC\lb \sum_k\la_k f(\mu_k x);s \rb =\sum_k\frac{\la_k}{\mu_k^s}\cdot f^*(s).\]
In the case that $f^*(s)$ can be meromorphically continued to
a strip $\langle a,\bar b\rangle$ with $\bar b>b$, information on the poles of $f^*(s)$ leads to asymptotic properties of $f(n), n \ra \II$.
This is called  the fundamental correspondence. In particular,
if there is a pole \[\frac{1}{(s-\xi)^{k+1}}\] of $f^*(s)$ at $\xi=\sig+ \ii t$ to the right of the fundamental strip, then this pole will contribute the term
\[-\frac{(-1)^k}{k!}  \ln(n)^k n^{-\sig} e^{-\ii t\ln(n)},\]
which is precisely the residue of $\frac{-n^{-s}}{(s-\xi)^{k+1}}$ at $s=\xi$,
to an asymptotic expansion of $f(n)$ at $\II$. Justification comes from residue calculus: If $f$ is smooth enough,
the inverse transform applies to yield
$f(n)=\frac1{2\pi\ii}\int_{c-\ii\infty}^{c+\ii\infty}n^{-s}f^*(s)ds$ with $a<c<b$. If for $\beta<\bar b$ the set of poles $z$ in $\langle a,\bar b\rangle$ satisfying $\Re z<\beta$ is denoted $\mathcal{S}_\beta$,
and there is no pole with real part $\beta$,
we have
$$f(n)=\frac1{2\pi\ii}\int_{\beta-\ii\infty}^{\beta+\ii\infty}n^{-s}f^*(s)ds-\sum_{z\in\mathcal{S}_\beta}\textup{Res}(f^*(s)n^{-s})_{s=z}
=-\sum_{z\in\mathcal{S}_\beta}\textup{Res}(f^*(s)n^{-s})_{s=z}+\BO(n^{-\beta})
,$$
with the integral being $\BO(n^{-\beta})$, provided that $f^*(s)$ decreases fast enough for $s=\beta+\ii t$ and $|t|\to \infty$. Equality of left and right hand side is established
by using a sequence of contours $\rho_k$ being the boundaries
of rectangles $\{z\in\mathbb{C}:c\le\Re z\le\beta,-h_k\le\Im z\le h_k\}$ with
$h_k\to\infty$, verifying that $f^*(s)$ decreases fast enough on the horizontal segments of $\rho_k$, as $k\to \infty$, and applying residue calculus.

Here is an illustration of fast enough decrease. $\Gam(s)$ decreases exponentially in the direction  $\ii\II$:
\[|\Gam( \sig+\ii t )| \sim \sqrt{2\pi}|t|^{\sig-1/2}e^{-\pi |t|/2}.\]
Also, similarly fast decrease can be observed for all other transforms we encounter.

In the following, recall the notations
\[L:=\ln\tfrac1q\mbox{ and }\chi:=\tfrac{2\ii\pi}{L}.\]

\subsection   {}                       \label{A1}
Let \[G(n):=\sum_{i\ge0}\lp 1- e^{-n q^{2i}}\rp .\]
The Mellin transform of this sum is
$G^*(s)=-\dfrac{1}{1-q^{-2s}}\Gam(s)$, with fundamental strip
$\langle -1,0 \rangle$, to the right of which the meromorphic extension of $G^*(s)$ has poles at $s=0$ and $s=\frac{\ell\chi}2$ for $\ell\in\Z\setminus\{0\}$, with singular expansions
\[G^*(s)\asymp\frac{1}{2Ls^2}-\lb\frac{\gam}{2L}+\frac12\rb\frac1s,\quad
\mbox{ and }G^*(s)\asymp\frac1{2L}\frac{\Gam\big(\frac{\ell\chi}2\big)}{\big(s-\frac{\ell\chi}2\big)},\mbox{ for }\ell\in\Z\setminus\{0\}.\]
Noting that there are no other singularities to the right of $\langle
-1,0 \rangle$, the error term in the following expansion can be chosen $\BO(n^{-\beta})$ with any fixed $\beta>0$.
\[G(n)\sim \frac1{2L} \ln(n)+\lb\frac{\gam}{2L}+\frac12\rb-\frac{1}{2L}\sum_{\ell \in \Z\setminus \{0\}}  \Gam\Big( \frac{\ell\chi}2\Big)             n^{-\ell\chi/2}.\]
\ignore{\subsection   {}                       \label{A11}
Let
\[G:=\sum_i i\lp  1- e^{-n q^{i}}\rp. \]
This  has the MT, with FS  $\langle -1,0 \rangle$ and  $\chi:=\frac{2\ii\pi}{L}$
\[F(s)=-\frac{Q^s}{(1-Q^s)^2}\Gam(s)\asymp -\frac{1}{L^2 s^3}+\frac{\gam}{L^2 s^2}
+\bigg(-\frac{\pi^2/12+\gam^2/2}{L^2}+\frac{1}{12}\bigg)\frac1s-\frac{\Gam(\ell\chi)}{L^2}\frac{1}{(s-\ell\chi)^2}
-\frac{\Gam'(\ell\chi)}{L^2(s-\ell\chi)},\]
\[G\sim \frac{\ln(n)^2}{2L^2}+\frac{\gam}{L^2}\ln(n)+\lb \frac{\pi^2/12+\gam^2/2}{L^2}-\frac{1}{12}\rb
+\frac{1}{L^2}\sum_{\ell \in \Z\setminus \{0\}}\lb -\Gam\lp \frac{2\ii \ell\pi}{L}\rp\ln(n)
+\Gam'\lp \frac{2\ii \ell\pi}{L}\rp \rb  e^{-2 \ii \ell\pi\ln(n)/L}.\]
\textbf{}
}%
\subsection   {}                       \label{A2}
Let
\[\tilde G(n):=\sum_{i,j\ge0} \lp  1- e^{-n q^{i+j}}\rp=\sum_{k\ge0} (k+1)\lp  1- e^{-n q^{k}}\rp. \]
Here we have $\tilde G^*(s)=-\dfrac{1}{(1-q^{-s})^2}\Gam(s)$, with fundamental strip
$\langle -1,0 \rangle$, and poles at $s=0$ and $s=\ell\chi$ for $\ell\in\Z\setminus\{0\}$, with singular expansions
\[\tilde G^*(s)\asymp -\frac{1}{L^2 s^3}+\lb\frac{\gam}{L^2}+\frac1L\rb\frac{1}{s^2}
-\lb\frac{\pi^2+6\gam^2}{12L^2}+\frac{5}{12}+\frac\gam L\rb\frac1s,\]
and \[\tilde G^*(s)\asymp-\frac{\Gam(\ell\chi)}{L^2(s-\ell\chi)^2}
-\frac{\Gam'(\ell\chi)-L\Gam(\ell\chi)}{L^2(s-\ell\chi)},\mbox{ for }\ell\in\Z\setminus\{0\},\]
leading to
\[\tilde G(n)\sim \frac{\ln(n)^2}{2L^2}+\lb\frac{\gam}{L^2}+\frac1L\rb\ln(n)+\lb \frac{\pi^2+6\gam^2}{12L^2}+\frac{5}{12}+\frac\gam L\rb
+\frac{1}{L^2}\!\!\!\sum_{\ell \in \Z\setminus \{0\}}\!\!\!\!\!\!\lb
\Gam'\lp \ell\chi\rp-(\ln(n)+L)\Gam\lp \ell\chi\rp \rb  n^{-\ell\chi},\]
again with error term $\BO(n^{-\beta})$ with any fixed $\beta>0$.
\subsection   {}                       \label{A3}
Set
\[\hat G(n):=\sum_{i,j,k}(e^{nP_iP_jP_k}-1)e^{-nP_i P_j-nP_j P_k}.\]
This leads to the Mellin transform, with fundamental strip $\langle -1,0 \rangle$,
\bals
\hat G^*(s)=&\sum_{i,j,k} \int_0^\II (e^{xP_iP_jP_k}-1) e^{-xP_iP_j-xP_jP_k} x^{s-1}dx\\
=&\sum_{j}P_j^{-s}\sum_{i,k} \int_0^\II \lb e^{-y(P_i+P_k-P_iP_k)}-e^{-y(P_i+P_k)} \rb  y^{s-1}dy\\
=&\frac{q^s}{p^s(q^s-1)}\Gam(s)\sum_{i,k}\lb(P_i+P_k-P_iP_k)^{-s}-(P_i+P_k)^{-s} \rb\\
=&\lp\frac{q}{p}\rp^{2s}\frac{\Gam(s)}{q^s-1}F_1(s),
\end{align*}
where
\bals
F_1(s)&=\sum_{i,k}\lb(q^i+q^k-pq^{i+k-1})^{-s}-(q^i+q^k)^{-s} \rb\\
&=\sum_{i\ge1}q^{-is}\left[(2-P_i  )^{-s}-2^{-s}+2
\sum_{j\ge1}\lb(1+q^j-pq^{i+j-1})^{-s}-(1+q^j)^{-s} \rb\right].
\end{align*}
Note that $F_1(s)$, being a general Dirichlet series in the variable $-s$, is analytic at least for $\sigma=\Re s<1$, since, using the Mean Value Theorem, we have
$$\left|(1+q^j-pq^{i+j-1})^{-\sigma}-(1+q^j)^{-\sigma}\right|
\le|\sigma|\frac{pq^{i+j-1}}{(1+q^j)^{1+\sigma}},
$$
and therefore
$$|F_1(s)|\le2|\sigma|\sum_{i\ge1}q^{i(1-\sigma)}
\sum_{j\ge0}\frac{P_j }{(1+q^j)^{1+\sigma}}<\infty.
$$
Moreover, $F_1(0)=0$,
so to the right of the fundamental strip we have the singular expansions
\[\hat G^*(s)\asymp -\frac{F'_1(0)}{Ls},\quad\textup{and}\
\hat G^*(s)\asymp-p^{-2\ell\chi}\frac{\Gam(\ell\chi)}{L}\frac{ F_1(\ell\chi)}{s-\ell\chi},\ \textup{for}\ \ell\in\mathbb{Z}\setminus\{0\}
.\]
This leads to
$$\hat G(n)=\frac{F'_1(0)}{L}+\frac1L\sum_{\ell \in \Z\setminus \{0\}}
  \Gam(\ell\chi)F_1(\ell\chi)(np^2)^{-\ell\chi}+\BO(n^{-\beta}),$$
with any fixed $\beta<1$, where the constant term simplifies to
$$\frac{F'_1(0)}{L}=\frac1L\ln\Bigg(\prod_{i,k\ge1}\frac{q^i+q^k}{q^i+q^k-pq^{i+k-1}}\Bigg).$$
\end{appendices}
\ignore{
\section   {Proofs of equations \eqref{iialpha} and \eqref{iralpha}}
\label{B}
The estimate \eqref{iialpha} follows from the following general result: If for some $c<1$ a set $\mathcal{P}=\{x_i:i\in\mathbb{N}\}$ satisfies $x_i>0$ and $\frac{x_{i+1}}{x_i}\le c$ for $i\in\mathbb{N}$, then $\sum_{x\in\mathcal{P}}x^\alpha e^{-x}<\infty$. For a proof observe that there is a constant $C_\alpha>0$ such that $x^\alpha e^{-x}\le\min(x^\alpha,C_\alpha x^{-\alpha})$ for $x>0$. Let $\bar x:=(C_\alpha)^{1/(2\alpha)}$. Then
$$\sum_{x\in\mathcal{P}}x^\alpha e^{-x}\le \sum_{x\in\mathcal{P}\cap\,]0,\bar x]}x^\alpha+\sum_{x\in\mathcal{P}\cap[\bar x,\infty[}C_\alpha x^{-\alpha}\le \bar x^\alpha\sum_{i\ge0}c^i+C_\alpha\bar x^{-\alpha}\sum_{i\ge0}c^i=2\frac{\sqrt{C_\alpha}}{1-c}.$$
The estimate \eqref{iralpha} can be deduced from \eqref{iialpha}, using
$\beta=1$, observing
$$\sum_{i,k\ge1}(nP_iP_k)^\alpha
e^{-nP_iP_k}=\sum_{\ell\ge2}(\ell-1)(n\tfrac pqP_\ell)^\alpha
e^{-n\tfrac pqP_\ell},$$
and furthermore
$$\sum_{\ell\ge2}\ell(nP_\ell)^\alpha
e^{-nP_\ell}=\sum_{\ell\ge2}\BO(\ln n-\ln(nP_\ell))(nP_\ell)^\alpha e^{-nP_\ell}=\BO(\ln n),$$
because of $x^\alpha\ln x=\BO(x^{\alpha/2})$.
}%
\end{document}